\documentclass[a4paper,12pt]{article}

\usepackage[utf8]{inputenc}
\usepackage[T1,OT2]{fontenc}
\usepackage[left=2.5cm,right=2.5cm,top=1.5cm,bottom=2cm]{geometry}
\usepackage{amsmath,amssymb,makeidx}
\usepackage[french]{babel}

\usepackage[language=french,backend=biber]{biblatex}
\usepackage{fourier}
\usepackage{graphics}
\usepackage{graphicx}
\usepackage{amsthm}
\usepackage{faktor}
\usepackage{thmtools}
\usepackage[hyperindex=true]{hyperref}
\usepackage{color}
\usepackage{enumitem}
\newlist{itemizeth}{itemize}{2}
\setlist[itemizeth]{label=\textbullet,noitemsep,topsep=0 mm}
\newlist{enumerateth}{enumerate}{2}
\setlist[enumerateth]{label*=\textbf{\alph*)},itemsep=0.2mm}
\setlist[itemize]{fullwidth,topsep=0 mm, leftmargin=0pt}
\setlist[enumerate]{label*=\textbf{\arabic*)},itemsep=0.2mm, leftmargin=0pt}

\declaretheoremstyle[
title=Démonstration,
numbered=no,
headfont=\normalfont\bfseries,
notefont=\bfseries, notebraces={}{},
postheadspace=17pt,
bodyfont=\normalfont,
headindent=0pt,
qed=\qedsymbol,
spacebelow=3,
]{demostyle}

\declaretheoremstyle[
headfont=\normalfont\bfseries,
notefont=\mdseries, notebraces={(}{)},
bodyfont=\normalfont,
thmbox=S
]{standard}

\declaretheoremstyle[
title=Théorème,
headfont=\normalfont\scshape,
notefont=\mdseries, notebraces={(}{)},
thmbox=S,
]{theo}

\declaretheoremstyle[
title=Théorème,
headfont=\scshape\bfseries,
notefont=\mdseries, notebraces={(}{)},
thmbox=L,
]{import}

\declaretheoremstyle[headfont=\normalfont\itshape, 
numbered=no
]{rem}

\declaretheorem[title=Définition,style=standard]{defi}
\declaretheorem[title=Proposition,style=standard]{prop}
\declaretheorem[title=Lemme,style=standard]{lem}
\declaretheorem[title=Corollaire,style=standard]{coro}
\declaretheorem[style=theo]{Th}

\declaretheorem[title=Remarque,style=rem]{rqu}
\declaretheorem[title=Remarques,style=rem]{rqus}
\declaretheorem[style=demostyle]{dem}

\makeatletter
\newcommand*{\house}[1]{%
   \mathord{%
     \mathpalette\@house{#1}%
   }%
}
\newcommand*{\@house}[2]{%
   \dimen@=\fontdimen8 %
       \ifx#1\scriptscriptstyle\scriptscriptfont
       \else\ifx#1\scriptstyle\scriptfont
       \else\textfont\fi\fi
       3 %
   \sbox0{%
     $#1%
       \vrule width\dimen@\relax
       \overline{%
         \kern2\dimen@
         \begingroup 
           #2%
         \endgroup
         \kern2\dimen@
       }%
       \vrule width\dimen@\relax
       \mathsurround=1.5\dimen@ 
     $%
   }%
   \ht0=\dimexpr\ht0-\dimen@\relax
   \dp0=\dimexpr\dp0+2\dimen@\relax
   \vbox{%
     \kern\dimen@ 
     \copy0 %
   }%
}

\newcommand{\R}{\mathbb{R}}
\newcommand{\N}{\mathbb{N}}

\newcommand{\Z}{\mathbb{Z}}
\newcommand{\Q}{\mathbb{Q}}
\newcommand{\C}{\mathbb{C}}
\newcommand{\K}{\mathbb{K}}
\newcommand{\Oal}{\mathcal{O}}
\newcommand{\Qbar}{\overline{\mathbb{Q}}}
\newcommand{\Spec}{\mathrm{Spec} \,}
\newcommand{\Gal}{\mathrm{Gal} \,}
\newcommand{\GL}{\mathrm{GL}}
\newcommand{\Pro}{\mathbb{P}}
\newcommand{\Sym}{\mathrm{Sym}}
\newcommand{\ord}{\mathrm{ord} \,}
\newcommand{\ssi}{\Leftrightarrow}
\newcommand{\gf}{\dfrac}

\newcommand{\Img}{\mathrm{Im} \,}
\newcommand{\den}{\mathrm{den}}
\newcommand{\ddz}{\gf{\mathrm{d}}{\mathrm{d}z}}

\newcommand{\fonction}[5]{\begin{array}[t]{lrcl} 
#1: & #2 & \longrightarrow & #3 \\
    & #4 & \longmapsto & #5 \end{array}}
\newcommand{\Zemla}{{\fontencoding{OT2}\selectfont\char'132}}

\title{Le théorème d'André-Chudnovsky-Katz \og{} au sens large \fg{}}
\date{\today}
\author{Gabriel Lepetit}
\makeindex
\begin{document}

\maketitle
\begin{abstract}
Les $E$- et $G$-fonctions de Siegel ont été définies en deux sens, strict et large, conjecturalement équivalents. En reprenant et complétant une esquisse d'André \cite{AndregevreyII}, nous énonçons et démontrons l'analogue \emph{au sens large} du théorème d'André-Chudnovsky-Katz, qui est un théorème de structure sur les $G$-opérateurs \emph{au sens strict} (il s'agit d'opérateurs différentiels annulant les $G$-fonctions \emph{au sens strict}). Nous en déduisons un théorème de structure sur les $E$-opérateurs \emph{au sens large}, qui sont des opérateurs différentiels annulant les $E$-fonctions \emph{au sens large}. En application de ce dernier théorème,  nous donnons une nouvelle preuve d'une généralisation par André \cite{Andre2014} du théorème de Siegel-Shidlovskii sur l'indépendance algébrique des valeurs des $E$-fonctions \emph{au sens large}.
\end{abstract}

\section{Introduction}

Le but de cet article est d'étudier la structure des $G$-opérateurs et des $E$-opérateurs \emph{au sens large}. Nous commençons par donner quelques éléments de contexte. Les notions de $E$- et de $G$-fonctions ont été introduites par Siegel dans \cite{Siegelarticle} pour généraliser le théorème de Lindemann-Weierstrass sur l'indépendance algébrique des valeurs de la fonction exponentielle.

\begin{defi}\label{def:gfonctionlarge}
Une \emph{$G$-fonction} au sens large est une série $f(z)=\sum\limits_{n=0}^{\infty} a_n z^n \in \Qbar\llbracket z\rrbracket$ telle que
\begin{enumerateth}
\item $f$ est solution d'une équation différentielle linéaire à coefficients dans $\Qbar(z)$ ;
\item Pour tout $\varepsilon >0$, il existe $n_1(\varepsilon) \in \N$ tel que $\forall n \geqslant n_1(\varepsilon), \house{a_n} \leqslant (n!)^{\varepsilon}$, où $\house{a_n}$ est la \emph{maison} de $a_n$, c'est-à-dire le maximum des modules des conjugués (au sens de Galois) de $a_n$   ;
\item Pour tout $\varepsilon >0$, il existe $n_2(\varepsilon) \in \N$ tel que $\forall n \geqslant n_2(\varepsilon), \mathrm{den}(a_0, \dots,a_n) \leqslant (n!)^{\varepsilon}$, où $\mathrm{den}(a_0, \dots, a_n)$ est le plus petit $d \in \N^*$ tel que $da_0, \dots, da_n$ sont des entiers algébriques.
\end{enumerateth}
\end{defi}

On dispose d'une notion de $G$-fonction \emph{au sens strict}, plus restrictive que la définition \ref{def:gfonctionlarge}, qui est en fait celle considérée par Siegel.

\begin{defi} \label{def:gfonctionstrict}
Une \emph{$G$-fonction} au sens strict est une série $f(z)=\sum\limits_{n=0}^{\infty} a_n z^n \in \Qbar\llbracket z\rrbracket$ telle que
\begin{enumerateth}
\item $f$ est solution d'une équation différentielle linéaire à coefficients dans $\Qbar(z)$ ;
\item Il existe $C_1 >0$ tel que $\forall n \in \N, \house{a_n} \leqslant C_1^{n+1}$ ;
\item Il existe $C_2 >0$ tel que $\forall n \in \N,  \mathrm{den}(a_0, \dots,a_n) \leqslant C_2^{n+1}$.
\end{enumerateth}
\end{defi}

On définit de la même manière les $E$-fonctions \emph{au sens strict} (resp. \emph{au sens large}) qui sont les séries $f(z)=\sum\limits_{n=0}^{\infty} \gf{a_n}{n!} z^n \in \Qbar\llbracket z\rrbracket$ vérifiant la condition \textbf{a)} de la définition \ref{def:gfonctionstrict} (resp. \ref{def:gfonctionlarge}) et telles que les $a_n$ vérifient les conditions \textbf{b)} et \textbf{c)} de la définition \ref{def:gfonctionstrict} (resp. \ref{def:gfonctionlarge}). Siegel a étudié les $E$-fonctions \emph{au sens large}.

L'étude des $E$- et $G$-fonctions a été développée, entre autres, par Shidlovskii \cite{Shidlovskii}, puis poursuivie par Nesterenko et Shidlovskii \cite{NesterenkoShidlovskii}, André \cite{AndregevreyI,AndregevreyII}, Bombieri, Galochkin \cite{Galochkin74} et Beukers \cite{Beukers2006}. 

Siegel a étudié les $E$-fonctions \emph{au sens large}, mais n'a fait qu'évoquer les $G$-fonctions \emph{au sens large}. 
Il est conjecturé que les définitions large et stricte sont équivalentes pour les $E$- et $G$-fonctions, mais cela n'a pas été prouvé à ce jour. Précisément, on sait que la condition \textbf{b)} de la définition \ref{def:gfonctionlarge} implique, sous la condition \textbf{a)}, la condition \textbf{b)} de la définition \ref{def:gfonctionstrict}, car on peut appliquer des estimation \og{} Gevrey \fg{} dues à Perron (cf \cite{Perron}, voir aussi \cite[pp. 85--86]{Ramis}) ; en revanche, on ne sait pas si la condition \textbf{c)} de la définition \ref{def:gfonctionlarge} implique la condition \textbf{c)} de la définition \ref{def:gfonctionstrict} (cf \cite[p. 715]{AndregevreyI}).

\medskip

Un théorème fondamental de D. et G. Chudnovsky (voir \cite{Chudnovsky}) affirme que l'équation différentielle minimale satisfaite par une $G$-fonction \emph{au sens strict} vérifie une condition de croissance modérée appelée \emph{condition de Galochkin}. Ceci implique, entre autres, qu'elle est fuchsienne. Par ailleurs, la condition de Galochkin est équivalente à une condition introduite par Bombieri, qui implique par un théorème de Katz \cite[p. 98]{Dwork} que l'équation différentielle minimale en question est à exposants rationnels en tout point de $\mathbb{P}^1(\C)$ (théorème d'André-Chudnovsky-Katz).

Dans \cite[pp. 746--747]{AndregevreyII}, André esquisse la preuve du fait que \textit{la singularité en l’infini d’un opérateur $\phi$ d’ordre minimal annulant
une $G$-fonction au sens large est régulière.} Son argument est le suivant :

\og{}
\textit{Pour établir le point ci-dessus, le critère $p$-adique de régularité de Katz montre qu’il suffit d’établir que $\sum_{p(v) \leqslant n} \log R_v(\phi, 1) = o(\log n)$. Avec les notations de \cite{Andre}, on voit facilement que cette condition découle d’une estimation $$\sum_{v \; \text{finie}} h_{v,n}(\phi) = \sum_{p(v) \leqslant n} h_{v,n}(\phi)=o(\log n) ;$$
or les estimations de \cite[p. 122]{Andre} donnent
$$\sum_{v \; \text{finie}} h_{v,n}(\phi) \leqslant C_1 \gf{1} {n} \sum_{v} \log \max(1, |a_0|_v, \dots, |a_{n C_2}|_v )+ C_3$$
(pour des constantes $C_i $ indépendantes de $n$), et la condition $\mathbf{(G^{-})}$ équivaut à $$\gf{1}{n} \sum_{v \; \text{finie}} \log \max(1, |a_0|_v,\dots,|a_n|_v) = o(\log n). \fg{}$$}
La condition $\mathbf{(G^{-})}$ correspond aux points \textbf{b)} et \textbf{c)} de la définition \ref{def:gfonctionlarge} (cf \cite[p. 714]{AndregevreyI}).

Dans cet article, nous commençons par donner les détails de l'esquisse d'André et nous compléterons ensuite ses résultats. Précisément, nous allons donner la démonstration du théorème suivant. Il est implicite dans l'esquisse ci-dessus d'André, même s'il ne l'énonce pas formellement.
\begin{Th} \label{th:acklarge}
Soit $f(z) \in \Qbar\llbracket z\rrbracket $ une $G$-fonction \emph{au sens large}, et $L \in \Qbar(z)\left[\gf{\mathrm{d}}{\mathrm{d}z}\right]$ un opérateur différentiel non nul d'ordre minimal $\mu$ pour $f$ tel que $L(f(z))=0$. Alors 
\begin{itemizeth}
\item L'opérateur $L$ est un $G$-opérateur \emph{au sens large}.
 \item L'opérateur  $L$ est globalement nilpotent.
 \item Tout point de $\mathbb{P}^{1}(\C)$ est un point singulier régulier de $L$ et les exposants de $L$ en tout point sont dans $\Q$.
\end{itemizeth}
\end{Th}

La preuve de ce résultat fera l'objet des sections \ref{sec:chudnovsky} et \ref{sec:acklarge}.
Dans la partie \ref{sec:baseACK}, nous
raffinerons le théorème \ref{th:acklarge} en précisant la forme d'une base de solutions d'un $G$-opérateur \emph{au sens large}. C'est l'objet du théorème suivant, qui constitue un analogue complet du théorème d'André-Chudnovsky-Katz.

\begin{Th}
 Soit $f(z)$ une $G$-fonction \emph{au sens large} et $L \in \Qbar(z)[\mathrm{d}/\mathrm{d}z] \setminus \{ 0 \}$ un opérateur différentiel tel que $L(f(z))=0$ et d'ordre minimal $\mu$ pour $f$. 

Alors au voisinage de tout $\alpha \in \mathbb{P}^1(\Qbar)$, il existe une base de solutions de $Ly(z)=0$ de la forme $$(f_1(z-\alpha), \dots, f_{\mu}(z-\alpha)) (z-\alpha)^{C_{\alpha}},$$ où  $C_{\alpha} \in \mathcal{M}_{\mu}(\Qbar)$ est triangulaire supérieure à valeurs propres dans $\Q$  et les $f_i(u) \in \Qbar\llbracket u \rrbracket$ sont des $G$-fonctions \emph{au sens large}.
\end{Th}

L'étude de la structure des $G$-opérateurs permet d'obtenir des informations sur les équations différentielles satisfaites par les $E$-fonctions. En effet, via la transformée de Fourier-Laplace des opérateurs différentiels, André (cf \cite{AndregevreyI, AndregevreyII}) en a déduit que toute $E$-fonction \emph{au sens strict} était annulée par un \emph{$E$-opérateur} dont les seules singularités sont $0$ et $\infty$, la première étant régulière. André a déduit du théorème d'André-Chudnovsky-Katz un théorème de structure sur les $E$-opérateurs. La section \ref{sec:csqACK} sera consacrée à la définition et à l'étude des $E$-opérateurs \emph{au sens large}, à l'aide du théorème \ref{th:acklarge}, et à ses conséquences diophantiennes sur les valeurs des $E$-fonctions \emph{au sens large}.

\medskip

\section{Un analogue \og{} large \fg{} du théorème des Chudnovsky} \label{sec:chudnovsky}

Soit $\mathbf{f}={}^t (f_1(z), \dots, f_{\mu}(z))\in \Qbar\llbracket z\rrbracket^{\mu}$ vérifiant $\mathbf{f}'=G\mathbf{f}$, avec $G \in \mathcal{M}_{\mu}\left(\Qbar(z) \right)$. Soit $G_s \in \mathcal{M}_{\mu}\left(\Qbar(z)\right)$ la matrice telle que $\mathbf{f}^{(s)}=G_s \mathbf{f}$. On montre par récurrence que les $G_s$, $s \in \N$, sont liées par la relation \begin{equation}\label{eq:recurrenceGs}
G_{s+1}=G_s G+G'_s,
\end{equation} où $G'_s$ désigne la matrice $G_s$ dérivée coefficient par coefficient. On prend $T(z) \in \Qbar[z]$ le plus petit dénominateur commun de tous les coefficients de la matrice $G(z)$. On montre également par récurrence sur $s$ que \begin{equation} \label{eq:tsgspoly}\forall s \in \N, \;\; T^s G_s \in \mathcal{M}_{\mu}\left(\Qbar[z]\right).\end{equation}

\subsection{Condition de Galochkin \emph{au sens large}}\label{subsec:galochkinlarge}

Rappelons tout d'abord la définition de la condition de Galochkin \emph{au sens strict}, introduite dans \cite{Galochkin74}.

\begin{defi}[Galochkin] \label{def:galochkinstrict}
On note, pour $s \in \N$, $q_s$ le plus petit dénominateur supérieur ou égal à $1$ de tous les coefficients des coefficients des matrices $T(z)^m \gf{G_m(z)}{m !}$, quand $m \in \{ 1, \dots, s \}$. On dit que le système $y'=Gy$ \emph{vérifie la condition de Galochkin au sens strict} si $$\exists C>0 : \; \forall s \in \N, \quad  q_s \leqslant C^{s+1}.$$
\end{defi}

On a alors le théorème fondamental suivant (cf  \cite[p. 17]{Chudnovsky}).

\begin{Th}[Chudnovsky] \label{th:chudnovsky}
Sous les mêmes hypothèses que ci-dessus, si pour tout $i \in \{ 1, \dots, \mu \}$, $f_i(z)$ est une $G$-fonction \emph{au sens strict} et $(f_1(z), \dots, f_{\mu}(z))$ est une famille libre sur $\Qbar(z)$, alors $G$ vérifie la condition de Galochkin \emph{au sens strict}.
\end{Th}

Dans \cite[p. 747]{AndregevreyII}, André a introduit, en la formulant différemment, la condition suivante, qui est adaptée au contexte des $G$-fonctions \emph{au sens large}.

\begin{defi}[André] \label{def:galochkinlarge}
Avec les notations de la définition précédente, on dit que le système $y'=Gy$ \emph{vérifie la condition de Galochkin \emph{au sens large}} si $$\forall \varepsilon >0, \; \exists s_0(\varepsilon) \in \N : \; \forall s \geqslant s_0(\varepsilon), \; q_s \leqslant (s!)^{\varepsilon}.$$
\end{defi}

Rappelons que si $L=\left(\gf{\mathrm{d}}{\mathrm{d}z}\right)^{\mu}+a_{1}(z) \left(\gf{\mathrm{d}}{\mathrm{d}z}\right)^{\mu-1}+\dots+a_n(z) \not\equiv 0$ est un opérateur différentiel d'ordre $\mu$ à coefficients dans $\Qbar(z)$, la matrice compagnon de $L$ est $$A_L=\begin{pmatrix}
0 & 1 &  & (0) \\ 
 & \ddots & \ddots &  \\ 
(0) &  & 0 & 1 \\ 
-a_{\mu} & \dots &  & -a_1
\end{pmatrix}.$$ On sait que les solutions du système différentiel $y'=Gy$ sont les vecteurs $\mathbf{f}={}^t (f, f', \dots, f^{(\mu-1)})$ quand $L(f(z))=0$.

Suivant la définition des $G$-opérateurs au sens strict (cf \cite[p. 718]{AndregevreyI}), on peut considérer une notion analogue \emph{au sens large}.

\begin{defi}
Soit $L \in \Qbar(z)\left[\mathrm{d}/\mathrm{d}z\right]$. On dit que $L$ est un \emph{$G$-opérateur \emph{au sens large}} (resp. \emph{au sens strict}) si la matrice compagnon de $L$ vérifie la condition de Galochkin \emph{au sens large} (resp. \emph{au sens strict}).
\end{defi}

La dénomination de \og{}$G$-opérateur \fg{} est justifiée par la proposition suivante. 

\begin{prop} \label{prop:baseGfonctionspointordinaire}
Soit $L \in \Qbar(z)\left[\mathrm{d}/\mathrm{d}z\right]$ un $G$-opérateur \emph{au sens large} non nul d'ordre $\mu$. Soit $\alpha \in \Qbar$ un point ordinaire de $L$, alors il existe une base de solutions de l'équation $L(y(z))=0$ au voisinage de $\alpha$ de la forme $(f_1(z-\alpha), \dots, f_{\mu}(z-\alpha))$, où les $f_i(u)$ sont des $G$-fonctions \emph{au sens large}.
\end{prop}

Il est bien connu que cette proposition est également vraie pour les $G$-opérateurs \emph{au sens strict}, la preuve ci-dessous s'adaptant \emph{mutatis mutandis}. 

Le théorème d'André-Chudnovsky-Katz affirme entre autres qu'un $G$-opérateur \emph{au sens strict} a au voisinage de toute singularité $\alpha$ une base de solutions de la forme $(f_1(z-\alpha), \dots, f_{\mu}(z-\alpha)) (z-\alpha)^{C_{\alpha}}$, où $C_{\alpha} \in \mathcal{M}_{\mu}(\Qbar)$ a ses valeurs propres rationnelles et les $f_i(u)$ sont des $G$-fonctions \emph{au sens strict}. On verra dans la section \ref{sec:baseACK} que ceci est également vrai \emph{au sens large}.

\begin{dem}[de la proposition \ref{prop:baseGfonctionspointordinaire}]
Notons $G(z)=A_L(z)$ la matrice compagnon de $L$. Comme $\alpha$ est un point ordinaire, on sait qu'il existe une base de solutions $(f_1(z-\alpha), \dots, f_{\mu}(z-\alpha))$ de l'équation $L(y(z))=0$, où les $f_i$ sont holomorphes au voisinage de $0$. On sait aussi que la matrice wronskienne de cette base $Y(z) \in \mathcal{M}_{\mu}(\Qbar\llbracket z\rrbracket)$ a un rayon de convergence non nul et est telle que $Y(\alpha) \in \GL_{\mu}(\Qbar)$ et $Y'(z)=G(z)Y(z)$, de sorte que $Y^{(s)}(z)=G_s(z) Y(z)$ pour tout entier $s$. D'où \begin{equation}\label{eq:YetGs} Y(z)=\sum\limits_{n=0}^{\infty} Y^{(n)}(\alpha)(z-\alpha)^n=\left(\sum\limits_{n=0}^{\infty}  \gf{G_n(\alpha)}{n!}(z-\alpha)^n\right)Y(\alpha).\end{equation} Puisque $\alpha$ est un point ordinaire, $G(z)$ n'a pas de pôle en $\alpha$ et la condition de Galochkin \emph{au sens large} implique qu'il existe une suite d'entiers positifs $(q_n)_{n \in \N}$ telle que $$\forall n \in \N,\; \forall k \leqslant n, \;\; q_n \gf{G_k(\alpha)}{k!} \in \mathcal{M}_{\mu}\left(\Oal_{\Qbar}\right) \quad \mathrm{et} \quad \forall \varepsilon >0,\; \exists n_0(\epsilon) \in \N, \;\; \forall n \geqslant n_0(\varepsilon), q_n \leqslant (n!)^{\varepsilon}.$$ Ainsi, selon \eqref{eq:YetGs}, on a $\forall n \in \N, \mathrm{den}(Y(\alpha)) q_n Y^{(n)}(\alpha) \in \mathcal{M}_{\mu}\left(\Oal_{\Qbar}\right)$. Ainsi, les $f_i(z)$ vérifient la condition \textbf{c)} de la définition \ref{def:gfonctionlarge}. 

Par ailleurs, soit $\K$ un corps de nombres galoisien contenant $\alpha$ tel que $L \in \K(z)\left[\mathrm{d}/\mathrm{d}z\right]$. Soit $\tau \in \Gal(\K/\Q)$. Si $L=\sum\limits_{k=0}^{\mu} a_k(z) \left(\gf{\mathrm{d}}{\mathrm{d}z}\right)^k$, $a_k(z) \in \K[z]$, on définit $L^{\tau} :=\sum\limits_{k=0}^{\mu} a_k^{\tau}(z) \left(\gf{\mathrm{d}}{\mathrm{d}z}\right)^k$ en étendant l'action de $\tau$ à $\Qbar\llbracket z\rrbracket$ coefficient par coefficient. 

\medskip

Alors pour tout $i \in \{1, \dots, \mu\}$, $L^{\tau}(f_i^{\tau}(z-\sigma(\alpha)))=0$. De plus, comme $a_{\mu}(\alpha) \neq 0$, on a $a_{\mu}^{\tau}(\tau(\alpha))=\tau(a_{\mu}(\alpha)) \neq 0$, de sorte que $\tau(\alpha)$ est un point ordinaire de $L^{\tau}$. Ainsi, $f_i^{\tau}$ est analytique au voisinage de $0$. Ceci valant pour tout $\tau$, on en déduit que, si $f_i(z)=\sum\limits_{n=0}^{\infty} b_{i,n} z^n$, il existe une constante $C>0$ telle que $\forall n \in \N, \house{b_{i,n}} \leqslant C^{n+1}$, ce qui prouve que les $f_i(z)$ vérifient la condition \textbf{b)} de la définition \ref{def:gfonctionlarge}.\end{dem}

L'ensemble des $G$-opérateurs \emph{au sens large} possède une structure algébrique analogue à celle de l'ensemble des $G$-opérateurs \emph{au sens strict}.  On a les propriétés suivantes (listées par André dans le cas strict \cite[p. 720]{AndregevreyI}) :
\begin{itemize}[label=\textbullet]
\item Un produit de $G$-opérateurs \emph{au sens large} est un $G$-opérateur \emph{au sens large}.
\item Tout diviseur à droite d'un $G$-opérateur \emph{au sens large} dans $\Qbar(z)[\mathrm{d}/\mathrm{d}z]$ est un $G$-opérateur \emph{au sens large}.
\item L'opérateur adjoint $L^*$ d'un $G$-opérateur \emph{au sens large} $L$ est un $G$-opérateur \emph{au sens large}.
\item Si $L$ et $L'$ sont deux $G$-opérateurs \emph{au sens large}, alors ils admettent un multiple commun à gauche qui est un $G$-opérateur \emph{au sens large} (\emph{propriété de Ore à gauche}).
\end{itemize}

La démonstration de ces propriétés est donnée dans la section \ref{sec:ch1senslarge}. Elle consiste en l'adaptation de propriétés des modules différentiels données dans le cas strict par André \cite[\S IV]{Andre}, qui sont détaillées dans \cite{LepetitSize}. 

\bigskip

Le but de la suite de cette partie est de démontrer un analogue \emph{au sens large} du théorème~\ref{th:chudnovsky}.

\begin{Th} \label{th:chudnovskylarge}
Le théorème \ref{th:chudnovsky} reste vrai si l'on remplace \og{} strict \fg{} par \og{} large \fg{}.
\end{Th}

Ceci implique en particulier que si $f$ est une $G$-fonction \emph{au sens large}, tout opérateur différentiel non nul $L$ à coefficients dans $\Qbar(z)$ tel que $L(f(z))=0$ et d'ordre minimal pour $f$ est un $G$-opérateur \emph{au sens large}. En effet, la condition de minimalité sur l'ordre $\mu$ de $L$ impose que $(f, \dots, f^{(\mu-1)})$ est libre sur $\Qbar(z)$. Le théorème \ref{th:chudnovskylarge} assure donc que la condition de Galochkin \emph{au sens large} est vérifiée pour $A_L$.

Notons que la proposition \ref{prop:baseGfonctionspointordinaire} constitue une réciproque partielle du théorème \ref{th:chudnovskylarge}.

\subsection{Démonstration du théorème \ref{th:chudnovskylarge}} \label{subsec:preuvechudnovsky}

La preuve que nous allons présenter est une adaptation de la preuve originale du théorème \ref{th:chudnovsky} \cite[pp. 38--50]{Chudnovsky} au cas des $G$-fonctions \emph{au sens large}. Les six premières étapes de la démonstration sont identiques dans les cas strict et large puisque les conditions \textbf{b)} et \textbf{c)} de la définition d'une $G$-fonction (définitions \ref{def:gfonctionlarge} ou \ref{def:gfonctionstrict}) n'y sont pas utilisées, y compris dans le lemme de Shidlovskii évoqué dans l'étape 3. Dans \cite[pp. 21--22]{Beukers}, Beukers a reformulé les idées des Chudnovsky en une trame condensée ; à des fins de clarté, nous suivrons cette trame dans les étapes 1 à 6, en la détaillant. La nouveauté de cette démonstration est l'étape 7, dans laquelle les conditions \textbf{b)} et \textbf{c)} de la définition \ref{def:gfonctionlarge}, spécifiques aux $G$-fonctions \emph{au sens large}, seront utilisées. André a également mentionné dans \cite[pp. 746--747] {AndregevreyII} un argument qui permet d'adapter \emph{au sens large} la preuve du théorème des Chudnovsky donnée dans \cite[pp. 112--123]{Andre}, mais sans donner de détails.

\medskip

Remarquons tout d'abord que si $\K$ est un corps de nombres contenant les coefficients des coefficients de $G(z)$ et les $f_i(0)$, $1 \leqslant i \leqslant \mu$, alors $G \in \mathcal{M}_{\mu}(\K(z))$ et $\mathbf{f} \in \K\llbracket z\rrbracket^\mu$. En effet, cela découle de l'équation $\mathbf{f}'=G \mathbf{f}$ en écrivant $G$ comme un élément de $\mathcal{M}_{\mu}\big(\K((z))\big)$ et identifiant les coefficients du développement en série de Laurent de part et d'autre.

\bigskip
\bigskip

\textit{Notations et hypothèses} :

On a $T(z) \in \K[z]$, mais quitte à multiplier par un entier adapté, on peut supposer que $T(z) \in \Oal_\K[z]$ et $T(z) G(z) \in \mathcal{M}_{\mu}(\Oal_\K [z])$.

On note $D=\gf{\mathrm{d}}{\mathrm{d}z}$ la dérivation usuelle sur $\K\llbracket z\rrbracket$.

Si $A \in \K\llbracket z\rrbracket$ et $\ell \in \N$, on note $A = O\left(z^{\ell} \right)$ s'il existe $B \in \K\llbracket z\rrbracket$ tel que $A=z^{\ell} B$.

On note $\delta=[\K:\Q]$ le degré du corps de nombres $\K$.

\medskip

\textbf{Étape 1} : Soient $N, M \in \N$. On introduit des approximants de Padé $(Q, \mathbf{P})$ de type II de paramètres $(N,M)$ associés à $\mathbf{f}$ dont on laisse les paramètres libres pour l'instant, c'est-à-dire des polynômes $Q, \mathbf{P}_1, \dots, \mathbf{P}_{\mu} \in \K[z]$ tels que $\deg(Q) \leqslant N$, $\max\limits_{1 \leqslant i \leqslant \mu} \deg(P_i) \leqslant N$ et si $\mathbf{P}=(P_1, \dots, P_n)$, alors $$ Q \mathbf{f} - \mathbf{P}=O\left(z^{N+M}\right).$$ On ne discutera des conditions d'existence de tels approximants de Padé que dans l'étape~7.

On a pour tout $m < N+M$, $\gf{T^m}{m!} (D-G)^m \mathbf{P} \in \K[z]$, ce qui est immédiat en utilisant la formule de Leibniz. De plus, on va montrer par récurrence sur $m$ que \begin{equation}\label{eq:derivpade}
\forall m \in \N^*, \quad \gf{T^m}{m!} Q^{(m)} \mathbf{f} - \gf{T^m}{m!} (D-G)^m \mathbf{P}=O\left(z^{N+M-m}\right).
\end{equation}

Pour $m=1$, c'est la définition. 

Supposons la relation vraie au rang $m$. Alors en dérivant \eqref{eq:derivpade} et en multipliant par $T$, on~a 
\begin{multline*}
(mT' T^m Q^{(m)}+T^{m+1} Q^{(m+1}) \mathbf{f}+T^{m+1} Q^{(m)} \mathbf{f}' \\ - mT' T^m (D-G)^m \mathbf{P} - T^{m+1} D(D-G)^m \mathbf{P} =O\left(z^{N+M-(m+1)}\right).
\end{multline*}
Or, $\mathbf{f}'=G\mathbf{f}$ donc en multipliant \eqref{eq:derivpade} par la matrice polynomiale $TG$ et en retranchant à l'équation précédente, on obtient
\begin{multline*}
(mT' T^m Q^{(m)}+T^{m+1} Q^{(m+1}) \mathbf{f}- mT' T^m (D-G)^m \mathbf{P}\\ - T^{m+1} (D-G)^{m+1} \mathbf{P}=O\left(z^{N+M-(m+1)}\right). 
\end{multline*}
Finalement, comme $T^m Q^{(m)}\mathbf{f}- T^m (D-G)^m \mathbf{P}=O\left(z^{N+M-m}\right)$, on a le résultat voulu.

\medskip

\textbf{Étape 2} : Montrons que
\begin{equation} \label{eq:similileibniz}
\forall s \in \N^*, \quad \gf{G_s}{s!} \mathbf{P}=\sum_{j=0}^s \gf{(-1)^j}{(s-j)! j!} D^{s-j}(D-G)^j \mathbf{P}.
\end{equation}
On procède par récurrence :
\begin{itemize}[label=\textbullet]
\item pour $s=1$, $G_s=G$ et $G \mathbf{P}=D(\mathbf{P})-(D-G)\mathbf{P}$.

\item Soit $s \in \N^*$, supposons la formule \eqref{eq:similileibniz} vraie pour $s$. 

Alors $G_{s+1}=G_s G+G'_s$, donc $\gf{G_{s+1}}{(s+1)!}=\gf{1}{s+1} \left( \gf{G_s}{s!} G+\gf{G'_s}{s!} \right)$. Donc en appliquant l'hypothèse de récurrence au vecteur $G \mathbf{P}$, on a $$ \gf{G_{s+1}}{(s+1)!} \mathbf{P}=\gf{1}{s+1} \left(\sum_{j=0}^s \gf{(-1)^j}{(s-j)! j!} D^{s-j}(D-G)^j (G\mathbf{P})+\gf{G'_s}{s!} \mathbf{P} \right).$$ Or, $$\gf{G'_s}{s!} \mathbf{P}=\left(\gf{G_s}{s!} \mathbf{P}\right)' - \gf{G_s}{s!} \mathbf{P}'=\sum_{j=0}^s \gf{(-1)^j}{(s-j)! j!} D^{s+1-j}(D-G)^j \mathbf{P}-\sum_{j=0}^s \gf{(-1)^j}{(s-j)! j!} D^{s-j}(D-G)^j D \mathbf{P}.$$ Donc, en remarquant que $\forall j \in \{ 0, \dots, s \}, (D-G)^j D=(D-G)^{j+1}+(D-G)^j G $, on a 
\begin{align*}
\gf{G_{s+1}}{(s+1)!} \mathbf{P} &= \gf{1}{s+1} \left(\sum_{j=0}^s \gf{(-1)^j}{(s-j)! j!} D^{s+1-j} (D-G)^j \mathbf{P} - \sum_{j=0}^s \gf{(-1)^j}{(s-j)! j!} D^{s-j}(D-G)^{j+1} \mathbf{P} \right) \\
&= \gf{1}{s+1} \left(\sum_{j=0}^s \gf{(-1)^j}{(s-j)! j!} D^{s+1-j} (D-G)^j \mathbf{P} + \sum_{k=1}^{s+1} \gf{(-1)^k}{(s-k+1)! (k-1)!} D^{s+1-k}(D-G)^{k} \mathbf{P} \right) \\
&= \gf{1}{s+1} \left(\sum_{j=1}^s \gf{(-1)^j (s+1-j)+j}{(s+1-j)! j!} D^{s+1-j} (D-G)^j \mathbf{P} + \right. \\ 
& \qquad \left. \sum_{k=1}^{s+1} \gf{(-1)^k}{(s-k+1)! (k-1)!} D^{s+1-k}(D-G)^{k} \mathbf{P} + \gf{D^{s+1}}{s!} \mathbf{P}+\gf{(-1)^{s+1}}{s!}(D-G)^{s+1} \mathbf{P} \right) \\
&= \sum\limits_{j=0}^{s+1} \gf{(-1)^j}{(s+1-j)! j!} D^{s+1-j} (D-G)^j \mathbf{P}. 
\end{align*}
Cela conclut la récurrence.
\end{itemize}

\bigskip

\textbf{Étape 3} : Utilisation du lemme de Shidlovskii.

On note pour $h \in \N$, $\mathbf{P}_h=\gf{1}{h!} (D-G)^h \mathbf{P}$ et $R_{(h)} \in \mathcal{M}_{\mu}(\K(z))$ la matrice dont la $j$\up{ème} colonne est $\binom{h+j-1}{j-1} \mathbf{P}_{h+j-1}$. Alors la formule \eqref{eq:similileibniz} implique immédiatement que

$$ \gf{G_s}{s!} R_{(0)}=\sum\limits_{j=0}^s \gf{(-1)^j}{(s-j)!} D^{s-j} R_{(j)}.$$ Selon le lemme de Shidlovskii pour les approximants de Padé de type II (voir \cite[p. 115]{Andre}), la matrice $R_{(0)}$ est inversible pourvu que $M$ soit assez grand ce qui sera réalisé quand on spécifiera $M$ et $N$ dans l'étape 7. Donc \begin{equation} \label{eq:conclusionetape3}
T^s \gf{G_s}{s!}=\sum\limits_{j=0}^s (-1)^j T^{s+\mu-1} \gf{D^{s-j} R_{(j)}}{(s-j)!} (T^{\mu-1} R_{(0)})^{-1}.
\end{equation}

\medskip

\textbf{Étape 4} :
Soit $d$ le dénominateur commun des coefficients d'ordres inférieurs à $N+M$ du développement en série entière de $\mathbf{f}$. On suppose trouvés des approximants de Padé $Q, \mathbf{P}_1, \dots, \mathbf{P}_{\mu} \in \K[X]$ de $d\mathbf{f}$, c'est-à-dire des polynômes tels que $\deg(Q) \leqslant N$, $\max\limits_{1 \leqslant i \leqslant \mu} \deg(P_i) \leqslant N$ et $$ Q(d\mathbf{f}) - \mathbf{P}=O\left(z^{N+M}\right),$$ avec $\mathbf{P}=(P_1, \dots, P_{\mu})$. On fait l'hypothèse supplémentaire que $Q$ est un polynôme à coefficients entiers algébriques. On peut alors appliquer les résultats des trois étapes précédentes, dont on conservera les notations, à $Q$ et $\mathbf{P}$ puisque $d\mathbf{f}$ est encore solution du système différentiel $y'=Gy$.

Selon \eqref{eq:derivpade}, on a
\begin{equation} \label{eq:derivpadedenom}
 \forall m \leqslant N+M, \quad  \gf{T^m}{m!} Q^{(m)} (d\mathbf{f}) - T^m \mathbf{P}_m=O\left(z^{N+M-m}\right).
\end{equation}
On remarque que $\gf{Q^{(m)}}{m!} \in \Oal_\K[z]$ puisque $Q$ est à coefficients dans $\Oal_\K$. On en déduit que si $N+M-m > \max\limits_{1 \leqslant i \leqslant \mu} \deg(\mathbf{P}_{i,m})$, les coefficients de $T^m \mathbf{P}_m$ sont des éléments de $\Oal_\K[z]$.

\medskip

Montrons par récurrence sur $m$ que \begin{equation} \label{eq:bornedegrePm}
\max\limits_{1 \leqslant i \leqslant \mu} \deg(T^m \mathbf{P}_{i,m}) \leqslant N+tm.
\end{equation}

\begin{itemize}[label=\textbullet]
\item Pour $m=0$, il s'agit simplement du fait que les composantes de $\mathbf{P}$ sont de degrés inférieurs à $N$.

\item Pour $m=1$, $T \mathbf{P}'-T G \mathbf{P}$ a des composantes de degrés inférieurs à $N+t$, car les coefficients de $T(z)G$ sont de degrés bornés par $t$.

\item Soit $m \in \N$, supposons le résultat vrai au rang $m$. Alors\begin{align*}
(D-G)(T^m (D-G)^m \mathbf{P}) &= m T' T^{m-1} (D-G)^m \mathbf{P}+T^m D(D-G)^m \mathbf{P}-T^m G(D-G)^m \mathbf{P} \\
&= m T' T^{m-1} (D-G)^m \mathbf{P}+T^m (D-G)^{m+1} \mathbf{P}.
\end{align*} Donc $$ T^{m+1} (D-G)^{m+1} \mathbf{P}=(D-G)(T^m (D-G)^m \mathbf{P})- m T' T^m (D-G)^m \mathbf{P}. $$
Or, en utilisant à la fois l'hypothèse de récurrence et le cas $m=1$, on voit que $T (D-G)(T^m (D-G)^m \mathbf{P})$ a ses composantes de degrés bornés par $N+mt+t=N+(m+1)t$ ; par ailleurs, l'hypothèse de récurrence nous assure que $m T' T^m (D-G)^m \mathbf{P}$ a des composantes de degrés inférieurs à $t+N+mt=N+(m+1)t$. On en déduit le résultat souhaité~\eqref{eq:bornedegrePm}.
\end{itemize}

On déduit, à condition que $N+M-m \geqslant N+tm$, c'est-à-dire 
\begin{equation} \label{eq:conditionentieral}
m \leqslant \gf{M}{t+1},
\end{equation}
que $T^m \mathbf{P}_m \in \Oal_\K[z]^\mu$. Ceci implique immédiatement que si $j$ est un entier naturel tel que $j+\mu-1 \leqslant \gf{M}{t+1}$, alors $T^{j+\mu-1} R_{(j)}$ est une matrice à coefficients dans $\Oal_\K [z]$. En particulier, $T^{\mu-1} R_{(0)} \in \mathcal{M}_{\mu}(\Oal_\K[z])$.

\medskip

Montrons à présent que \begin{equation} \label{eq:derivrjentier}
\forall s \in \N, \;\; s+\mu-1 \leqslant \gf{M}{t+1}, \;\; \forall j \in \{ 0, \dots, s \}, \quad \gf{T^{s+\mu-1}}{(s-j)!} D^{s-j} R_{(j)} \in \mathcal{M}_{\mu}(\Oal_\K [z]).
\end{equation}
Rappelons la formule de Leibniz généralisée : si $\ell \in \N^*$ et $f_1, \dots, f_{\ell}$ sont $l$ fonctions dérivables $k$ fois, alors $$ (f_1 \dots f_{\ell})^{(k)}=\sum_{i_1+\dots+i_{\ell}=k} \binom{k}{i_1, \dots, i_{\ell}} \prod_{1 \leqslant t \leqslant k} f_t^{(i_t)}.$$ En particulier, considérons un élément quelconque $W \in \Oal_\K[z]$. Alors $$\gf{(W^\ell)^{(k)}}{k!}=\sum\limits_{i_1+\dots+i_{\ell}=k} \prod\limits_{1 \leqslant t \leqslant \ell} \gf{W^{(i_t)}}{i_t!}.$$ Si $k <\ell$, pour tout $(i_1, \dots, i_{\ell})$ intervenant dans la somme, chaque terme $\prod\limits_{1 \leqslant t \leqslant \ell} \gf{W^{(i_t)}}{i_t!}$ contient au moins $\ell-k$ indices de dérivation nuls. Ainsi, comme pour tout entier $s$,$\gf{W^{(s)}}{s!} \in \Oal_\K[z]$, on a $\gf{(W^\ell)^{(k)}}{k!} \in W^{\ell-k} \Oal_\K[z]$.

\medskip

Déduisons de cela par récurrence le résultat \eqref{eq:derivrjentier}. Pour $s=0$, c'est évident.

Soit $s \in \N^*$, supposons \eqref{eq:derivrjentier} vrai pour $s' \in \{ 0, \dots, s-1 \}$. Par la formule de Leibniz, on a, pour $j \in \{0, \dots, s \}$, \begin{align*}
\gf{D^{s-j}\big(T^{s+\mu-1} R_{(j)}\big)}{(s-j)!} &= \sum_{k=0}^{s-j} \binom{s-j}{k} \times \gf{1}{(s-j)!} \big(T^{s+\mu-1}\big)^{(k)} D^{s-j-k} R_{(j)} \\
&= T^{s+\mu-1} \gf{D^{s-j} R_{(j)}}{(s-j)!}+ \sum_{k=1}^{s-j} \gf{\big(T^{s+\mu-1}\big)^{(k)}}{k!} \gf{D^{s-j-k} R_{(j)}}{(s-j-k)!} \\
&= T^{s+\mu-1} \gf{D^{s-j} R_{(j)}}{(s-j)!}+ \sum_{k=1}^{s-j} U_k T^{s+\mu-1-k} \gf{D^{s-k-j} R_{(j)}}{(s-k-j)!},
\end{align*} avec $U_k \in \Oal_\K[z]$, en utilisant la remarque précédente.

Or, d'une part, $$\gf{D^{s-j}\big(T^{s+\mu-1} R_{(j)}\big)}{(s-j)!} \in \mathcal{M}_{\mu}(\Oal_\K[z]),$$ puisque $T^{s+\mu-1}=T^{s-j} T^{j+\mu-1} R_{(j)} \in \mathcal{M}_{\mu}(\Oal_\K[z])$, et d'autre part, par hypothèse de récurrence, $$\forall k \in \{ 1, \dots, s-j \}, \quad T^{s-k+\mu-1} \gf{D^{s-k-j} R_{(j)}}{(s-k-j)!} \in \mathcal{M}_{\mu}(\Oal_\K[z]).$$ Par conséquent, $$T^{s+\mu-1} \gf{D^{s-j} R_{(j)}}{(s-j)!} \in \mathcal{M}_{\mu}(\Oal_\K[z]),$$ ce qu'il fallait démontrer.

\bigskip

\textbf{Étape 5} : Lien entre $q_s$ et taille des coefficients de $\det(T^{\mu-1} R_{(0)})$.

Selon \eqref{eq:conclusionetape3}, on a pour tout $s \in \N$, \begin{equation} \label{eq:formuleetape5} T^s \gf{G_s}{s!}=\gf{1}{V} \sum\limits_{j=0}^s (-1)^j T^{s+\mu-1} \gf{D^{s-j} R_{(j)}}{(s-j)!}  \left(\mathrm{com}\left(T^{\mu-1} R_{(0)}\right)\right)^{^T},\end{equation} où $V=\det(T^{\mu-1} R_{(0)}) \in \Oal_\K[z]$, et selon \eqref{eq:derivrjentier}, tous les termes de la somme sont à coefficients entiers algébriques à condition que $s+\mu-1 \leqslant \gf{M}{t+1}$.

Prenons pour $s \in \N$, $q_s$ le dénominateur des coefficients des coefficients des matrices $TG, T^2 \gf{G_2}{2!}, \dots, T^s \gf{G_s}{s!}$, comme dans la définition \ref{def:galochkinlarge}. Pour estimer $q_s$ sous la condition $s+\mu-1 \leqslant \gf{M}{t+1}$, il suffit donc d'obtenir une estimation de la maison des coefficients de $V_s$. C'est ce que permet de faire ce lemme :

\begin{lem} \label{lem:etape5}
Soient $U \in \Oal_\K[z]$, $V \in \Oal_\K[z]$, $W \in \K[z]$ tels que $U=VW$. Notons $V=\sum\limits_{i=0}^\ell v_i z^i$, alors pour tout $k$ tel que $v_k \neq 0$, on a $N_{\K/\Q}(v_k) W \in \Oal_\K[z]$.
\end{lem}

\begin{dem}
Introduisons la valuation de Gauss associée à un premier $\mathfrak{p}$ de $\Oal_\K$ : $$v_{\mathfrak{p}}\left(\sum\limits_{i=0}^{q} a_i z^i \right) :=\min\limits_{0 \leqslant i \leqslant q}(v_{\mathfrak{p}}(a_i)),$$ où $v_{\mathfrak{p}}$ est la valuation $\mathfrak{p}$-adique associée à l'idéal premier $\mathfrak{p}$ de $\Oal_\K$. En utilisant les propriétés de valuation, on a $v_{\mathfrak{p}}(U)=v_{\mathfrak{p}}(V)+v_{\mathfrak{p}}(W)$ donc $v_{\mathfrak{p}}(W) \geqslant - v_{\mathfrak{p}}(V)$ car, comme $U$ est à coefficients entiers algébriques, $v_{\mathfrak{p}}(U) \geqslant 0$.

Notons $S$ l'ensemble fini des premiers divisant tous les coefficients de $V$. Alors si $W=\sum\limits_{i=0}^{d} w_i z^i$, pour tout $i \in \{ 0, \dots, d \}$ et $\mathfrak{p} \in S$, on a $v_{\mathfrak{p}}(w_i) + v_{\mathfrak{p}}(V) \geqslant v_{\mathfrak{p}}(V)+v_{\mathfrak{p}}(W) \geqslant 0$, donc $ \prod\limits_{\mathfrak{p} \in S} \mathfrak{p}^{v_{\mathfrak{p}}(V)} (w_i) \subset \Oal_\K$. En particulier, si $v_{k} \neq 0$, comme $v_{k} \in \prod\limits_{\mathfrak{p} \in S} \mathfrak{p}^{v_{\mathfrak{p}}(V)}=\mathrm{pgcd}((v_0), \dots, (v_{\ell}))$, on a $$\forall i \in \{ 0, \dots, d \}, \quad (v_k) (w_i) \subset \Oal_\K.$$ D'où comme $N_{\K/\Q}(v_k) \in (v_k) \cap \Z$, on a $N_{\K/\Q}(v_k) W \in \Oal_\K[z]$.\end{dem}

Pour $W=\sum\limits_{i=0}^\ell w_i z^i \in \K[z]$, on définit $\sigma(W)=\max\limits_{0 \leqslant i \leqslant \ell} \house{w_i}$, la \textit{maison} de $W$. Selon le lemme \ref{lem:etape5} appliqué à la formule \eqref{eq:formuleetape5} avec $$U=\sum\limits_{j=0}^s (-1)^j T^{s+\mu-1} \gf{D^{s-j} R_{(j)}}{(s-j)!}  \left(\mathrm{com}\left(T^{\mu-1} R_{(0)}\right)\right)^{^T},$$ $V=\det(T^{\mu-1} R_{(0)})$ et $W=T^s \gf{G_s}{s!}$, on a donc ici \begin{equation} \label{eq:conclusionetape5}
\forall s \in \N \;\; \text{tel que} \;\; s+\mu-1 \leqslant \gf{M}{t+1}, \;  \; q_s \leqslant \sigma(V)^{\delta}.
\end{equation}

\bigskip

\textbf{Étape 6} : Majoration de la taille des coefficients de $\det(T^{\mu-1} R_{(0)}) \in \Oal_\K [z]$ en fonction de la maison de $Q$. 

Par commodité, on s'intéresse à $$\widetilde{V}=\det(\mathbf{P}, T \mathbf{P}_1, \dots, T^{\mu-1} \mathbf{P}_{\mu-1})=T^{\frac{\mu(\mu-1)}{2}} \det(R_{(0)})= T^{-\frac{\mu(\mu-1) }{2}} V.$$ Le lemme suivant nous assure que ce changement n'introduit qu'une constante multiplicative dépendant seulement de $G$ dans la majoration recherchée.

\begin{lem} \label{lem:etape6}
Soient $A, B \in \K[z]$ et $C=AB$, alors $\sigma(C) \leqslant (\deg(A)+\deg(B)+1) \sigma(A) \sigma(B)$.
\end{lem}

\begin{dem}
On écrit $A=\sum\limits_{i=0}^{p} a_i z^i$ et $B=\sum\limits_{i=0}^q b_i z^i$, de sorte que $C=\sum\limits_{j=0}^{p+q} c_j z^j$, avec, pour tout $j \in \{ 0, \dots, p+q \}$, $c_j=\sum\limits_{i=0}^{j} a_i b_{j-i}$.

Soit $\tau : \K \hookrightarrow \Qbar$ un plongement. Alors $$ |\tau(c_j)| \leqslant \sum_{i=0}^{j} |\tau(a_i)| |\tau(b_{j-i})| \leqslant (j+1) \sigma(A) \sigma(B) \leqslant (p+q+1) \sigma(A) \sigma(B),$$ si bien qu'en prenant le maximum sur $j$ et sur $\tau$, on obtient l'inégalité voulue.\end{dem}

Soit $m \in \{0, \dots, \mu-1\}$. Si $Q=\sum\limits_{i=0}^N q_i z^i$, alors $$\gf{Q^{(m)}}{m!}=\sum\limits_{i=0}^{N-m} \gf{(i+1) \dots (i+m)}{m!} q_{m+i} z^m=\sum\limits_{i=0}^{N-m} \binom{m+i}{i} q_{m+i} z^m.$$
De la majoration $\binom{m+i}{i} \leqslant 2^{m+i} \leqslant 2^{N}$, il s'ensuit que $\sigma\left(\gf{Q^{(m)}}{m!}\right) \leqslant 2^N \sigma(Q)$. Selon le lemme \ref{lem:etape6} appliqué à $A=T^m$ et $B=Q^{(m)}/m!$,

\begin{multline*}\sigma\left(T^m \gf{Q^{(m)}}{m!} \right) \leqslant 2^N \sigma(Q) \sigma(T^m)(mt+N-m+1) \\ \leqslant 2^N \sigma(Q) \sigma(T^m)((\mu-1)(t-1)+N+1) \leqslant c_1^N \sigma(Q) \sigma(T^m), 
\end{multline*} avec $c_1$ constante dépendant seulement de $G$ pour $N$ suffisamment grand. 

En appliquant $m$ fois le le lemme \ref{lem:etape6}, on obtient $$\sigma(T^m) \leqslant c_2^N \sigma(T) \sigma(T^{m-1}) \leqslant \dots \leqslant c_3^N \sigma(T)^m,$$ où $c_2$ et $c_3$ sont des constantes. Dans ce qui suit, les $c_i$ désigneront des constantes.

Ainsi, pour $N$ suffisamment grand, $$\sigma\left(T^m \gf{Q^{(m)}}{m!} \right) \leqslant c_4^N \sigma(Q) \sigma(T)^m \leqslant c_5^N \sigma(Q).$$

Soit $\theta_{N+M}$ le maximum des maisons des $N+M$ premiers coefficients du développement en série entière des $f_i$, et $d_{N+M}$ leur dénominateur commun. En répétant le raisonnement de la preuve du lemme \ref{lem:etape6}, on voit que la  maison de la partie polynomiale tronquée à l'ordre $N+tm$ de $T^m \gf{Q^{(m)}}{m!} (d_{N+M} f_i)$ est majorée par $$c_5^N \sigma(Q) (d_{N+M} \theta_{N+M})(N+tm+1) \leqslant c_6^N \sigma(Q) (d_{N+M} \theta_{N+M}).$$ Or, si $N+M-(\mu-1) \geqslant N+t(\mu-1)$, selon \eqref{eq:derivpadedenom}, cette partie polynomiale est $T^m \mathbf{P}_m$, donc, avec une extension de la notation $\sigma$ aux vecteurs colonnes, $$\forall m \in \{ 0, \dots, \mu-1 \}, \quad \sigma(T^m \mathbf{P}_m) \leqslant c_7^N \sigma(Q) d_{N+M} \theta_{N+M}.$$ On a $\widetilde{V}=\sum\limits_{\tau \in \mathfrak{S}_{\mu}} \varepsilon(\tau) \prod\limits_{j=0}^{\mu-1} T^j \mathbf{P}_{\tau(j),j}$. Pour $\tau \in \mathfrak{S}_{\mu}$, en appliquant le lemme \ref{lem:etape6} à $A=\mathbf{P}_{\tau(0),0}$ et $B=\prod\limits_{j=1}^{\mu-1} T^j \mathbf{P}_{\tau(j),j}$, on a 

$$\sigma\left(\prod\limits_{j=0}^{\mu-1} T^j \mathbf{P}_{\tau(j),j}\right) \leqslant \sigma(\mathbf{P}_{\tau(0),0}) \sigma\left(\prod\limits_{j=1}^{\mu-1} T^j \mathbf{P}_{\tau(j),j}\right) (\mu(N+t(\mu-1)+1)$$ en utilisant \eqref{eq:bornedegrePm}. En itérant le procédé, on obtient une constante $c_8$ telle que

$$\sigma\left(\prod\limits_{j=0}^{\mu-1} T^j \mathbf{P}_{\tau(j),j}\right) \leqslant c_7^{\mu N} \sigma(Q)^\mu d_{N+M}^\mu \theta_{N+M}^\mu c_8^N \leqslant c_9^N \sigma(Q)^\mu d_{N+M}^\mu \theta_{N+M}^\mu.$$ Donc par inégalité triangulaire $$\sigma(\widetilde{V}) \leqslant c_{10}^N \sigma(Q)^\mu d_{N+M}^\mu \theta_{N+M}^\mu,$$ si bien que, comme $V=T^{\mu(\mu-1)/2} \widetilde{V}$, $\sigma(V) \leqslant c_{11}^N \sigma(Q)^\mu d_{N+M}^\mu \theta_{N+M}^\mu$. D'où selon \eqref{eq:conclusionetape5},

\begin{equation} \label{eq:conclusionetape6}
\forall s \in \N, \;\; s+\mu-1 \leqslant \gf{M}{t+1}, \quad q_s \leqslant \sigma(V)^{\delta} \leqslant c_{12}^N \sigma(Q)^{\mu \delta} d_{N+M}^{\mu \delta} \theta_{N+M}^{\mu \delta}.
\end{equation}

\bigskip

\textbf{Étape 7} : Conclusion à l'aide d'un lemme diophantien.

Rappelons le lemme classique suivant dont une preuve peut être trouvée dans \cite[p. 37]{Siegel}.

\begin{lem}[Lemme de Siegel] \label{lem:siegel}
Soit $\K$ un corps de nombres. Considérons un système de $m$ équations linéaires
\begin{equation}\label{eq:systemesiegel}
 \forall 1 \leqslant i \leqslant m, \quad \sum\limits_{j=1}^{n} a_{ij} x_j=0, 
\end{equation} où $\forall i,j, a_{ij} \in \Oal_\K$. On note $A=\max\limits_{i,j} \house{a_{ij}}$. Alors si $n>m$, \eqref{eq:systemesiegel} a une solution non nulle $(x_j)_{1 \leqslant j \leqslant n} \in \Oal_\K^n$ vérifiant$$ \max\limits_{1 \leqslant j \leqslant n} \house{x_j} \leqslant c_1 (c_1 n A)^{\frac{m}{n-m}},$$ où $c_1 >0$ est une constante dépendant uniquement de $\K$. 
\end{lem}

Soit $s \in \N^*$. On choisit dorénavant $N$ et $M$ de la forme suivante : $N :=2\mu(t+1)(s+\mu-1)$ et $M :=N/(2\mu)=(t+1)(s+\mu-1) \in \N^*$. En particulier, on a bien $\gf{M}{t+1} \geqslant s+\mu-1$.

Alors l'équation de Padé $Q (d_{N+M} \mathbf{f})-\mathbf{P}=O\left(z^{N+M} \right)$ se traduit par un système linéaire de $\gf{\mu N}{2\mu}=\gf{N}{2}$ équations à $N+1$ inconnues (les coefficients de $Q$). Selon le lemme \ref{lem:siegel}, il existe une solution $Q \in \Oal_\K[z]$ telle que $$\sigma(Q) \leqslant c_{10} (c_{10} (N+1) \theta_{N+M})^{\frac{N/2}{N+1-N/2}} \leqslant c_{13}^N \theta_{N+M},$$ car $\gf{N}{2} \leqslant N+1-\gf{N}{2}$. 

\medskip

C'est à partir de maintenant que l'on va se servir des propriétés \textbf{b)} et \textbf{c)} de la définition \ref{def:gfonctionlarge}, qui sont propres aux $G$-fonctions \emph{au sens large}.

Soit $\varepsilon >0$. Puisque les composantes de $\mathbf{f}$ sont des $G$-fonctions \emph{au sens large}, on peut trouver une constante $c_{14}(\varepsilon)$ telles que $\theta_{k} \leqslant \left(k!\right)^{\varepsilon}$ et $d_{k} \leqslant \left(k!\right)^{\varepsilon}$ pour $k \geqslant c_{14}(\varepsilon)$.

On a $s=\left\lfloor N/(2n(t+1)) \right\rfloor-(\mu-1)$, de sorte que $M/(t+1) \geqslant s+\mu-1$ donc selon \eqref{eq:conclusionetape6}, \begin{equation} \label{eq:etape7} q_s \leqslant c_{12}^N (c_{13}^N \theta_{N+M})^{\mu \delta} d_{N+M}^{\delta \mu} \theta_{N+M}^{\delta \mu}. \end{equation}

D'une part, les termes géométriques $c_{i}^N$ peuvent être dominés à partir d'un certain rang qui dépend de $\varepsilon$ par $(N!)^{\varepsilon}$. D'autre part, on a $(N+M) \leqslant 2 N$. Or, si $\alpha >1$, selon la formule de Stirling, \begin{equation}\label{eq:stirling} \gf{(\alpha k)!}{(k!)^{\alpha}} \sim \gf{\sqrt{2 \pi \alpha k}}{(\sqrt{2 \pi k})^{\alpha}} \gf{\left(\frac{\alpha k}{e}\right)^{\alpha k}}{\left(\frac{k}{e}\right)^{\alpha k}}=r_k \alpha^{\alpha k}, \end{equation} avec $(r_N)_{N \in \N^*}$ une suite tendant vers $0$. Ainsi, si $k$ est assez grand, $(\alpha k)! \leqslant (k!)^{\alpha+1}$.

Par conséquent, en reprenant l'équation \eqref{eq:etape7}, on obtient que $q_s \leqslant (N!)^{c_{15} \times \varepsilon}$ pour $s$ plus grand qu'un certain rang dépendant de $\varepsilon$ et de la constante $c_{15}$.

Mais $N \leqslant 2\mu (t+1)(s+\mu-1) \leqslant 4\mu (t+1)s$ pour $s$ assez grand, donc à partir d'un certain rang, selon \eqref{eq:stirling}, $N! \leqslant (s!)^{4\mu (t+1)+1}$. Finalement, pour $s$ assez grand (relativement à $\varepsilon$), on a

\begin{equation} \label{eq:finalechudnovsky7} \forall s \geqslant s_0(\varepsilon), \;\; q_s \leqslant (s!)^{c_{16} \varepsilon}.\end{equation}

Il suffit d'appliquer ce résultat à $\varepsilon'=\gf{\varepsilon}{c_{16}}$ pour obtenir la majoration désirée.

\section{Démonstration du théorème \ref{th:acklarge} } \label{sec:acklarge}

Dans cette section, nous allons démontrer le théorème \ref{th:acklarge} en suivant l'esquisse fournie par André \cite[pp. 746--747]{AndregevreyII}.

\subsection{Reformulation de la condition de Galochkin}

On peut reformuler les conditions arithmétiques et analytiques définissant une $G$-fonction \emph{au sens large} à l'aide des valeurs absolues sur le corps de nombres $\K$.

\bigskip

\begin{prop}[André, \cite{AndregevreyII}] \label{prop:reformulationdefGfon}
Soit $f(z)=\sum\limits_{n=0}^{\infty} a_n z^n \in \K\llbracket z \rrbracket$. Les conditions \textbf{b)} et \textbf{c)} de la définition \ref{def:gfonctionlarge} équivalent à \begin{equation}\label{eq:G-alt}\gf{1}{n} \sum_{v} \log \max(1, |a_0|_v,\dots,|a_n|_v) = o(\log n),
 \end{equation}
 où la somme porte sur toutes les valeurs absolues (à équivalence près) $| \cdot |_v$ sur $\K$, finies et infinies.
\end{prop}

On rappelle que deux valeurs absolues $| \cdot |_v$ et $| \cdot |_{v'}$  sont dites équivalentes s'il existe $c >0$ tel que $| \cdot |_{v'}=| \cdot |_v^{c}$. Les valeurs absolues sur un corps de nombres $\K$ sont, à équivalence près, les valeurs absolues $\mathfrak{p}$-adiques (dites \emph{finies}) $|\cdot|_{\mathfrak{p}}$ pour $\mathfrak{p} \in \Spec(\Oal_{\K})$ et les valeurs absolues \emph{infinies} $|\cdot|_{\tau}$ pour tout plongement $\tau : \K \hookrightarrow \C$ définie par $$ |\zeta|_{\tau}=\begin{cases} |\tau(\zeta)|^{1/[\K:\Q]} \quad \mathrm{si} \; \tau(\K) \subset \R \\ |\tau(\zeta)|^{2/[\K:\Q]} \quad \mathrm{sinon}. \end{cases}$$
Pour la démonstration de la proposition \ref{prop:reformulationdefGfon}, on aura besoin du lemme technique suivant :

\begin{lem} \label{lem:techniquesuite}
\begin{enumerateth}
\item Soit $(u_n)_{n \in \N}$ une suite de réels et $g : \R_{+} \rightarrow \R_{+}$  une fonction croissante tendant vers $+\infty$ en l'infini. Alors la condition \begin{equation}
\max(0,u_n)=o(g(n))
\end{equation} est équivalente à la condition \begin{equation}\label{eq:techniquesuite2}
\max(0, u_0, \dots, u_n)=o(g(n)).
\end{equation}

\item Étant donnée $(u_n)_{n \in \N}$ une suite de réels positifs, on a \begin{equation} \label{eq:techniquesuite3} \forall \varepsilon>0,\; \exists n_0(\varepsilon) \in \N :\; \forall n \geqslant n_0(\varepsilon), \quad u_n \leqslant (n!)^{\varepsilon}\end{equation} si et seulement si $\max(0,\log u_n)=o(n \log n)$.
\end{enumerateth}
\end{lem}

\begin{dem}
\textbf{a)} Supposons que $\max(0,u_n)=o(g(n))$. 

Soit $\varepsilon >0$ et $n_0(\varepsilon) \in \N$ tel que $\forall n \geqslant n_0(\varepsilon), \; |u_n| \leqslant \varepsilon g(n)$. Comme $g$ est croissante, on a $$\forall n \geqslant n_0(\varepsilon), \;\; \max(0, u_{n_0(\varepsilon)}, \dots, u_{n}) \leqslant \varepsilon \max(g(n_0(\varepsilon)), \dots, g(n)) \leqslant \varepsilon g(n).$$
Soit $n_1(\varepsilon) \in \N$ tel que $$\forall n \geqslant n_1(\varepsilon), \;\; \max(0, u_0, \dots, u_{n_0(\varepsilon)-1}) \leqslant \varepsilon g(n).$$ Son existence est assurée par le fait que $g(x) \xrightarrow[x \rightarrow +\infty]{} + \infty$. Posons $n_2(\varepsilon) := \max(n_0(\varepsilon), n_1(\varepsilon))$. Alors pour tout $n \geqslant n_2(\varepsilon)$, $$\max(0, u_0, \dots, u_n) \leqslant \varepsilon g(n)$$ ce qui prouve  \eqref{eq:techniquesuite2}

La réciproque est évidente.

\textbf{b)} On a $\max(0, \log u_n)=o(n \log n)$ si et seulement si $\max(0, \log u_n)=o(\log n!)$ d'après la formule de Stirling, ce qui équivaut à $$\forall \varepsilon >0, \; \exists n_0(\varepsilon) \in \N : \; \forall n \geqslant n_0(\varepsilon), \;\; \log u_n \leqslant \varepsilon \log(n!),$$ ce qui est équivalent à \eqref{eq:techniquesuite3} en passant à l'exponentielle de part et d'autre de l'inégalité.
\end{dem}

\begin{dem}[de la proposition \ref{prop:reformulationdefGfon}]
On suppose sans perte de généralité que $\K$ est galoisien. Il est prouvé dans \cite[p. 225]{Dwork} que si $a_0, \dots, a_n \in \K$, alors \begin{equation} \label{eq:4denomhauteur}\sum\limits_{\mathfrak{p} \in \Spec(\Oal_{\K})}\sup\limits_{0 \leqslant i \leqslant n} \log^{+} |a_i |_{\mathfrak{p}} = \gf{1}{[\K:\Q]} \log\left( \den'(a_0, \dots, a_n)\right),\end{equation} où $\log^{+}(x)=\log \max(1,x)$ et $\den'(a_0,\dots,a_n)$ est la norme $N_{\K/\Q}(\mathfrak{a})$ du plus petit multiple commun  $\mathfrak{a}$ de $\mathfrak{a}_0, \dots, \mathfrak{a}_n$ dans le sens des anneaux de Dedekind. On a alors  
\begin{equation} \label{eq:equivdenetden'}\den(a_0,\dots,a_n) \leqslant \den'(a_0, \dots, a_n) \leqslant \den(a_0,\dots,a_n)^{[\K:\Q]}\end{equation} (cf \cite[p. 320]{Lepetit2}).  Ainsi, comme $$\sum_{v \; \text{valeur absolue}} \sup_{0 \leqslant i \leqslant n} \log^{+} |a_i|_v=\sum_{\mathfrak{p} \in \Spec(\Oal_{\K})} \sup_{0 \leqslant i \leqslant n} \log^{+}  |a_i|_{\mathfrak{p}}+\sum_{\tau \in \Gal(\K/\Q)}  \sup_{0 \leqslant i \leqslant n} \log^{+}|a_i|_{\tau},$$ la condition \eqref{eq:G-alt} est vérifiée si et seulement si $$\gf{1}{n}\sum_{\mathfrak{p} \in \Spec(\Oal_{\K})} \log^+ \sup_{0 \leqslant i \leqslant n}  |a_i|_{\mathfrak{p}}=o(\log n) \;\; \text{et} \;\; \gf{1}{n} \sum_{\tau \in \Gal(\K/\Q)} \log^{+} \sup_{0 \leqslant i \leqslant n}  |a_i|_{\tau}=o(\log n),$$ c'est-à-dire, selon \eqref{eq:4denomhauteur} et \eqref{eq:equivdenetden'}, si \begin{equation}\label{eq:estimlogden}\log \den(a_0,\dots,a_n)=o(n \log n) \;\; \text{et} \;\; \forall \tau \in \Gal(\K/\Q),\; \log 
\max(1, |\tau(a_0), \dots, |\tau(a_n)|)=o(n \log n).\end{equation} Selon le lemme \ref{lem:techniquesuite} \textbf{a)}, comme $$\log 
\max(1, |\tau(a_0), \dots, |\tau(a_n)|)=\max(0, \log|\tau(a_0)|, \dots, \log |\tau(a_n)|),$$ la condition \eqref{eq:estimlogden} équivaut à \begin{equation}\label{eq:estimlogden2}  \log \den(a_0,\dots,a_n)=o(n \log n) \;\; \text{et} \;\; \forall \tau \in \Gal(\K/\Q),\; 
\max(0, \log |\tau(a_n)|)=o(n \log n). \end{equation}

 Donc selon \eqref{eq:estimlogden2} et le lemme \ref{lem:techniquesuite} \textbf{b)} la condition \eqref{eq:G-alt} est bien équivalente aux conditions \textbf{b)} et \textbf{c)} de la définition~\ref{def:gfonctionlarge}.
\end{dem}

\begin{rqu}
De la même manière, on montre que $f(z)$ vérifie les conditions \textbf{b)} et \textbf{c)} de la définition \ref{def:gfonctionstrict} si et seulement si $$\gf{1}{n} \sum_{v} \log \max(1, |a_0|_v,\dots,|a_n|_v) = O(1).$$
\end{rqu}

\medskip

Par un raisonnement analogue, on peut reformuler la condition de Galochkin introduite dans la sous-section \ref{subsec:galochkinlarge} à l'aide d'une quantité $h(n,\mathfrak{p},G)$, introduite dans la définition \ref{def:hnpG} ci-dessous, définie en termes de valeurs absolues $p$-adiques sur $\K$.

Pour $\mathfrak{p} \in \Spec(\Oal_{\K})$, on rappelle que la \emph{valeur absolue de Gauss} associée à $|\cdot|_{\mathfrak{p}}$ est la valeur absolue non-archimédienne $$\fonction{|\cdot|_{\mathfrak{p},\mathrm{Gauss}}}{\K(z)}{\R}{\gf{\sum\limits_{i=0}^{N} a_i z^i}{\sum\limits_{j=0}^{M} b_j z^j}}{\gf{\max\limits_{0 \leqslant i \leqslant N} |a_i|_{\mathfrak{p}}}{\max\limits_{0 \leqslant j \leqslant M} |b_j|_{\mathfrak{p}}}.}$$ 

La valeur absolue $|\cdot|_{\mathfrak{p},\mathrm{Gauss}}$ induit naturellement une norme sur $\mathcal{M}_{\mu,\nu}(\K(z))$, définie pour tout $H=(h_{i,j})_{i,j} \in \mathcal{M}_{\mu,\nu}(\K(z))$ comme $\|H\|_{\mathfrak{p},\mathrm{Gauss}}=\max_{i,j} |h_{i,j}|_{\mathfrak{p},\mathrm{Gauss}}$. On l'appelle \emph{norme de Gauss}. Si $\mu=\nu$, $\mathcal{M}_{\mu}(\K(z))$ munie de $\|\cdot\|_{\mathfrak{p},\mathrm{Gauss}}$ est une algèbre normée.

\begin{defi}\label{def:hnpG} Pour tout idéal premier $\mathfrak{p}$ de $\Oal_{\K}$, la quantité $h(n,\mathfrak{p},G)$ est définie (cf \cite[p. 70]{Andre}) par $$ \forall n \in \N, \quad h(n,\mathfrak{p},G):= \gf{1}{n} \max_{m \leqslant n} \log^{+} \left\| \gf{G_m}{m!} \right\|_{\mathfrak{p}, \mathrm{Gauss}}.$$
Pour un opérateur différentiel $L \in \K(z)[\mathrm{d}/\mathrm{d}z]$, on pose $h(n,\mathfrak{p},L):=h(n,\mathfrak{p},A_{L})$, où $A_{L}$ est la matrice compagnon de $L$.
\end{defi}
' 

\begin{prop} \label{prop:galochkinlargereformulation}
Soit $G \in \mathcal{M}_{\mu}(\K(z))$. La condition de Galochkin \emph{au sens large} introduite dans la définition~\ref{def:galochkinlarge} est équivalente à $\sigma_n(G)=o(\log n)$, où $$\sigma_n(G) :=\sum_{\mathfrak{p} \in \Spec(\Oal_{\K})} h(n,\mathfrak{p},G).$$
\end{prop}

Dans le cas strict, la condition de Galochkin \emph{au sens strict} est équivalente au fait que $\sigma(G) <\infty$, où $$\sigma(G) := \limsup_{n \rightarrow +\infty} \sum_{\mathfrak{p} \in \Spec(\Oal_{\K})} h(n,\mathfrak{p},G)$$ est la \emph{taille} de $G$ (cf \cite[pp. 227--228]{Dwork}).

\begin{dem}[de la proposition \ref{prop:galochkinlargereformulation}]

On a vu dans la preuve de la proposition \ref{prop:reformulationdefGfon} que si $a_0, \dots, a_n \in \K$, alors $$\gf{1}{n}\sum\limits_{\mathfrak{p} \in \Spec(\Oal_{\K})}\sup\limits_{0 \leqslant i \leqslant n} \log^{+} |a_i |_{\mathfrak{p}}=o(\log n) \Leftrightarrow \log \den(a_0, \dots, a_n)= o(n \log n).$$ Donc, ici, en notant $q_s$ le plus petit dénominateur commun des coefficients des coefficients des matrices $\gf{G_m}{m!}$ pour $m \leqslant s$, on a $$\gf{1}{n} \sum_{\mathfrak{p} \in \Spec(\Oal_{\K})} \max_{m \leqslant n} \log^{+} \left\| \gf{G_m}{m!} \right\|_{\mathfrak{p}, \mathrm{Gauss}}=o(n \log n)$$ si et seulement si $\log q_s=o(n \log n)$, c'est-à-dire, comme $q_s \geqslant 1$, selon le lemme \ref{lem:techniquesuite} \textbf{b)}, si $$\forall \varepsilon >0, \exists s_0(\varepsilon) \in \N : \forall s \geqslant s_0(\varepsilon), \;\; q_s \leqslant s!^{\varepsilon},$$ ce qui n'est autre que la condition de Galochkin au sens large (définition \ref{def:gfonctionlarge}).
\end{dem}

\bigskip

\begin{rqu}
La preuve du théorème \ref{th:chudnovskylarge} présentée dans la partie \ref{subsec:preuvechudnovsky} permet de relier quantitativement $\sigma_n(G)$ et les quantités $$\sigma_{n}(\mathbf{f})=\gf{1}{n} \sum_{\mathfrak{p} \in \Spec(\Oal_{\K})} \max_{ m \leqslant n \atop 1 \leqslant i \leqslant \mu} |f_{i,m}|_{\mathfrak{p}} \quad \text{et} \quad \sigma_{n,\infty}(\mathbf{f})=\gf{1}{n} \sum_{\tau : \K \hookrightarrow \C} \max_{ m \leqslant n \atop 1 \leqslant i \leqslant \mu} |f_{i,m}|_{\tau} $$ encodant respectivement les conditions \textbf{b)} et \textbf{c)} de la définition \ref{def:gfonctionlarge} (Proposition \ref{prop:reformulationdefGfon}), les composantes du vecteur $\mathbf{f}={}^t (f_1(z), \dots, f_{\mu}(z))$ étant écrites sous la forme $f_i(z)=\sum\limits_{m=0}^{\infty} f_{i,m} z^m \in \K\llbracket z \rrbracket$.

On commence par appliquer le logarithme de part et d'autre de l'inégalité \eqref{eq:etape7} dans la partie \ref{subsec:preuvechudnovsky}. On obtient 
$$ \gf{\log q_s}{s} \leqslant N \gf{\log c_{12}}{s} + N \mu [\K:\Q] \gf{\log c_{13}}{s}+2 \mu [\K:\Q] \gf{\theta_{N+M}}{s}+\mu [\K:\Q] \gf{\log d_{N+M}}{s},$$ où $c_{12}$ et $c_{13}$ sont des constantes explicitables. 

Or, $N=2\mu(t+1)(s+\mu-1)$ et $M=(t+1)(s+\mu-1)$ donc en posant, pour tout $k$, $U_k=\gf{\log \theta_k}{k}$ et $V_k=\gf{\log d_k}{k}$, on a  \begin{align*}\gf{\log q_s}{s} &\leqslant 2 \mu(t+1)(1+o(1)) \log c_{12} + 2 \mu^2 [\K:\Q](t+1)(1+o(1)) \log c_{13}+ \\
& \qquad\qquad 2 \mu [\K:\Q] (2\mu+1)(t+1)(1+o(1)) U_{N+M}+\mu [\K:\Q](2\mu+1)(t+1)(1+o(1)) V_{N+M} \\
&\leqslant O(1)+\mu [\K:\Q] (2\mu+1)(t+1)(1+o(1))[2 U_{N+M}+V_{N+M}]. \end{align*}

Or, l'équation \eqref{eq:4denomhauteur} ci-dessus
implique que $$U_k=\gf{1}{k}\log \den(f_{i,m}, 1 \leqslant i \leqslant \mu, m \leqslant k) \leqslant \gf{[\K:\Q]}{k} \sum\limits_{\mathfrak{p} \in \Spec(\Oal_{\K})} \sup_{m \leqslant k \atop 1 \leqslant i \leqslant \mu} |f_{i,k}|_{\mathfrak{p}}=[\K:\Q] \sigma_k(\mathbf{f}).$$
De même, on montre que $\sigma_s(G) \leqslant \gf{1}{s} \log q_s$ et par ailleurs que $V_k \leqslant [\K:\Q] \sigma_{\infty,k}(\mathbf{f})$. Finalement, on obtient l'inégalité explicite \begin{multline}\label{eq:chudlargequantitatif}\sigma_s(G) \leqslant O(1)+ [\K:\Q]^2 \mu(2\mu+1)(t+1)(1+o(1))[2 \sigma_{(2\mu+1)(t+1)(s+\mu-1)}(f)+\\ \sigma_{\infty,(2\mu+1)(t+1)(s+\mu-1)}(f)],\end{multline} où les termes $O(1)$ et $o(1)$ peuvent être explicités.

Dans le cas des $G$-fonctions \emph{au sens strict}, une majoration analogue est déjà connue : voir la preuve de Dwork du théorème des Chudnovsky \cite[p. 299]{Dwork}.
\end{rqu}

\subsection{Condition de Bombieri \emph{au sens large}} \label{subsec:bombierilarge}

Le théorème de Katz affirme qu'un opérateur \emph{globalement nilpotent} (selon la terminologie de \cite[p. 95, p. 98]{Dwork}) est fuchsien à exposants rationnels en tout point de $\Pro^1(\Qbar)$ (voir \cite[p. 98]{Dwork}). Dans cette section, nous allons introduire une \emph{condition de Bombieri} au sens large (définition \ref{def:bombierifaible}) qui implique la nilpotence globale des opérateurs qui la vérifient. Nous verrons dans la sous-section \ref{subsec:andrebombierilarge} que les $G$-opérateurs \emph{au sens large} vérifient la condition de Bombieri \emph{au sens large}. Cette condition est donc adaptée au contexte des $G$-fonctions \emph{au sens large}.
 
 On se donne un corps de nombres $\K$. Dans toute la suite, pour tout idéal premier $\mathfrak{p}$  de $\Oal_{\K}$, on note $p(\mathfrak{p})$ le nombre premier engendrant l'idéal $\mathfrak{p} \cap \Z$.

Définissons d'abord une notion de densité qui sera utilisée par la suite (cf \cite[pp. 255-257]{Descombes}).

\begin{defi} \label{def:densites}
Soit $\mathcal{S} \subset \mathrm{Spec}(\Oal_\K)$. On note $N_{\K/\Q}(\mathfrak{p}) =\mathrm{Card}\left(\Oal_{\K}/\mathfrak{p}\right)$ la norme d'un idéal $\mathfrak{p}$ de $\Oal_{\K}$.

On dit que $\mathcal{S}$ admet $d \geqslant 0$ pour \emph{densité de Dirichlet} (ou \emph{densité analytique}) lorsque $$ \gf{-1}{\log(s-1)} \sum_{\mathfrak{p} \in \mathcal{S}} \gf{1}{N_{\K/\Q}(\mathfrak{p})^s} \longrightarrow d$$ quand $s$ tend vers $1$ pour $s$ réel, $s >1$. 
\end{defi}

\begin{defi}
Soit $\K$ un corps de nombres. Un système différentiel $y'=Gy$, $G \in \mathcal{M}_{\mu}(\K(z))$, est dit \emph{globalement nilpotent} si on a $R_{\mathfrak{p}}(G) > |p(\mathfrak{p})|_{\mathfrak{p}}^{1/(p(\mathfrak{p})-1)}$ pour tout idéal premier $\mathfrak{p}$ de $\Oal_{\K}$ en dehors d'un ensemble de densité de Dirichlet nulle.  
\end{defi}

En suivant la constuction faite au sens strict dans \cite{Dwork}, nous allons maintenant définir, en suivant \cite[p. 747]{AndregevreyII}, une condition de Bombieri adaptée à notre situation.

Si $G$ est une matrice à coefficients dans $\K(z)$ et $\mathfrak{p}$ est un idéal premier de $\Oal_{\K}$, on définit $R_{\mathfrak{p}}(G)$ comme le rayon de convergence d'une matrice fondamentale de solutions du système $y'=G y$ au voisinage du \emph{point générique $\mathfrak{p}$-adique} tel que construit dans \cite[pp. 92--97]{Dwork}.

\begin{defi}[André, \cite{AndregevreyII}, p. 747]\label{def:bombierifaible}
Soit, pour $n \in \N^*$, $\rho_n(G)=\sum\limits_{\mathfrak{p} \in \Spec(\Oal_{\K}) \atop p(\mathfrak{p}) \leqslant n} \log^+\left(\gf{1}{R_{\mathfrak{p}}(G)}\right)$. On dit que le système $y'=Gy$ satisfait la \emph{condition de Bombieri} au sens large quand $$\rho_n(G)=o(\log n).$$
\end{defi}

Dans le cas strict, la condition de Bombieri \emph{au sens strict} est $\rho(G) <\infty$, où $$\rho(G) := \sum\limits_{\mathfrak{p} \in \Spec(\Oal_{\K})} \log^+\left(\gf{1}{R_{\mathfrak{p}}(G)}\right) $$ est le \emph{rayon de convergence générique global inverse} de $G$ (cf \cite[pp. 226--227]{Dwork}).

\medskip

La proposition suivante illustre, comme annoncé au début de la section, l'intérêt de la condition de Bombieri : les systèmes différentiels la satisfaisant font partie de ceux auxquels on peut appliquer le théorème de Katz sur la rationalité des exposants.

\begin{prop}\label{prop:bombierifaible}
Si la matrice $G$ satisfait la condition de Bombieri \emph{au sens large}, alors $y'=Gy$ est un système différentiel globalement nilpotent.
\end{prop}

La démonstration fera usage du résultat suivant : 
\begin{prop}\label{prop:densdirnulle}
 Soit $\mathcal{S} \subset \Spec(\Z)$. Si $\sum\limits_{p \in \mathcal{S} \atop p \leqslant n} \gf{\log p}{p}=o(\log n)$, alors $\mathcal{S}$ a une densité de Dirichlet nulle.
\end{prop}

\begin{dem}
La preuve repose sur une transformation d'Abel de la série $\sum\limits_{p \in \mathcal{S}} p^{-s}$.

Pour tout $s >1$, on définit la fonction de classe $\mathcal{C}^1$ $$\fonction{g_s}{[1,+\infty[}{\R}{t}{\gf{1}{t^{s-1} \log t}.}$$

Selon \cite[Theorem 1, p. 3]{Tenenbaum}, on a pour tout $x \geqslant 0$, \begin{align}\label{eq:IPPabel}\sum_{p \in \mathcal{S} \atop p \leqslant x}\gf{1}{p^s} &=\sum_{p \in \mathcal{S} \atop p \leqslant x} \gf{\log p}{p} \gf{1}{p^{s-1} \log p}=T(x) g_s(x)-\displaystyle\int_{2}^x T(t) g_s'(t) \mathrm{d}t \notag\\
&=T(x) g_s(x)+(s-1)\displaystyle\int_{2}^x \gf{T(t)}{t^s \log(t)} \mathrm{d}t+\displaystyle\int_{2}^x \gf{T(t)}{t^s (\log(t))^2} \mathrm{d}t,\end{align} où $T(t)=\sum\limits_{p \in \mathcal{S} \atop p \leqslant t}  \gf{\log p}{p} = o(\log t)$ par hypothèse.

En particulier, $T(x) g_s(x) \xrightarrow[x \rightarrow +\infty]{} 0$ et $T g'_s$ est intégrable sur $[2,+\infty[$. Ainsi, en passant à la limite $x \rightarrow +\infty$, on obtient
\begin{equation}\label{eq:exprZetaS} Z(s) :=\sum_{p \in \mathcal{S}}\gf{1}{p^s}=(s-1)\displaystyle\int_{2}^{+\infty} \gf{E(t)}{t^s} \mathrm{d}t+\displaystyle\int_{2}^{+\infty} \gf{E(t)}{t^s \log(t)} \mathrm{d}t,\end{equation} où $E(t) := \gf{T(t)}{\log t}=o(1)$. On veut montrer que $\gf{Z(s)}{-\log(s-1)} \xrightarrow[s \rightarrow 1^{+}]{} 0$.

\begin{itemize}[label=\textbullet]
\item D'une part, soit $M >0$ tel que $\forall t \geqslant 2, \;\; 0 \leqslant E(t) \leqslant M$. Alors $$\forall s \in ]1,2[, \;\; 0 \leqslant \gf{(s-1)}{-\log(s-1)}\displaystyle\int_{2}^{+\infty} \gf{E(t)}{t^s} \mathrm{d}t \leqslant \gf{M(s-1)}{-\log(s-1)} \left[ \gf{t^{-s+1}}{-s+1} \right]_2^{+\infty}=\gf{-M}{2^{s-1} \log(s-1)} \xrightarrow[s \rightarrow 1^{+}]{} 0.$$

\item D'autre part, on va montrer que $$G(s) := \gf{1}{-\log(s-1)} \displaystyle\int_{2}^{+\infty} \gf{E(t)}{t^s \log(t)} \mathrm{d}t \xrightarrow[s \rightarrow 1^{+}]{} 0.$$
Par un changement de variable $u=(s-1)\log t$, on voit que pour tout $h \geqslant 2$, $$\displaystyle\int_{h}^{+\infty} \gf{\mathrm{d}t }{t^s \log(t)}=\mathrm{Ei}\left(-(s-1)\log h\right),$$ où $\mathrm{Ei} : x\mapsto - \displaystyle\int_{-x}^{+\infty} \gf{e^{-t}}{t} \mathrm{d}t$ est la fonction exponentielle intégrale. Par suite, il découle de l'équivalent (cf \cite[pp. 229--230]{AbramowitzStegun}) $\mathrm{Ei}(u) \sim_{u \rightarrow 0^{+}} \log(-u)$ que $$\forall h \geqslant 2, \quad \gf{1}{-\log(s-1)}\displaystyle\int_{h}^{+\infty} \gf{\mathrm{d}t}{t^s \log(t)} \xrightarrow[s \rightarrow 1^{+}]{} 1.$$

Soit $\varepsilon >0$. Comme $E$ tend vers $0$ en l'infini, on peut fixer $h \geqslant 2$ tel que $\forall t \geqslant h, E(t) \leqslant \gf{\varepsilon}{3}$. D'où \begin{align*}0 \leqslant G(s) &\leqslant \gf{M}{-\log(s-1)} \displaystyle\int_{2}^{h} \gf{\mathrm{d}t}{t^s \log(t)}+\gf{\varepsilon}{-\log(s-1)} \displaystyle\int_{h}^{+\infty} \gf{\mathrm{d}t}{t^s \log(t)} \\
& \leqslant \gf{M}{-\log(s-1)} \displaystyle\int_{2}^{h} \gf{\mathrm{d}t}{t \log(t)}+\varepsilon \gf{\mathrm{Ei}\left(-(s-1)\log h\right)}{-\log(s-1)},
\end{align*} car $t^s \geqslant t$ pour $t>1$. Le premier terme de la dernière somme tend vers $0$ et le second tend vers $1$ quand $s$ tend vers $1$, $s>1$, on peut donc fixer $\eta >0$ tel que pour tout $s \in ]1,1+\eta[$, $$\gf{M}{-\log(s-1)} \displaystyle\int_{2}^{h} \gf{\mathrm{d}t}{t \log(t)} \leqslant \gf{\varepsilon}{3} \;\; \text{et} \;\;  \gf{\mathrm{Ei}\left(-(s-1)\log h\right)}{-\log(s-1)} \leqslant 2$$ de sorte que $$\forall s \in ]1,1+\eta[, \quad 0 \leqslant G(s) \leqslant \varepsilon.$$ En d'autres termes, $G(s)$ tend vers $0$ quand $s$ tend vers $1^{+}$.
\end{itemize}
Selon ces deux points et \eqref{eq:exprZetaS}, on a bien $\lim\limits_{s \rightarrow 1^{+}}\gf{Z(s)}{\log(s-1)}=0$, soit selon la définition \ref{def:densites}, $\mathcal{S}$ a une densité de Dirichlet nulle.
\end{dem}

On peut à présent démontrer la proposition \ref{prop:bombierifaible}, en s'inspirant de la preuve du fait que la \emph{condition de Bombieri} au sens strict $\rho(G) < \infty$ implique que $y'=Gy$ est globalement nilpotent dans \cite[pp. 226--227]{Dwork} (voir aussi celle présentée par André, \cite[p. 77]{Andre}).

\begin{dem}[de la proposition \ref{prop:bombierifaible}]
Soit $\mathcal{S}$ l'ensemble des nombres premiers $p$ de $\Z$ tels qu'il existe $\mathfrak{p} \in \Spec(\Oal_\K)$ au-dessus de $p$ tel que $R_{\mathfrak{p}}(G) \leqslant |p|_{\mathfrak{p}}^{\frac{1}{p-1}}$.   Selon \cite[Proposition 5.1, p. 95]{Dwork}, $\mathcal{S}$ est l'ensemble des $p$ tels que $y'=(G \mod \mathfrak{p}) y$ est non nilpotent pour au moins un premier $\mathfrak{p}$ au-dessus de $p$. Fixons pour chaque $p \in \mathcal{S}$ un tel premier $\mathfrak{p}(p)$. Alors 
$$\rho_n(G) \geqslant \sum\limits_{p \in \mathcal{S} \atop p \leqslant n} \log^{+} \left(\gf{1}{R_{\mathfrak{p}(p)}(G)} \right) \geqslant \sum\limits_{p \in \mathcal{S} \atop p \leqslant n} \gf{1}{\delta} \gf{\log p}{p-1}.$$ En effet, $|p|_{\mathfrak{p}}=p^{-f_{\mathfrak{p}}/\delta}$ avec la normalisation choisie.

Grâce à l'hypothèse $\rho_n(G)=o(\log n)$, on a donc $\sum\limits_{p \in \mathcal{S} \atop p \leqslant n} \gf{\log p}{p-1}=o(\log n)$. Donc selon la proposition \ref{prop:densdirnulle}, $\mathcal{S}$ a une densité de Dirichlet nulle.

Or, si $\mathfrak{p} \in \Spec(\Oal_\K)$ est tel que $y'=(G \mod \mathfrak{p}) y$ est non nilpotent, alors $\mathfrak{p} \cap \Z \in \mathcal{S}$ donc l'ensemble $\mathcal{S}'$ des premiers vérifiant cette propriété est de densité de Dirichlet nulle. En effet, si $s>1$,
$$ \sum\limits_{\mathfrak{p} \in \mathcal{S}'} \gf{1}{N_{\K/\Q}(\mathfrak{p})^s}=\sum\limits_{p \in \Spec(\Z)} \sum\limits_{\mathfrak{p} \in \mathcal{S} \atop \mathfrak{p} \mid p} \gf{1}{p^{f_{\mathfrak{p}}s}} \leqslant \sum\limits_{p \in \mathcal{S}} \gf{\mathrm{Card}\{ \mathfrak{p} \in \mathcal{S}' : \mathfrak{p} \mid p\}}{p^s} \leqslant \delta \sum\limits_{p \in \mathcal{S}} \gf{1}{p^s},$$ car le nombre de premiers $\mathfrak{p}$ au-dessus de $p$ est borné par $\delta$ selon la formule $\sum\limits_{\mathfrak{p} \mid p} d_{\mathfrak{p}}=\delta$. Donc en divisant par $-\log(s-1)$ de part et d'autre de l'inégalité et en passant à la limite, on obtient que $\mathcal{S}'$ a une densité de Dirichlet nulle. Par définition, le système $y'=Gy$ est donc globalement nilpotent.\end{dem}

\begin{rqu}
L'assertion minimale pour que la preuve de la proposition \ref{prop:bombierifaible} fonctionne est la suivante. 

\textbf{Assertion \textbf{(A)}} : \og{} Soit $\mathcal{S} \subset \Spec(\Z)$. Si $\sum\limits_{p \in \mathcal{S} \atop p \leqslant n} \gf{\log p}{p}=o(\log n)$, alors $\mathcal{S}$ a une densité de Dirichlet $d < \gf{1}{2}$ \fg{}

En effet, en reprenant les notations de la preuve, si $\mathcal{S}$ vérifie \textbf{(A)}, alors $\mathcal{S}'$ est un ensemble de premiers de $\K$ divisant un ensemble de premiers de $\Z$ de densité de Dirichlet $<1/2$ et selon \cite[Remark 6.3, p. 100]{Dwork}, on peut appliquer le théorème de Katz au système $y'=Gy$, bien qu'il ne soit pas \emph{stricto sensu} globalement nilpotent, puisque $\mathcal{S'}$ n'a pas nécessairement dans ce cas une densité de Dirichlet nulle. 
\end{rqu}

\bigskip

\subsection{Théorème d'André-Bombieri \emph{au sens large}} \label{subsec:andrebombierilarge}

L'objectif de cette section est de démontrer la proposition \ref{prop:andrebombierilarge} ci-dessous, qui est l'analogue du théorème d'André-Bombieri \cite[pp. 228--234]{Dwork} pour les $G$-fonctions \emph{au sens large}. Ce résultat fait le lien entre les conditions de Galochkin et de Bombieri, présentées respectivement dans les parties \ref{subsec:galochkinlarge} et  \ref{subsec:bombierilarge}. Voir la proposition \ref{prop:ABbroad} dans la sous-section \ref{subsec:taillelarge} pour une version quantitative.

 \begin{prop} \label{prop:andrebombierilarge}
 Soit $G \in \mathcal{M}_{\mu}(\K(z))$. La condition de Galochkin \emph{au sens large} est vérifiée si et seulement si la condition de Bombieri \emph{au sens large} $\rho_n(G)=o(\log n)$ est vérifiée.
 \end{prop}

 Nous aurons besoin pour la démonstration du lemme suivant. On rappelle que si $\mathfrak{p} \in \Spec(\Oal_{\K})$, $p(\mathfrak{p})$ est défini comme l'unique premier positif de $\Z$ engendrant $\mathfrak{p} \cap \Z$.
 
 \begin{lem}\label{lem:lemmeABfaible}
 Si $G \in \mathcal{M}_{\mu}(\K(z))$, on a 
 $$\sum\limits_{\mathfrak{p} \in \Spec(\Oal_{\K})} h(n,\mathfrak{p},G) = \gf{1}{n}\sum\limits_{p(\mathfrak{p}) \leqslant n} h(n,\mathfrak{p},G)+C_{G},$$ où $C_{G}$ est une constante indépendante de $n$.
 \end{lem}
 
\begin{dem}
Soit $\mathfrak{p} \in \Spec(\Oal_{\K})$. On suppose que $\|G\|_{\mathfrak{p},\mathrm{Gauss}} \leqslant 1$. 
 A l'aide de la formule de récurrence $G_{s+1}=G_sG+G'_s$, comme $| \cdot |_{\mathfrak{p}}$ est une valeur absolue non archimédienne, on montre que $\|G_s\|_{\mathfrak{p},\mathrm{Gauss}} \leqslant 1$ pour tout $s \in \N$. Donc $$\gf{1}{n} \sum_{\mathfrak{p} \in \Spec(\Oal_{\K})} \left\|\gf{G_m}{m!}\right\|_{\mathfrak{p},\mathrm{Gauss}} = h(n,\mathfrak{p},G) \leqslant \gf{1}{n} \log^{+}\max_{s \leqslant n} \left\vert \gf{1}{s!}\right\vert_{\mathfrak{p}} \leqslant \gf{-1}{n} \log\min\left(1,\min_{s \leqslant n} |s!|_\mathfrak{p}\right).$$ Or si $s \leqslant n$, les diviseurs premiers $p$ de $s!$ vérifient $p \leqslant n$. Par suite, si $p(\mathfrak{p}) > n$, on a $|s!|_{\mathfrak{p}} =1$ et donc $h(n,\mathfrak{p},G)=0$.
 
 Si $\|G\|_{\mathfrak{p},\mathrm{Gauss}} >1$, on prend $\alpha \in \K$ tel que $|\alpha|_\mathfrak{p} \geqslant \|G\|_{\mathfrak{p},\mathrm{Gauss}}$. Ainsi, $\left\| G/\alpha\right\|_{\mathfrak{p},\mathrm{Gauss}} \leqslant 1$. De plus, si $\left\| G_s/\alpha^s \right\|_{\mathfrak{p},\mathrm{Gauss}} \leqslant 1$, alors $$\gf{G_{s+1}}{\alpha^{s+1}}=\left(\gf{G_s}{\alpha^s}\right)\left(\gf{G}{\alpha}\right)+\left(\gf{G_s}{\alpha^s} \right)' \gf{1}{\alpha}$$ est une somme de produits de termes dont la norme de Gauss associée à $\mathfrak{p}$ est inférieure à $1$, si bien que $\left\| G_{s+1}/\alpha^{s+1}\right\|_{\mathfrak{p},\mathrm{Gauss}} \leqslant 1$. Ceci montre par récurrence que $\|G_s/\alpha^s\|_{\mathfrak{p}, \mathrm{Gauss}}$ pour tout $s \in \N$.
 
Donc comme ci-dessus, il en découle que si $p(\mathfrak{p}) >n$, on a $\log^{+} \max\limits_{s \leqslant n} \left\|\gf{G_{s}}{\alpha^{s} s!}\right\|_{\mathfrak{p},\mathrm{Gauss}}=0$. Donc $$\log^{+} \max\limits_{s \leqslant n} \left\|\gf{G_{s}}{s!}\right\|_{\mathfrak{p},\mathrm{Gauss}} \leqslant \log^{+} \max_{s \leqslant n} \left\|\gf{G_{s}}{\alpha^{s} s!}\right\|_{\mathfrak{p},\mathrm{Gauss}}+s\log^{+}|\alpha|_{\mathfrak{p}},$$ de sorte que $h(n,\mathfrak{p},G) \leqslant \log |\alpha|_{\mathfrak{p},\mathrm{Gauss}}$, car $|\alpha|_{\mathfrak{p}} >1$.
 
 Comme on a $\|G\|_{\mathfrak{p},\mathrm{Gauss}} >1$ pour un nombre fini de premiers de $\K$, on peut prendre $\alpha$ tel que $|\alpha|_\mathfrak{p} \geqslant \|G\|_{\mathfrak{p},\mathrm{Gauss}}$ pour tout tel premier. On a alors $$\gf{1}{n} \sum_{\mathfrak{p} \in \Spec(\Oal_\K)} h(n,\mathfrak{p},G)=\gf{1}{n} \sum_{p(\mathfrak{p}) \leqslant n} h(n,\mathfrak{p},G)+\sum_{|G|_{\mathfrak{p},\mathrm{Gauss}} >1} \log|\alpha|_{\mathfrak{p},\mathrm{Gauss}},$$ ce qui est le résultat voulu, c'est-à-dire $$ \sum\limits_{\mathfrak{p} \in \Spec(\Oal_\K)} h(n,\mathfrak{p},G)=\sum\limits_{p(\mathfrak{p}) \leqslant n} h(n,\mathfrak{p},G)+O(1).$$
\end{dem}

\begin{dem}[de la proposition \ref{prop:andrebombierilarge}]

Selon la proposition \ref{prop:galochkinlargereformulation} et le  lemme \ref{lem:lemmeABfaible}, la condition de Galochkin \emph{au sens large} est vérifiée si et seulement si $$\sum\limits_{p(\mathfrak{p}) \leqslant n} h(n,\mathfrak{p},G)=o(\log n), \;\; n \rightarrow +\infty.$$

\begin{itemize}[label=\textbullet]
\item Selon \cite[p. 234]{Dwork}, on a pour tout $n \in \N^*$, $$\log^{+}\left(\gf{1}{R_{\mathfrak{p}}(G)}\right)=\limsup_{N \rightarrow +\infty} h(n,\mathfrak{p},G) \leqslant h(n,\mathfrak{p},G)+\gf{\log p(\mathfrak{p})}{n}+\gf{\log n}{n}.$$ Donc $$\rho_n(G) \leqslant \sum_{p(\mathfrak{p}) \leqslant n} h(n,\mathfrak{p},G)+ C_1 \sum_{p \in \Spec(\Z) \atop p \leqslant n} \gf{\log p}{n} + C_1 \gf{\log n}{n} \pi(n),$$ où $C_1$ est une constante majorant le nombre d'idéaux premiers $\mathfrak{p}$ de $\Oal_{\K}$ au-dessus d'un nombre premier $p$ donné, et $\pi$ est la fonction de comptage des nombres premiers.

Par le théorème des nombres premiers, on a donc $$\rho_n(G) \leqslant \sum_{p(\mathfrak{p}) \leqslant n} h(n,\mathfrak{p},G)+ 2C_1(1+o(1)),$$ de sorte que si $\sum\limits_{p(\mathfrak{p}) \leqslant n} h(n,\mathfrak{p},G)=o(\log n)$, alors $\rho_n(G)=o(\log n)$.

\item Réciproquement, l'équation (2.6) de \cite[p. 229]{Dwork} donne $$h(n,\mathfrak{p},G) \leqslant \log^{+}\left(\gf{1}{R_{\mathfrak{p}}(G)}\right)+ \mathbb{1}_{p(\mathfrak{p}) \leqslant n} \gf{d_{\mathfrak{p}}}{[\K:\Q]} (\mu-1) \gf{\log n}{n} + \gf{C_{\mathfrak{p}}}{n},$$ où $d_{\mathfrak{p}}$ est un entier tel que pour tout premier $p$, $\sum\limits_{\mathfrak{p} \mid p} d_{\mathfrak{p}} = [\K :\Q]$, et $C_{\mathfrak{p}}$ est une constante nulle pour tous les premiers $\mathfrak{p}$ sauf un nombre fini. 

On peut donc trouver une constante $C$ telle que, pour tout $n$, \begin{align*}
\sum\limits_{p(\mathfrak{p}) \leqslant n} h(n,\mathfrak{p},G) &\leqslant \rho_n(G)+\sum_{p \leqslant s} \left(\sum_{\mathfrak{p} \mid p} \gf{d_{\mathfrak{p}}}{[\K:\Q]} \right) (\mu-1)  \gf{\log n}{n} + \gf{C}{n} \\
&=\rho_n(G)+(\mu-1)  \gf{\log n}{n} \pi(n)+ \gf{C}{n}=\rho_n(G)+(\mu-1)(1+o(1)),
\end{align*}
\end{itemize}
puisque selon le théorème des nombres premiers, $\pi(n) \sim \gf{n}{\log n}$. Ainsi, si $\rho_n(G)=o(\log n)$, alors $\sum\limits_{p(\mathfrak{p}) \leqslant n} h(n,\mathfrak{p},G)=o(\log n)$.
\end{dem}

\subsection{Conclusion  de la démonstration du théorème \ref{th:acklarge}} \label{subsec:preuveacklarge}

Nous pouvons à présent synthétiser les différents résultats des sections précédentes pour prouver le théorème \ref{th:acklarge}.

Soit $\K$ un corps de nombres, $f(z)=\sum\limits_{n=0}^{\infty} a_n z^n \in \K\llbracket z \rrbracket$ une $G$-fonction et $L \in \K\left[z,\gf{\mathrm{d}}{\mathrm{d}z}\right]$ un opérateur d'ordre minimal non nul tel que $L(f(z))=0$.

\begin{itemize}[label=\textbullet]
 \item Selon le théorème \ref{th:chudnovskylarge} (partie \ref{subsec:galochkinlarge}), la matrice compagnon $A_L$ de $L$ vérifie la condition de Galochkin \emph{au sens large}.
 \item La proposition \ref{prop:andrebombierilarge} (partie \ref{subsec:andrebombierilarge}) implique alors que $A_L$ vérifie la condition de Bombieri \emph{au sens large} (définition \ref{def:bombierifaible}).
 \item Selon la proposition \ref{prop:bombierifaible} (partie \ref{subsec:bombierilarge}), $L$ est donc un opérateur différentiel globalement nilpotent. Par le théorème de Katz \cite[p. 98]{Dwork}, on en déduit que $L$ est donc un opérateur fuchsien à exposants rationnels en tout point de $\Pro^1(\Qbar)$. Ceci conclut la preuve du théorème \ref{th:acklarge}.
\end{itemize}

\section{Existence d'une base de type ACK au voisinage d'un point singulier} \label{sec:baseACK}

Le théorème d'André-Chudnovsky-Katz dans le cas des $G$-fonctions \emph{au sens strict} est un résultat plus fort que l'analogue strict du théorème \ref{th:acklarge}.

\begin{Th}[André-Chudnovsky-Katz]\label{th:ack}
 Soit $f(z)$ une $G$-fonction \emph{au sens strict} et $L \in \Qbar(z)[\mathrm{d}/\mathrm{d}z] \setminus \{ 0 \}$ un opérateur différentiel tel que $L(f(z))=0$ et d'ordre minimal $\mu$ pour $f$. 

Alors au voisinage de tout $\alpha \in \mathbb{P}^1(\Qbar)$, il existe une base de solutions de $Ly(z)=0$ de la forme $$(f_1(z-\alpha), \dots, f_{\mu}(z-\alpha)) (z-\alpha)^{C_{\alpha}},$$ où  $C_{\alpha} \in \mathcal{M}_{\mu}(\Qbar)$ est triangulaire supérieure à valeurs propres dans $\Q$  et les $f_i(u) \in \Qbar\llbracket u \rrbracket$ sont des $G$-fonctions \emph{au sens strict}.
\end{Th}

Le but de cette partie est d'en montrer l'analogue \emph{au sens large}, qui n'avait pas été mentionné par André dans \cite{AndregevreyII}. 

\begin{Th}\label{th:complementACKlarge}
Le théorème \ref{th:ack} reste vrai si l'on remplace partout \og strict \fg{} par \og large \fg{}.
\end{Th}

Les valeurs propres de la matrice $C_{\alpha}$ sont, modulo $\Z$, les exposants de $L$ en $\alpha$, le théorème~\ref{th:acklarge} implique donc qu'elles sont rationnelles. Il reste donc à prouver que les $f_i$ sont des $G$-fonctions.

La partie suivante permet de montrer un point essentiel de la preuve du théorème \ref{th:complementACKlarge} (section \ref{subsec:preuvecomplementACK}), puisqu'il permet de ramener l'étude au voisinage du point $0$ dans l'étape 4.

\subsection{Invariance des $G$-opérateurs \emph{au sens large} par changement de variable} \label{subsec:GoplargeInvchgvar}

Le but de cette partie est de démontrer la proposition suivante affirmant que la condition de Galochkin \emph{au sens large} est invariante par changement de variable $u=z-\alpha$ si $\alpha \in \Qbar$ ou $u=z^{-1}$ en l'infini. 

\begin{prop} \label{prop:galochkinLargeInvChgeVar}
\begin{enumerateth}
\item Soit $G \in \mathcal{M}_{\mu}(\K(z))$ et $\alpha \in \Pro^1(\Qbar)$. On définit $G_{\alpha} \in \mathcal{M}_{\mu}(\K(z))$ la matrice telle que, pour tout vecteur colonne $y$, $$y'(z)=G(z)y(z) \ssi \widetilde{y}'(u)=G_{\alpha}(u) \widetilde{y}(u),$$ où $\widetilde{y}(u) :=y(u-\alpha)$ si $\alpha \neq \infty$ et $\widetilde{y}(u)=y(u^{-1})$ si $\alpha=\infty$.  Alors la condition de Galochkin \emph{au sens large} est vérifiée par $G_{\alpha}$ si et seulement si elle est vérifiée par $G$.
\item Soit $L \in \K(z)\left[\mathrm{d}/\mathrm{d}z\right]$. On définit $L_{\alpha} \in \K(u)\left[\mathrm{d}/\mathrm{d}u\right]$ tel que, pour toute fonction $f(z)$, $L(f(z))=0 \ssi L_{\alpha}(\widetilde{f}(u))$, où $\widetilde{f}(u)=f(u-\alpha)$ si $\alpha \neq \infty$ et $\widetilde{f}(u)=f(u^{-1})$ si $\alpha=\infty$.

Alors $L_{\alpha}$ est un $G$-opérateur \emph{au sens large} si et seulement si $L$ est un $G$-opérateur \emph{au sens large}
\end{enumerateth}
\end{prop}

\begin{rqu}
En reproduisant la même preuve \emph{mutatis mutandis}, on prouve l'énoncé analogue concernant la condition de Galochkin \emph{au sens strict}.
\end{rqu}

Nous aurons besoin pour prouver la proposition \ref{prop:galochkinLargeInvChgeVar} du lemme \ref{lem:galochkininvariantequivalence} ci-dessous, dû à André, qui montre que la condition de Galochkin est invariante par équivalence de systèmes différentiels. Commençons par rappeler cette notion d'équivalence (cf \cite[p. 120]{Sauloy}) :
\begin{defi}
Soient $A, B \in \mathcal{M}_{\mu} \left(\Qbar((z))\right)$. On dit que les systèmes différentiels $y'=Ay$ et $y'=By$ sont équivalents sur $\Qbar(z)$ s'il existe $P \in \GL_{\mu}(\Qbar(z))$ tel que $B=P[A]$, où $P[A]=PAP^{-1}+P'P^{-1}$. Ceci est équivalent à dire que $y \mapsto Py$ réalise une bijection de l'ensemble des solutions de $y'=Ay$ sur l'ensemble des solutions de $y'=By$.
\end{defi}

\begin{lem}[André, \cite{Andre}, Lemma 1 p. 71] \label{lem:galochkininvariantequivalence}
Soient $A, B \in \mathcal{M}_{\mu}\left(\Qbar(z)\right)$ définissant deux systèmes différentiels équivalents $y'=Ay$ et $y'=By$ sur $\Qbar(z)$. Alors $A$ satisfait la condition de Galochkin \emph{au sens strict} (resp. large) si et seulement si $B$ satisfait la condition de Galochkin \emph{au sens strict} (resp. large).
\end{lem}

André a déjà prouvé le cas strict et nous donnons ici une preuve différente de celle de \cite[pp. 71--72]{Andre}, qui utilise la notion de module différentiel. 

\begin{dem}[du lemme \ref{lem:galochkininvariantequivalence}]
Soit $P \in \mathrm{GL}_{\mu}(\Qbar(z))$ tel que $y \mapsto P y$ induit une bijection de l'ensemble des solutions de $y'=Ay$ sur l'ensemble des solutions de $y'=By$.

Si $y'=Ay$ et $s \in \N$, on a, $y^{(s)}=A_s y$ et $(Py)^{(s)}=B_s Py$. Or, $$B_s Py= (Py)^{(s)}=\sum_{k=0}^s \binom{s}{k} P^{(s-k)} y^{(k)}=\sum_{k=0}^s \binom{s}{k} P^{(s-k)} A_k y.$$ L'égalité étant valable pour toute solution $y$ de $y'=Ay$, on en déduit $$B_s=\sum_{k=0}^s \binom{s}{k} P^{(s-k)} A_k P^{-1}.$$

Soit $T(z) \in \Oal_{\Qbar}[z]$ tel que $T(z) A(z) \in \mathcal{M}_{\mu}(\Oal_{\Qbar}[z])$. Sans perte de généralité, on peut supposer que $T(z) P(z) \in \mathcal{M}_{\mu}(\Oal_{\Qbar}[z])$ et $T(z) P^{-1}(z) \in \mathcal{M}_{\mu}(\Oal_{\Qbar}[z])$. On montre par récurrence sur $\ell$ que \begin{equation}\label{eq:denomderiveelieme} \forall \ell \in \N, \quad T^{\ell+1} \gf{P^{(\ell)}}{\ell!} \in \mathcal{M}_{\mu}(\Oal_{\Qbar}[z]).\end{equation}

En effet, c'est clair pour $\ell=0$ et si \eqref{eq:denomderiveelieme} est vraie au rang $\ell-1$ pour un certain $\ell \geqslant 1$, alors par la formule de Leibniz,

$$T^\ell \gf{(TP)^{(\ell)}}{\ell!}=T^{\ell+1} \gf{P^{(\ell)}}{\ell!}+\sum_{k=0}^{\ell-1} T^{\ell-(k+1)} \gf{T^{(\ell-k)}}{(\ell-k)!} T^{k+1} \gf{P^{(k)}}{k!}.$$
Le membre de gauche est à coefficients dans $\Oal_{\Qbar}[z]$ par le cas $\ell=0$ et les termes d'ordre $k \leqslant \ell-1$ du membre de droite le sont également par hypothèse de récurrence. On obtient donc $T^{\ell+1} \gf{P^{(\ell)}}{\ell!} \in \mathcal{M}_{\mu}(\Oal_{\Qbar}[z])$, ce qui conclut la récurrence.

\bigskip

Ainsi, si $q_s$ est le dénominateur des coefficients des coefficients de $TA, \gf{T^2 A_2}{2!}, \dots, T^{s} \gf{A_s}{s!}$, et $\widetilde{T}=T^3$, on a $$\forall 1 \leqslant \ell \leqslant s, \quad q_s \gf{\widetilde{T}^\ell B_\ell}{s!}=\sum_{k=0}^\ell T^{2\ell-2} \gf{T^{\ell-k+1} P^{(\ell-k)}}{(\ell-k)!} \left( q_\ell \gf{T^k A_k}{k!} \right) TP^{-1} \in \mathcal{M}_{\mu}(\Oal_{\Qbar}[z]).$$ Donc, en notant $\widetilde{q}_s$ le dénominateur des coefficients des coefficients des matrices $\widetilde{T} B, \widetilde{T}^2 \gf{B_2}{2!}, \dots, \widetilde{T}^{s} \gf{B_s}{s!}$, on obtient que $\widetilde{q}_s$ divise $q_s$ pour tout entier $s$.

On en déduit que si $A$ vérifie la condition de Galochkin \emph{au sens strict} (resp. \emph{large}), il en va de même de $B$. Par symétrie de la relation d'équivalence des systèmes différentiels, on obtient l'implication réciproque, ce qui achève la preuve. \end{dem}

La preuve de la proposition \ref{prop:galochkinLargeInvChgeVar} a déjà été donnée par André \cite[p. 83]{Andre} dans le cas où $\alpha \in \Oal_{\Qbar}$, nous donnons une preuve complète par souci d'exhaustivité.

\begin{dem}[de la proposition \ref{prop:galochkinLargeInvChgeVar}]

\textbf{a)} Soit $\alpha \in \Qbar$. On a $\forall s \in \N, G_{\alpha,s}(u)=G_s(u+\alpha)$. De plus, le polynôme $T_{\alpha}(u):=d_{\alpha}^t T(u+\alpha) \in \Oal_{\Qbar}[u]$, où $d_{\alpha} \alpha \in \Oal_{\Qbar}$, vérifie $T_{\alpha}(u) G_{\alpha}(u) \in \Qbar[u]$.

On voit que pour tout $s \in \N^*$ et $p \leqslant s$, $$q_s \gf{T_{\alpha}(u)^p G_{\alpha,p}(u)}{p!}=d_{\alpha}^{pt} \left(\gf{q_s T^p G_p}{p!}\right)(u+\alpha) \in \mathcal{M}_{\mu}\left(\Oal_{\Qbar}[u]\right)$$ puisque $q_s T(z)^p G_p(z)/p! \in \mathcal{M}_{\mu}\left(\Oal_{\Qbar}[z]\right)$ et on montre par récurrence que les coefficients de $T^p G_p$ sont de degrés au plus $pt$.

Donc si $G$ satisfait la condition de Galochkin \emph{au sens large}, $G_{\alpha}$ également. En remarquant que $G=\left(G_{\alpha}\right)_{-\alpha}$, on en déduit l'implication réciproque.

\bigskip

Le résultat au point à l'infini découle essentiellement de \cite[Lemma 7, p. 322]{Lepetit2}. Remarquons que si $T(z) G(z) \in \mathcal{M}_{\mu}\left(\Oal_{\Qbar}[z]\right)$, alors $T_{\infty}(u):=u^{t+2} T(u^{-1}) \in \Oal_{\Qbar}[u]$ vérifie $T_{\infty}(u) G_{\infty}(u) \in \mathcal{M}_{\mu}\left(\Oal_{\Qbar}(u)\right)$, de sorte que, selon \eqref{eq:tsgspoly}, $T_{\infty}^s G_{\infty,s} \in \mathcal{M}_{\mu}\left(\Oal_{\Qbar}[u]\right)$ pour tout entier $s$.

On a montré dans \cite[Lemma 7, p. 322]{Lepetit2} que
 \begin{equation} \label{eq:4formuleGinfty} \forall s \in \N^*, \quad G_{\infty,s}(u)=(-1)^{s} \sum\limits_{k=1}^{s} \gf{c_{s,k}}{u^{s+k}} G_k\left(\gf{1}{u}\right),\end{equation} où \begin{equation} \label{eq:4formulecoeffsAtilde} \forall s \in \N^*,\quad  \forall 1 \leqslant k \leqslant s, \quad c_{s,k}=\binom{s-1}{s-k} \gf{s!}{k!} \in \Z.\end{equation}

Ainsi, si $s \in \N^*$ et $p \leqslant s$, on a, selon \eqref{eq:4formuleGinfty},  \begin{align*}
q_s \gf{T_{\infty}(u)^p G_{\infty,p}(u)}{p!} &=(-1)^p \sum_{\ell=1}^p \gf{c_{p,\ell}}{p!} u^{-(p+\ell)} \left(u^2 u^t T\left(\gf{1}{u} \right)\right)^p q_s G_\ell\left(\gf{1}{u}\right) \\
&= (-1)^p \sum_{\ell=1}^p \binom{p-1}{p-\ell} u^{p-\ell} u^{(p-\ell)t} T\left(\gf{1}{u} \right)^{p-\ell} q_s u^{\ell t} \left(\gf{T^\ell G_\ell}{\ell!}\right)\left(\gf{1}{u} \right) \in \mathcal{M}_{\mu}\left(\Oal_{\Qbar}[u]\right)
\end{align*} de sorte qu'en notant $q_{\infty,s}$ le dénominateur des coefficients des coefficients des matrices $$T_{\infty} G_{\infty}, \gf{T_{\infty}^2 G_{\infty,2}}{2!}, \dots, \gf{T_{\infty}^s G_{\infty,s}}{s!},$$ on a $\forall s \in \N^*, q_{\infty,s} \leqslant q_s$. Donc si la condition de Galochkin \emph{au sens large} est vérifiée par $G$, elle l'est également par $G_{\infty}$. Puisque $\left(G_{\infty}\right)_{\infty}=G$, on en déduit l'implication réciproque.

\bigskip

\textbf{b)} Soit $L \in \Qbar(z)[\mathrm{d}/\mathrm{d}z]$, on note $G=A_{L}$ la matrice compagnon de $L$.

Si $\alpha \in \Qbar$, la matrice $G$ vérifie $G_{\alpha}=A_{L_{\alpha}}$, donc selon \textbf{a)}, $L$ est un $G$-opérateur \emph{au sens large} si et seulement si $L_{\alpha}$ l'est.

\medskip

Si $\alpha=\infty$, selon le point \textbf{b)} de \cite[Lemma 7, p. 322]{Lepetit2}, le système $h'=G_{\infty} h$ est équivalent au système $h'=A_{L_{\infty}} h$. La conclusion voulue découle ainsi du lemme \ref{lem:galochkininvariantequivalence}.

Supposons que $L$ est un $G$-opérateur \emph{au sens large}. Selon \textbf{a)}, $G_{\infty}$ vérifie  alors la condition de Galochkin \emph{au sens large}. Il en va donc de même de $A_{L_{\infty}}$ selon le lemme \ref{lem:galochkininvariantequivalence}. Par conséquent, $L_{\infty}$ est un $G$-opérateur \emph{au sens large}. 
\end{dem}

\subsection{Démonstration du théorème \ref{th:complementACKlarge}} \label{subsec:preuvecomplementACK}

Le but de cette section est de prouver le théorème \ref{th:complementACKlarge} en adaptant la démonstration donnée dans le cas strict par Dwork dans \cite[pp. 234--248]{Dwork}, et en utilisant également les résultats de la section \ref{subsec:GoplargeInvchgvar}.

On considère $\K$ un corps de nombres tel que $G \in \mathcal{M}_{\mu}(\K(z))$ et on se donne $C \in \mathcal{M}_{\mu}(\Qbar)$ une matrice à valeurs propres rationnelles, et $Y(z)=\sum\limits_{m=0}^{\infty} Y_m z^m \in \mathcal{M}_{\mu}(\K\llbracket z \rrbracket) \cap \mathrm{GL}_{\mu}\big(\K((z))\big)$ telle que $Y(z) z^{C}$ est une matrice de solutions du système $y'=Gy$. Quitte à remplacer $\K$ par une extension de degré supérieur, on peut supposer que $C \in \mathcal{M}_{\mu}(\K)$.

Commençons par démontrer que les coefficients de $Y$ vérifient les points \textbf{a)} et \textbf{b)} de la définition \ref{def:gfonctionlarge}.

\begin{itemize}[label=\textbullet]
\item D'abord, selon \cite[Corollary, p. 109]{Andre}, si $\widetilde{Y}(z)=Y(z)z^C$, avec $Y(z) \in \mathcal{M}_{\mu}(\Qbar\llbracket z \rrbracket)$ et $C \in \mathcal{M}_{\mu}(\Qbar)$, vérifie $\widetilde{Y}'=G\widetilde{Y}$, $G \in \mathcal{M}_{\mu}(\Qbar(z))$, alors chaque coefficient de $Y(z)$ est solution d'une équation différentielle non nulle à coefficients dans $\Qbar(z)$ (qui plus est d'ordre inférieur à $\mu^2$) . Autrement dit, les coefficients de $Y(z)$ satisfont le point \textbf{a)} de la définition \ref{def:gfonctionlarge}. Ce fait est indépendant du contexte des $G$-fonctions (voir \cite[Remark 2, p. 53]{Andre}).

\item Soit $\tau : \K \hookrightarrow \C$  un plongement. On note $R_{\tau}(Y)$ le rayon de convergence de la série $Y^{\tau} := \sum\limits_{m=0}^{\infty} \tau(Y_m) z^m$. Alors \cite[Equation (3.5), p. 242]{Dwork} implique que $$R_{\tau}(Y) \geqslant \rho:= \min_{\zeta \neq 0, \infty} |\zeta|_{\tau} >0,$$ où le minimum porte sur les singularités non apparentes du système $y'=Gy$. D'où selon la formule d'Hadamard-Cauchy, $$\limsup\limits_{n \rightarrow +\infty} \|\tau(Y_m)\|_{\infty}^{1/n} = \gf{1}{R_{\tau}(Y)} \leqslant \gf{1}{\rho}.$$ Ainsi, pour tout entier $m$, les maisons des coefficients de $Y_m$ sont bornées par $\kappa^{m+1}$, où $\kappa >0$ est une constante. On rappelle que la maison de $\alpha \in \Qbar$ est $\house{\alpha}=\max\limits_{\sigma \in \Gal(\Qbar/\Q)} |\sigma(\alpha)|$. Ceci prouve que les coefficients de $Y(z)$ vérifient la condition \textbf{b)} de la définition \ref{def:gfonctionlarge}.
\end{itemize}

\medskip

Le reste de la démonstration est consacré à la vérification de la condition arithmétique \textbf{c)} de la définition \ref{def:gfonctionlarge}.

\textbf{Étape 1.} Pour tout idéal premier $\mathfrak{p}$ de $\Oal_{\K}$, on désigne par $R_{\mathfrak{p}}(Y)$ le rayon de convergence pour la norme $\| \cdot \|_{\mathfrak{p}}$ de la série $Y(z)=\sum\limits_{m=0}^{\infty} Y_m z^m$. On note alors $$\rho_n(Y):=\sum\limits_{\mathfrak{p} \in \Spec(\Oal_{\K}) \atop p(\mathfrak{p}) \leqslant n} \log^{+}\left(\gf{1}{R_{\mathfrak{p}}(Y)}\right),$$ où $p(\mathfrak{p})$ est l'unique premier positif de $\Z$ engendrant $\mathfrak{p} \cap \Z$.

Le but de cette étape est d'établir un lien entre $\rho_n(Y)$ et de la quantité $\rho_n(G)$ introduite dans la définition \ref{def:bombierifaible} en suivant la démonstration de \cite[Theorem 3.3, p. 238]{Dwork}.

Dwork a obtenu dans \cite[Equation (3.4), p. 240]{Dwork}, sous les hypothèses sur $C$ et $Y$ ci-dessus, que pour tout premier $\mathfrak{p}$ de $\Oal_{\K}$,  $$\log^{+}\left(\gf{1}{R_{\mathfrak{p}}(Y)} \right) \leqslant \mu^2 \log\left(\gf{1}{R_{\mathfrak{p}}(G)}\right)+(\mu^2+1) \sum\limits_{\zeta \neq 0,\infty} \log^{+}\left| \gf{1}{\zeta}\right|_{\mathfrak{p}},$$ où la dernière somme porte sur toutes les singularités non apparentes $\zeta \not\in \{0,\infty\}$ du système $y'=Gy$.

Donc en sommant sur tous les premiers $\mathfrak{p}$ tels que $p(\mathfrak{p}) \leqslant n$, on obtient $$\rho_n(Y) \leqslant \mu^2 \rho_n(G)+(\mu^2+1) \sum\limits_{\zeta \neq 0,\infty} \sum_{\mathfrak{p} \in \Spec(\Oal_{\K}) \atop p(\mathfrak{p}) \leqslant n} \log^{+}\left| \gf{1}{\zeta}\right|_{\mathfrak{p}} \leqslant  \mu^2 \rho_n(G)+\beta,$$ où $\beta :=(\mu^2+1) \sum\limits_{\zeta \neq 0,\infty} \sum\limits_{\mathfrak{p} \in \Spec(\Oal_{\K})} \log^{+}\left| \gf{1}{\zeta}\right|_{\mathfrak{p}}$ est une constante ne dépendant que de $G$. La constante $\beta$ est finie, car pour tout $\alpha \in \K$, $|\alpha|_{\mathfrak{p}}=0$ pour tout premier $\mathfrak{p}$ sauf un nombre fini.

Ainsi, si $G$ vérifie la condition de Bombieri \emph{au sens large} $\rho_n(G)=o(\log n)$, on a $\rho_n(Y)=o(\log n)$.

\bigskip

\textbf{Étape 2.} On note $$\sigma_n(Y) :=\gf{1}{n} \sum\limits_{\mathfrak{p} \in \Spec(\Oal_{\K}) \atop p(\mathfrak{p}) \leqslant n} \sup\limits_{m \leqslant n} \log^{+} |Y_m|_{\mathfrak{p}}.$$ Le but de cette étape est de lier $\sigma_n(Y)$ et $\rho_n(Y)$.

Soit $N \in \N^*$ un dénominateur commun des valeurs propres de $C$, qui sont par hypothèse rationnelles. Alors Dwork prouve dans \cite[p. 247]{Dwork} qu'on peut trouver des constantes positives $\ell_0$ et $h_0$ telles que si $\mathfrak{p}$ est en dehors d'un ensemble fini de premiers $\mathcal{S} \subset \Spec(\Oal_{\K})$, on a pour tout $n \in \N^*$, \begin{multline}\label{eq:Dworkp247}\gf{1}{n} \sup\limits_{m \leqslant n} \log^{+} |Y_m|_{\mathfrak{p}} \leqslant \gf{1}{n}\left(n+\gf{\ell_0+h_0}{N}\right) \log^{+}\left(\gf{1}{R_{\mathfrak{p}}(Y)}\right) \\ +\gf{\mu-1}{n} \gf{d_{\mathfrak{p}}}{[\K:\Q]} \mathbb{1}_{nN+h_0+\ell_0 \geqslant p(\mathfrak{p})} \log(nN+\ell_0+h_0), \end{multline} où $d_{\mathfrak{p}}=[\K_{\mathfrak{p}}:\Q_{p(\mathfrak{p})}]$ est le degré de la complétion $\mathfrak{p}$-adique de $\K$ sur le corps $\Q_p$ adapté. Cette inégalité découle essentiellement de l'application du théorème de Christol-Dwork \cite[Theorem 2.1, p. 159]{Dwork}.

On obtient alors en sommant l'inégalité \eqref{eq:Dworkp247} sur tous les premiers $\mathfrak{p}$ tels que $p(\mathfrak{p}) \leqslant n$ et $\mathfrak{p} \not\in \mathcal{S}$
$$\sigma_n(Y)\leqslant \left(1+\gf{\ell_0+h_0}{nN}\right) \rho_n(Y)+\gf{\mu-1}{n} \log(nN+\ell_0+h_0) \sum\limits_{u<p \leqslant n} \left(\sum\limits_{\mathfrak{p} \mid p} \gf{d_{\mathfrak{p}}}{[\K:\Q]}\right) +\sum_{\mathfrak{p} \in \mathcal{S}} \gf{1}{n} \sup\limits_{m \leqslant n} \log^{+} |Y_m|_{\mathfrak{p}},$$ où $u$ est un entier tel que $\forall \mathfrak{p} \in \mathcal{S}, p(\mathfrak{p}) \leqslant u$.

Or, selon la formule d'Hadamard, on a pour tout $\mathfrak{p}$, $\limsup |Y_s|_{\mathfrak{p}}^{1/s}=\gf{1}{R_{\mathfrak{p}}(Y)}$, donc $|Y_m|_{\mathfrak{p}}^{1/m}  \leqslant \gf{1}{R_{\mathfrak{p}}(Y)}+o(1)$, d'où si $m \leqslant n$, $$\log |Y_m|_{\mathfrak{p}} \leqslant m\log\left(\gf{1}{R_{\mathfrak{p}}(Y)}+o(1)\right) \leqslant n\log\left(\gf{1}{R_{\mathfrak{p}}(Y)}\right)+n \nu,$$ où $\nu >0$ est une constante. Ainsi, $$\gf{1}{n} \sup\limits_{m \leqslant n} \log^{+} |Y_m|_{\mathfrak{p}} \leqslant \max\left(0,n \log\left(\gf{1}{R_{\mathfrak{p}}(Y)} \right)+ \nu n\right) \leqslant \log^{+}\left(\gf{1}{R_{\mathfrak{p}}(Y)}\right) + \nu n.$$
Donc $$\sum_{\mathfrak{p} \in \mathcal{S}} \gf{1}{n} \sup\limits_{m \leqslant n} \log^{+} |Y_m|_{\mathfrak{p}} \leqslant \mathrm{Card}(\mathcal{S}) \nu + \sum\limits_{\mathfrak{p} \in \mathcal{S}} \log^{+}\left(\gf{1}{R_{\mathfrak{p}}(Y)}\right) := \gamma$$ et $\gamma$ est une constante ne dépendant que de $Y$.

De plus, on a, pour tout nombre premier $p$, $\sum\limits_{\mathfrak{p} \mid p} d_{\mathfrak{p}}=[\K:\Q]$ (cf \cite[p. 222]{Dwork}), d'où, en notant $\pi$ la fonction de comptage des nombres premiers, \begin{align*}\sigma_n(Y) &\leqslant \left(1+\gf{\ell_0+h_0}{nN}\right) \rho_n(Y)+\gf{\mu-1}{n} \log(nN+\ell_0+h_0)\left(\pi(n)-\pi(u)\right) +\gamma \\
 &\leqslant (1+o(1)) \rho_n(Y)+\gf{\mu-1}{n} \log(n)(1+o(1))\gf{n}{\log n} (1+o(1)) +\gamma\end{align*} selon le théorème des nombres premiers. Finalement, $$\sigma_n(Y) \leqslant (1+o(1)) \rho_n(Y)+(\mu-1)(1+o(1))+\gamma,$$ de sorte que si $\rho_n(Y)=o(\log n)$, alors $\sigma_n(Y)=o(\log n)$ également. 
 
\bigskip

Donc selon les étapes 1 et 2 précédentes, si $G$ vérifie la condition de Bombieri \emph{au sens large}, alors $\sigma_n(Y)=o(\log n)$. Or, en reprenant la preuve de la proposition \ref{prop:reformulationdefGfon}, on voit que $\sigma_n(Y)=o(\log n)$ si et seulement si tous les coefficients de $Y$ vérifient la condition \textbf{c)} de la définition \ref{def:gfonctionlarge}. Ainsi, les coefficients de la matrice $Y$ sont des $G$-fonctions \emph{au sens large}.

\bigskip

\textbf{Étape 3}. Soit $L \in \K(z)[\mathrm{d}/\mathrm{d}z]$ un $G$-opérateur \emph{au sens large} d'ordre $\mu$. Selon le théorème \ref{th:acklarge} et le théorème de Frobenius \cite[Theorem 3.5.2, p. 349]{Hille}, on peut trouver une base de solutions de l'équation $L(y(z))=0$ au voisinage de $0$ de la forme $$\widetilde{F}=(\widetilde{f}_1(z),\dots,\widetilde{f}_{\mu}(z))=(f_1(z), \dots, f_{\mu}(z))z^{C}=F z^C,$$ où $C \in \mathcal{M}_{\mu}(\K)$ est à valeurs propres rationnelles. 

Donc si $G=A_L$ est la matrice compagnon de $L$, et $U(z)$ est la matrice wronskienne de la famille $(g_1, \dots, g_{\mu})$, on a $U={}^t \left(\widetilde{F},\widetilde{F}', \dots, \widetilde{F}^{(\mu-1)}\right)$ et $U$ est inversible, car $\widetilde{F}$ est une base de solutions de $L(y(z))=0$. De plus, $U$ vérifie $U'=A_L U$. Or, selon la formule de Leibniz, pour tout $s \leqslant \mu-1$, \begin{align*}
\widetilde{F}^{(s)}&=\sum_{k=0}^s \dbinom{s}{k} F^{(k)} (z^C)^{(s-k)}=\sum_{k=0}^s \dbinom{s}{k} F^{(k)} C(C-I_n)\dots (C-(s-k-1)I_n) z^{C-(s-k)I_n}\\
&=\left(\sum_{k=0}^s \dbinom{s}{k} z^{\mu-1-s+k} F^{(k)} C(C-I_n)\dots (C-(s-k-1)I_n) \right)z^{C-(\mu-1)I_n} := H_s z^{C-(\mu-1)I_n}
\end{align*}
donc la matrice $U$ s'écrit sous la forme $U(z)=H(z) z^{\widetilde{C}}$, où $\widetilde{C}=C-(\mu-1)I_n$ est à valeurs propres rationnelles et $H \in \mathcal{M}_{\mu}(\K\llbracket z \rrbracket)$ est la matrice de $s$\up{ème} ligne $H_{s-1}$. Comme $U$ est inversible, $H=U z^{-\widetilde{C}}$ l'est également. 

Donc selon la conclusion de l'étape 2, les coefficients de $H$ sont des $G$-fonctions \emph{au sens large}. En effet, la proposition \ref{prop:andrebombierilarge} nous assure que $G$ vérifie la condition de Bombieri \emph{au sens large}, puisque $L$ est un $G$-opérateur au sens large.

De plus, la première ligne de $H$ est $H_0=z^{\mu-1} F$, si bien que les $f_i(z)$ sont des $G$-fonctions \emph{au sens large}.

\medskip

\textbf{Étape 4 : changement de variable.} Soit $L \in \K(z)[\mathrm{d}/\mathrm{d}z]$ un $G$-opérateur \emph{au sens large} d'ordre $\mu$ et $\alpha \in \Pro^1(\Qbar)$. Selon la proposition \ref{prop:galochkinLargeInvChgeVar} \textbf{b)}, l'opérateur obtenu par changement de variable $u=z-\alpha$ (avec la convention que $z-\alpha=z^{-1}$ si $\alpha=\infty$) est un $G$-opérateur \emph{au sens large}. Donc selon l'étape 3, il existe une base de solutions de $L_{\alpha}(y(u)=0$ au voisinage de $0$ de la forme $\left(f_1(u), \dots, f_{\mu}(u)\right) u^{C}$, où les $f_i(u) \in \K\llbracket u \rrbracket$ sont des $G$-fonctions \emph{au sens large} et $C \in \mathcal{M}_{\mu}(\K)$ est à valeurs propres rationnelles. Par changement de variable $z=u+\alpha$, on obtient donc pour base de solutions de $L(y(z))=0$ au voisinage de $\alpha$ la famille $$
\left(f_1(z-\alpha), \dots, f_{\mu}(z-\alpha)\right) (z-\alpha)^C,$$ qui est bien de la forme annoncée dans l'énoncé du théorème \ref{th:complementACKlarge}.

\section{Taille \emph{au sens large} d'un opérateur différentiel} \label{sec:ch1senslarge}

Le but de cette partie est d'adapter les résultats de \cite{LepetitSize} aux $G$-opérateurs \emph{au sens large}, ce qui permettra de prouver des propriétés algébriques des $G$-opérateurs \emph{au sens large} cruciales pour les applications aux $E$-opérateurs \emph{au sens large} définis et étudiés dans la section \ref{sec:csqACK}.

\subsection{Taille \emph{au sens large} d'un module différentiel}\label{subsec:taillelarge} 

Commençons par définir une notion de taille \emph{au sens large} pour les modules différentiels. Pour cela, il convient d'étudier le comportement de la suite $(\sigma_s(G))_{s \in \N}$ par équivalence de systèmes différentiels.

\begin{lem}\label{lem:4sizeinvariantequivsystdiff}
Soit $\K$ un corps de nombres, $G \in M_n(\K(z))$ et $P \in \GL_n(\K(z))$, soit $H=P[G]=PGP^{-1}+P'P^{-1}$ une matrice définissant un système $y'=Hy$ équivalent à $y'=Gy$ sur $\Qbar(z)$. Alors $\sigma_s(H)=\sigma_s(G)+o(1)$, de sorte que $H$ satisfait la condition de Galochkin \emph{au sens large} si et seulement si $G$ la satisfait.
\end{lem}

\begin{dem}
    \item Selon \cite[p. 7]{LepetitSize}, on a $$\max_{0 \leqslant m \leqslant s} \left|\gf{H_m}{m!}\right|_{\mathfrak{p},\mathrm{Gauss}} \leqslant C(\mathfrak{p}) \max_{0 \leqslant m \leqslant s} \left|\gf{G_m}{m!}\right|_{\mathfrak{p}, \mathrm{Gauss}}.$$ où $C(\mathfrak{p})=0$ pour tout premier sauf un nombre fini.
    Donc $$\sigma_s(H) \leqslant \gf{1}{s} \sum_{\mathfrak{p} \in \Spec(\Oal_{\K})} \log^{+} C(\mathfrak{p}) + \sigma_s(G)=\sigma_s(G)+o(1).$$
    
    Symétriquement, puisque $G=P^{-1}[H]$, on a $\sigma_s(G) \leqslant \sigma_s(H)+o(1)$. Ceci donne l'égalité voulue.\end{dem}

Le lemme \ref{lem:4sizeinvariantequivsystdiff} implique que l'on peut définir sans ambiguïté la \emph{taille au sens large} d'un module différentiel $\mathcal{M}$ par la suite $\sigma_s(\mathcal{M}) := \sigma_s(A)$, car tout module différentiel $\mathcal{M}$ peut être associé à une unique classe d'équivalence de systèmes différentiels $[A]$. Précisément, $\left(\sigma_s(\mathcal{M})\right)_{s \in \N^*}$ est une classe d'équivalence dans $\R^{\N^*}$ pour la relation d'équivalence $$(u_n) \sim (v_n) \ssi u_n-v_n=o(1).$$

On dit que $u_n \lesssim v_n$ s'il existe $(\varepsilon_n)$ tendant vers $0$ quand $n$ tend vers l'infini et tel que $u_n \leqslant v_n + \varepsilon_n$. Cette notation sera utile, car toutes les inégalités énoncées ci-dessous sont vraies à un terme négligeable en $o(1)$ près.

On rappelle qu'un opérateur $L$ est un $G$-opérateur \emph{au sens large} si et seulement si $$\sigma_s(L)=o(\log s).$$

Comme, pour tout idéal premier $\mathfrak{p}$ de $\K$ et $G \in \mathcal{M}_n(\K(z))$, la quantité $R_{\mathfrak{p}}(G)$ est invariante par équivalence de systèmes différentiels, on peut associer à un module différentiel $\mathcal{M}$ représentant un système $y'=Gy$ une quantité $\rho_s(\mathcal{M})=\rho_s(G)$(voir \cite[p. 67]{Andre}). On dit que $(\rho_s(G))_{s \in \N}$ est le rayon de convergence global \emph{au sens large} de $\mathcal{M}$.

Le résultat suivant est l'analogue \emph{au sens large} du théorème d'André-Bombieri, voir la démonstration de la proposition \ref{prop:andrebombierilarge}.
C'est une version quantitative de la proposition \ref{prop:andrebombierilarge}.

\begin{prop} \label{prop:ABbroad}
On a, pour $G \in M_{\mu}(\K(z))$, \begin{equation}\label{eq:ABbroad1}\rho_s(G) \lesssim \sigma_s(G) \lesssim \rho_s(G)+\mu-1.\end{equation}

Par conséquent, la condition de Bombieri \emph{au sens large} $\rho_s(G)=o(\log s)$ est vérifiée si et seulement la condition de Galochkin \emph{au sens large} est satisfaite par $G$.
\end{prop}
La preuve découle immédiatement de celle de la proposition \ref{prop:andrebombierilarge}. 

\begin{rqu}
Si $\limsup\limits_{s \rightarrow +\infty} \sigma_s(G) < +\infty$, la matrice $G$ satisfait la condition de Galochkin \emph{au sens strict} et de manière équivalente, par le théorème d'André-Bombieri, la condition de Bombieri \emph{au sens strict} $\rho_s(G) \xrightarrow[s \rightarrow +\infty]{} \rho(G) <+\infty$. Cette situation a déjà été étudiée dans \cite{LepetitSize}, si bien que l'on peut supposer que $G$ ne satisfait pas la condition de Galochkin \emph{au sens strict}. 
\end{rqu}

\medskip

Le résultat suivant montre le comportement de la notion de taille \emph{au sens large} sous les opérations usuelles sur les modules différentiels.

\begin{prop} \label{prop:4tailleetopmodules}
Soient $\mathcal{M}_1, \mathcal{M}_2, \mathcal{M}_3$ des modules différentiels sur $\Qbar(z)$. Alors
\begin{enumerateth}
\item Si $\mathcal{M}_2$ est un sous-module différentiel de $\mathcal{M}_1$, alors $$\sigma_s(\mathcal{M}_2) \lesssim \sigma_s(\mathcal{M}_1) \quad \text{et} \quad \sigma_s\left(\mathcal{M}_1 \middle/\mathcal{M}_2\right) \lesssim \sigma_s(\mathcal{M}_1).$$
\item On a $\sigma_s(\mathcal{M}_1 \times \mathcal{M}_2)=\max\left(\sigma_s(\mathcal{M}_1), \sigma_s(\mathcal{M}_2)\right) \lesssim \sigma_s(\mathcal{M}_1) + \sigma_s(\mathcal{M}_2)$.
\item On a  $\sigma_s(\mathcal{M}_1^*) \lesssim \sigma_s(\mathcal{M}_1)+\mu_1-1$, où $\mu_1=\dim_{\Qbar(z)} \mathcal{M}_1$.
\item Si la suite $0 \rightarrow \mathcal{M}_2 \rightarrow \mathcal{M}_1 \rightarrow \mathcal{M}_3 \rightarrow 0$ est exacte, alors on a \begin{equation} \label{eq:4seqexacteeq1}\sigma_s(\mathcal{M}_1) \lesssim 1+ \sigma_s(\mathcal{M}_2)+ \sigma_s(\mathcal{M}_3)+\sigma_s (\mathcal{M}_3^*).\end{equation}

\item On a $\sigma_s(\Sym^N \mathcal{M}_1)  \lesssim N \sigma_s(\mathcal{M}_1)$.
\end{enumerateth}
\end{prop}

La preuve est l'adaptation de celle donnée dans \cite[pp. 8--11]{LepetitSize} pour la proposition 3.

\begin{dem}[de la proposition \ref{prop:4tailleetopmodules}]

On rappelle que si $A \in M_n(\Qbar(z))$, la suite de matrices $(A_s)_{s \in \N}$ est définie par $A_0=I_n$ et $A_{s+1}=A_sA+A'_s$ pour tout $s$. 

\medskip

\textbf{a)} est une conséquence directe de \cite[Lemma 3 p. 8]{LepetitSize}, et \textbf{b)} suit directement de la définition. 

Prouvons l'inégalité \textbf{c)} avec l'aide de la proposition \ref{prop:ABbroad}.

La proposition \ref{prop:ABbroad} implique que $\rho_s(\mathcal{M}_1)=o(\log s)$ si et seulement si $ \sigma_s(\mathcal{M}_1)=o(\log s)$. 

Maintenant, on a $\rho_s(\mathcal{M}_1^*) =\rho_s(\mathcal{M}_1)$, puisque $R_{\mathfrak{p}}(\mathcal{M}_1^*)=R_{\mathfrak{p}}(\mathcal{M}_1)$ pour tout idéal premier $\mathfrak{p}$ d'un corps de nombres $\K$ \cite[Proposition 1, p. 67]{Andre}.

 Ainsi, l'application de la proposition \ref{prop:ABbroad} au module adjoint de $\mathcal{M}_1^*$ donne $$\sigma_s(\mathcal{M}_1^*) \lesssim \rho_s(\mathcal{M}_1^*)+\mu_1-1=\rho_s(\mathcal{M}_1)+\mu_1-1 \lesssim \sigma_s(\mathcal{M}_1)+\mu_1-1,$$ où $\mu_1=\dim\limits_{\Qbar(z)} \mathcal{M}_1$. Ceci démontre \textbf{b)}.

\bigskip

On passe à la preuve de \textbf{d)}.

Par \cite[Lemma 3, p. 8]{LepetitSize}, on peut trouver des bases adaptées de $\mathcal{M}_1$, $\mathcal{M}_2$, $\mathcal{M}_3$ telles que, dans ces bases, $\mathcal{M}_1$ (resp. $\mathcal{M}_2$, $\mathcal{M}_3$) représente le système $\partial y=Gy$ (resp. $\partial y=G^{(2)} y$ et $\partial y=G^{(3)} y$), où $$G=\begin{pmatrix} G^{(2)} & G^{(0)} \\ 0 & G^{(3)} \end{pmatrix}.$$

Soit $\K$ un corps de nombres tel que $G \in M_n(\K(z))$.
On fixe $\mathfrak{p} \in \Spec(\Oal_{\K})$. On considère $t_{\mathfrak{p}}$ une variable libre sur $\K$ appelée \emph{point générique}. On peut alors construire une extension $\Omega_{\mathfrak{p}}$ complète et algébriquement close de $\left(\K(t_{\mathfrak{p}}), | \cdot |_{\mathfrak{p}, \mathrm{Gauss}}\right)$.

On pose $$X^{(2)}_{t_{\mathfrak{p}}}(z)=\sum\limits_{s=0}^{\infty} \gf{G^{(2)}_s(t_{\mathfrak{p}})}{s!} (z-t_{\mathfrak{p}})^s \in \GL_n(\Omega_{\mathfrak{p}}\llbracket z-t_{\mathfrak{p}}\rrbracket)$$ (resp. $X^{(3)}_{t_{\mathfrak{p}}} \in \GL_m(\Omega_{\mathfrak{p}}\llbracket z-t_{\mathfrak{p}}\rrbracket)$) une matrice fondamentale de solutions de $\partial y=G^{(2)} y$ (resp. $\partial y=G^{(3)}y$) au point générique $t_{\mathfrak{p}}$ telle que $X^{(2)}_{t_{\mathfrak{p}}}(t_{\mathfrak{p}})=I_n$ (resp. $X^{(3)}_{t_{\mathfrak{p}}}(t_{\mathfrak{p}})=I_m$). On considère $X^{(0)}_{t_{\mathfrak{p}}} \in M_{n,m}(\Omega_{\mathfrak{p}}\llbracket z-t_{\mathfrak{p}}\rrbracket)$ une solution de \begin{equation} \label{eq:4EDX0}\partial X_{t_{\mathfrak{p}}}^{(0)}=G^{(3)}X_{t_{\mathfrak{p}}}^{(0)}+G^{(0)} X_{t_{\mathfrak{p}}}^{(2)} \end{equation} telle que $X^{(0)}_{t_{\mathfrak{p}}}(t_{\mathfrak{p}})=0$.

Alors $X_{t_{\mathfrak{p}}}=\begin{pmatrix} X_{t_{\mathfrak{p}}}^{(2)} & X_{t_{\mathfrak{p}}}^{(0)} \\ 0 & X_{t_{\mathfrak{p}}}^{(3)} \end{pmatrix}$ est une matrice fondamentale de solutions de $\partial y=Gy$ telle que $X_{t_{\mathfrak{p}}}(t_{\mathfrak{p}})=I_{n+m}$. Ainsi, on a $\forall s \in \N,\; \partial^s(X_{t_{\mathfrak{p}}})(t_{\mathfrak{p}})=G_s(t_{\mathfrak{p}})$. Puisque l'on sait que $\partial^s(X^{(2)}_{t_{\mathfrak{p}}})(t_{\mathfrak{p}})=G^{(2)}_s(t_{\mathfrak{p}})$ (et de même pour $X_{t_{\mathfrak{p}}}^{(3)}$), il suffit d'estimer la norme de Gauss de $\partial^s(X^{(0)}_{t_{\mathfrak{p}}})(t_{\mathfrak{p}})$ pour obtenir une estimation sur la taille de $G$.

Par souci de simplicité, on omettra l'indice $t_{\mathfrak{p}}$ dans ce qui suit.

Selon \cite[p. 10]{LepetitSize}, on a, pour tout $\ell \in \N^*$ et $N \geqslant \ell$, \begin{multline*}
    \log^{+} \left\|\gf{\partial^{\ell}(X^{(0)})(t_{\mathfrak{p}})}{\ell!}\right\|_{\mathfrak{p}, \mathrm{Gauss}} \lesssim h(N,\mathfrak{p},G^{(3)})+h(N,\mathfrak{p},-{}^t G^{(3)})+ h(N,\mathfrak{p},G^{(2)}) \\ + \log^{+} \| G^{(0)} \|_{\mathfrak{p}, \mathrm{Gauss}}+ \log \max_{0 \leqslant i \leqslant N} \gf{1}{|i|_{p}}.
\end{multline*}

De plus, \begin{align*}\left\|\gf{G_\ell}{\ell!}\right\|_{\mathfrak{p}, \mathrm{Gauss}}&=\left\|\gf{\partial^{\ell}(X)(t_{\mathfrak{p}})}{\ell!}\right\|_{\mathfrak{p}, \mathrm{Gauss}} \\ &=\max\left(\left\|\gf{\partial^{\ell}(X^{(0)})(t_{\mathfrak{p}})}{\ell!}\right\|_{\mathfrak{p}, \mathrm{Gauss}}, \left\|\gf{\partial^{\ell}(X^{(2)})(t_{\mathfrak{p}})}{\ell!}\right\|_{\mathfrak{p}, \mathrm{Gauss}}, \left\|\gf{\partial^{\ell}(X^{(3)})(t_{\mathfrak{p}})}{\ell!}\right\|_{\mathfrak{p}, \mathrm{Gauss}} \right),\end{align*} donc finalement, pour tout $s \in \N^*$, \begin{equation}\label{eq:4result1tailleopmod}\sigma_s(\mathcal{M}_1) \lesssim 1+\sigma_s(\mathcal{M}_3)+\sigma_s(\mathcal{M}_3^*)+\sigma_s(\mathcal{M}_2)+\gf{1}{s} \sum_{\mathfrak{p} \in \Spec(\Oal_{\K})} \log^{+} \left\| G^{(0)} \right\|_{\mathfrak{p}, \mathrm{Gauss}}\end{equation} et le dernier terme de la somme tend vers $0$ quand $s$ tend vers l'infini.

\medskip

Passons à présent à la preuve de \textbf{e)}. On notera $a_1 \dots a_N$ la classe de $a_1 \otimes \dots \otimes a_N \in \mathcal{M}^{\otimes N}$ dans le quotient $\mathrm{Sym}^N(\mathcal{M})$.

Soit une base $\mathcal{B}=(e_1, \dots, e_n)$ de $\mathcal{M}_1$ telle que $\mathcal{M}_1$ représente le système $\partial y=Gy$, $G \in M_{n}(\K(z))$. Alors on peut trouver une matrice $G_{\mathrm{sym}}$ à coefficients dans $\K(z)$ telle que dans la base $\mathcal{B}_{\mathrm{sym}}=\{ e_{i_1} \dots e_{i_N}, 1 \leqslant  i_1 \leqslant \dots \leqslant i_N \leqslant n \}$, $\mathrm{Sym}^N(\mathcal{M}_1)$ représente le système $\partial y=G_{\mathrm{sym}} y$ . Suivant André dans \cite[p. 72]{Andre}, introduisons, pour $\mathfrak{p} \in \Spec(\Oal_{\K})$, $X_{t_{\mathfrak{p}}}$ une matrice fondamentale de solutions de $\partial y=Gy$ au voisinage du point générique $t_{\mathfrak{p}}$ telle que  $X_{t_{\mathfrak{p}}}(t_{\mathfrak{p}})=I_n$. Alors on a que  $\partial \mathrm{Sym}^N(X_{t_{\mathfrak{p}}})=G_{\mathrm{sym}} \mathrm{Sym}^N(X_{t_{\mathfrak{p}}}) $.

\medskip

On omettra l'indice $t_{\mathfrak{p}}$ dans ce qui suit. Rappelons que si $\mathcal{C}=(f_1, \dots, f_n)$ est la base canonique de $\Omega_{\mathfrak{p}}^n$, $X=\mathrm{Mat}_{\mathcal{C}}(u)$, alors $\mathrm{Sym}^N(X)$ est la matrice représentant l'endomorphisme $\mathrm{Sym}^N(u)$ dans la  base $\mathcal{C}_{\mathrm{sym}} :=\{ f_{j_1} \dots f_{j_N}, 1 \leqslant j_1 \leqslant \dots \leqslant j_N \leqslant m \}$, où \begin{align*}
\mathrm{Sym}^N(u)(f_{j_1} \dots f_{j_N}) &= u(f_{j_1}) \dots u(f_{j_N}) = \left(\sum\limits_{i_1=1}^{n} X_{i_1 j_1} f_{i_1}\right) \dots \left(\sum\limits_{i_N=1}^{n} X_{i_N j_N} f_{i_N}\right) \\
&=\sum\limits_{1 \leqslant i_1, \dots, i_N \leqslant n} X_{i_1 j_1} \dots X_{i_N j_N} f_{i_1} \dots f_{i_N} \\
&= \sum_{1 \leqslant i_1 \leqslant \dots \leqslant i_N \leqslant n} \;\; \sum\limits_{\sigma \in \mathfrak{S}(\{ i_1, \dots i_N\})} X_{\sigma(i_1) j_1} \dots X_{\sigma(i_N) j_N} f_{i_1} \dots f_{i_N}.
\end{align*}

Donc les coefficients de $\mathrm{Sym}^N(X)$ sont les $\sum\limits_{\sigma \in \mathfrak{S}(\{ i_1, \dots i_N\})} X_{\sigma(i_1) j_1} \dots X_{\sigma(i_N) j_N}$. On constate en particulier, en regardant les coefficients pour $(i_1, \dots, i_N)=(j_1, \dots, j_N)$, que $\Sym^N(X)(t_{\mathfrak{p}})$ est la matrice identité, de sorte que pour tout $s$, $\left(\partial^s \mathrm{Sym}^N(X)\right)(t_{\mathfrak{p}})=(G_{\mathrm{sym}})_{s}(t_{\mathfrak{p}})$.

Comme $|\cdot |_{\mathfrak{p},\mathrm{Gauss}}$ est non archimédienne, il s'agit donc de majorer, à $(j_1, \dots, j_N)$ fixé, $|X_{i_1 j_1} \dots X_{i_N j_N}|_{\mathfrak{p},\mathrm{Gauss}}$ pour tout $(i_1, \dots, i_N)$.

Soit $s \in \N^*$ et $m \leqslant s$. Selon la formule de Leibniz généralisée, on a $$\gf{1}{m!}\partial^m (X_{i_1 j_1} \dots X_{i_N j_N})=\sum_{k_1+\dots+k_N=m} \gf{\partial^{k_1} X_{i_1 j_1}}{k_1!} \dots \gf{\partial^{k_N} X_{i_N j_N}}{k_N!}.$$
On a donc \begin{align}\log^+ \left\vert\gf{1}{m!}\partial^m (X_{i_1 j_1} \dots X_{i_N j_N})(t_{\mathfrak{p}})\right\vert_{\mathfrak{p},\mathrm{Gauss}} &\leqslant \log^+ \max\limits_{k_1 \leqslant m} \left\vert\gf{\partial^{k_1} X_{i_1 j_1}}{k_1!}(t_{\mathfrak{p}})\right\vert_{\mathfrak{p},\mathrm{Gauss}} + \dots \notag\\
& \qquad\qquad\qquad\qquad+ \log^+ \max\limits_{k_N \leqslant m} \left\vert \gf{\partial^{k_N} X_{i_N j_N}}{k_N!}(t_{\mathfrak{p}})\right\vert_{\mathfrak{p},\mathrm{Gauss}} \label{eq:majorationXsymGauss} \\ & \leqslant N \log^+ \max_{k \leqslant s} \left\| \gf{\partial^{k} X}{k!}(t_{\mathfrak{p}})\right\|=N \log^+ \max_{k \leqslant s} \left\| \gf{G_k}{k!} \right\|_{\mathfrak{p},\mathrm{Gauss}}.\notag \end{align}
Ainsi, si $c$ est un coefficient de la matrice $\mathrm{Sym}^N X$, et $m \leqslant s$ pour tout $m \leqslant s$, on a  $$\log^+ \left\vert\gf{\partial^m c}{m!}(t_{\mathfrak{p}})\right\vert_{\mathfrak{p},\mathrm{Gauss}}  \leqslant N \log^+ \max_{k \leqslant s} \left\| \gf{G_k}{k!} \right\|_{\mathfrak{p},\mathrm{Gauss}}.$$ En prenant le maximum sur tous les coefficients de $\mathrm{Sym}^N X$ et sur tous les $m \leqslant s$, on obtient donc $$\log^+ \max_{m \leqslant s} \left\|\gf{\partial^m \mathrm{Sym}^N X}{m!}(t_{\mathfrak{p}})\right\|_{\mathfrak{p},\mathrm{Gauss}}  \leqslant N \log^+ \max_{k \leqslant s} \left\| \gf{G_k}{k!} \right\|_{\mathfrak{p},\mathrm{Gauss}},$$ c'est-à-dire $h(s,\mathfrak{p},G_{\mathrm{sym}}) \leqslant N h(s,\mathfrak{p},G)$. Donc en sommant sur tous les $\mathfrak{p} \in \Spec(\Oal_{\K})$, on a $\sigma_s(G_{\mathrm{sym}}) \leqslant N\sigma_s(G)$. D'où finalement $$\sigma_s(\mathrm{Sym}^N \mathcal{M}_1) \lesssim N \sigma_s(\mathcal{M}_1).$$
Ceci conclut la preuve de la proposition \ref{prop:4tailleetopmodules}.
\end{dem}

\begin{rqu}
 Dans le cas strict, on obtient une inégalité $\sigma(\mathrm{Sym}^N \mathcal{M}_1) \leqslant (1+\log(N)) \sigma(\mathcal{M}_1)$.
 
 En effet, dans l'équation \eqref{eq:majorationXsymGauss}, on peut supposer quitte à réordonner les termes du produit que $k_1 \geqslant k_2 \geqslant \dots \geqslant k_N$, de sorte que, pour tout $\ell$, $m \geqslant k_1+ \dots +k_{\ell} \geqslant \ell k_{\ell} $ et donc $k_{\ell} \leqslant \gf{m}{\ell}$.
  
On obtient ainsi $$\log^+ \left\vert\gf{1}{m!}\partial^m (X_{i_1 j_1} \dots X_{i_N j_N})(t_{\mathfrak{p}})\right\vert_{\mathfrak{p},\mathrm{Gauss}} \leqslant  \log^+ \max_{k \leqslant s} \left\| \gf{G_k}{k!}(t_{\mathfrak{p}})\right\|+\dots+\log^+ \max_{k \leqslant \lfloor s/N \rfloor} \left\| \gf{G_k}{k!}(t_{\mathfrak{p}})\right\|,$$ de sorte que $$\gf{1}{s} \log^+ \max_{m \leqslant s} \left\|\gf{(G_{\mathrm{sym}})_m}{m!}\right\|_{\mathfrak{p},\mathrm{Gauss}}  \leqslant  \gf{1}{s}\log^+ \max_{k \leqslant s} \left\| \gf{G_k}{k!}(t_{\mathfrak{p}})\right\|+\dots+\gf{1}{s}\log^+ \max_{k \leqslant \lfloor s/N \rfloor} \left\| \gf{G_k}{k!}(t_{\mathfrak{p}})\right\|.$$ En d'autres termes, $$h(s,\mathfrak{p},G_{\mathrm{sym}}) \leqslant h(s,\mathfrak{p},G)+\left\lfloor\gf{s}{2}\right\rfloor \gf{1}{s} h\left(\left\lfloor\gf{s}{2}\right\rfloor,\mathfrak{p},G\right)+\dots+ \left\lfloor\gf{s}{N}\right\rfloor \gf{1}{s} h\left(\left\lfloor\gf{s}{N}\right\rfloor,\mathfrak{p},G\right).$$
En sommant sur tous les $\mathfrak{p} \in \Spec(\Oal_{\K})$, on obtient
\begin{equation} \label{eq:conclusionaltE}
    \sigma_s(G_{\mathrm{sym}}) \leqslant \sigma_s(G)+\left\lfloor\gf{s}{2}\right\rfloor \gf{1}{s} \sigma_{\lfloor s/2 \rfloor}(G)+\dots+\left\lfloor\gf{s}{N}\right\rfloor \gf{1}{s} \sigma_{\lfloor s/N \rfloor}(G),
\end{equation}
ce qui est une inégalité alternative au point \textbf{e)} de la proposition \ref{prop:4tailleetopmodules}, mais est peu exploitable.

En revanche, en passant à la limite supérieure $s \rightarrow +\infty$ dans \eqref{eq:conclusionaltE}, on obtient 
\begin{equation}\sigma(G_{\mathrm{sym}}) \leqslant \sigma(G)+\gf{1}{2} \sigma(G) + \dots + \gf{1}{N} \sigma(G) \leqslant (1+\log(N))\sigma(G), \end{equation}
 ce qui est le résultat annoncé pour la taille \emph{au sens strict} dans \cite[Lemma 2 \textbf{c)}, p. 72]{Andre}.
\end{rqu}

\subsection{Propriétés algébriques des $G$-opérateurs \emph{au sens large}} 

Dans ce qui suit, on note $\partial$ la dérivation standard $\mathrm{d}/\mathrm{d}z$ sur $\Qbar(z)$.

Nous allons utiliser les résultats de la section précédente pour en déduire des propriétés de stabilité algébriques sur l'ensemble des $G$-opérateurs \emph{au sens large}, qui ont été énoncées dans la partie \ref{subsec:galochkinlarge} Toutes les propriétés listées ci-dessous sont également vraies pour les $G$-opérateurs \emph{au sens strict}.

\begin{prop}\label{prop:algpropGop}
Soit $L$ un $G$-opérateur \emph{au sens large}. Alors : 
\begin{enumerateth}[label=\textbf{\alph*)}]
  \item Tout diviseur à droite de $L$ dans $\Qbar(z)[\partial]$ est un $G$-opérateur \emph{au sens large}.
  \item L'opérateur adjoint $L^*$ de $L$ est un $G$-opérateur \emph{au sens large}.
  \item L'ensemble des $G$-opérateurs \emph{au sens large} satisfait la \emph{propriété de Ore à gauche} : si $L$ et $M$ sont des $G$-opérateurs  \emph{au sens large}, il existe un multiple commun à gauche de $L$ et $M$ qui est un $G$-opérateur \emph{au sens large}. 
\end{enumerateth}
\end{prop}

\begin{dem}
Les assertions \textbf{a)} et \textbf{b)} découlent directement des points \textbf{a)} et \textbf{c)} de la proposition  \ref{prop:4tailleetopmodules}, puisque, d'une part, si $N$ est un diviseur à droite de $L$, le module différentiel $\mathcal{M}_N$ est un sous-module différentiel de $\mathcal{M}_L$ et d'autre part, $\mathcal{M}_{L^*} \simeq \left(\mathcal{M}_L\right)^*$.

Il reste à prouver \textbf{c)}. Soit $N=\mathrm{LCLM}(L,M)$ le générateur de l'idéal à gauche $\Qbar(z)[\partial]L \cap \Qbar(z)[\partial]N$. Alors $N$ est un multiple commun à gauche de $L$ et $M$ et il y a un morphisme injectif de modules différentiels $$\fonction{\phi}{\Qbar(z)[\partial]/\Qbar(z)[\partial]N}{\Qbar(z)[\partial]/\Qbar(z)[\partial]L \times \Qbar(z)[\partial]/\Qbar(z)[\partial]M}{P \mod N}{(P \mod L, P \mod M)},$$ qui fait de $\Qbar(z)[\partial]/\Qbar(z)[\partial]N$ un sous-module différentiel du module produit $\Qbar(z)[\partial]/\Qbar(z)[\partial]L \times \Qbar(z)[\partial]/\Qbar(z)[\partial]M $. De plus, $\Qbar(z)[\partial]/\Qbar(z)[\partial]N\simeq \mathcal{M}_{N^*}$ et de même pour $L$ et $M$.

D'où, par la proposition \ref{prop:4tailleetopmodules} \textbf{a)} et \textbf{b)}, on a $\sigma_s(N^*) \lesssim \max(\sigma_s(L^*), \sigma_s(M^*))$.

Le point \textbf{b)} nous assure alors, puisque $L$ et $M$ sont des $G$-opérateurs \emph{au sens large} que $N$ est un $G$-opérateur \emph{au sens large}. 
\end{dem}

\begin{coro} \label{prop:4andreproduitsommegfonctions}
Si $y_1$ (resp. $y_2$) est une solution d'un $G$-opérateur \emph{au sens large} $L_1$ (resp. $L_2$), alors $y_1+y_2$ est solution d'un $G$-opérateur \emph{au sens large} $L_3$;
\end{coro}

\begin{dem}
En utilisant la proposition \ref{prop:algpropGop} \textbf{c)}, on considère $M=M_1 L_1=M_2 L_2$ un multiple commun à gauche de $L_1$ et $L_2$ dans $\Qbar(z)[\partial]$ qui est un $G$-opérateur \emph{au sens large}. Alors on a $M(y_1)=M(y_2)=0$, de sorte que $M(y_1+y_2)=0$.
\end{dem}

La proposition \ref{prop:4tailleetopmodules} implique aussi le résultat suivant :
\bigskip
\bigskip
\begin{prop} \label{prop:tailleproduitgop}
Pour tout $L_1, L_2 \in \Qbar(z)[\partial]$, on a les inégalités suivantes : \begin{align*} \max\big(\sigma_s(L_1), \sigma_s(L_2)\big) \lesssim \sigma_s(L_1 L_2) &\lesssim \sigma_s(L_1)+2\sigma_s(L_2)+\ord(L_2)-1.\end{align*}

Ainsi, un produit de $G$-opérateurs \emph{au sens large} est un $G$-opérateur \emph{au sens large}.
\end{prop}

\begin{dem}[de la proposition \ref{prop:tailleproduitgop}]
Selon \cite[p. 47]{Singer}, la suite $$ 0 \rightarrow \mathcal{M}_{L_2} \rightarrow \mathcal{M}_{L_1 L_2} \rightarrow \mathcal{M}_{L_1} \rightarrow 0$$ est exacte. Par conséquent, $\mathcal{M}_{L_2}$ est un sous-module différentiel de $\mathcal{M}_{L_1 L_2}$ et on a $\mathcal{M}_{L_1} \simeq \left.\mathcal{M}_{L_1 L_2} \middle/ \mathcal{M}_{L_2}\right.$. Il découle alors de la proposition \ref{prop:4tailleetopmodules} \textbf{a)} que $$\max(\sigma_s(L_1),\sigma_s(L_2))=\max\left(\sigma_s(\mathcal{M}_{L_1}),\sigma_s(\mathcal{M}_{L_2})\right) \lesssim \sigma_s(\mathcal{M}_{L_1L_2})=\sigma_s(L_1 L_2),$$ ce qui est la première inégalité que l'on voulait prouver.

D'autre part, la  proposition \ref{prop:4tailleetopmodules} \textbf{d)} donne $\sigma_s(L_1 L_2) \lesssim \sigma_s(L_1)+ \sigma_s(L_2)+\sigma_s(L_2^*)$. En utilisant finalement la proposition \ref{prop:4tailleetopmodules} \textbf{b)}, on obtient $$\sigma_s(L_1 L_2) \lesssim \sigma_s(L_1)+ 2\sigma_s(L_2)+\ord(L_2)-1,$$ ce qui est le résultat voulu.
\end{dem}

\subsection{Produit de solutions d'un $G$-opérateur \emph{au sens large}}

Le résultat suivant est une conséquence des points \textbf{a)} et \textbf{e)} du théorème \ref{prop:4tailleetopmodules}. La preuve est, \emph{mutatis mutandis}, la même que celle de \cite[Proposition 4, p. 11]{LepetitSize}.

\begin{prop}\label{prop:tailleetproduittensoriellarge}
Soit $\left(\mathcal{M}_i\right)_{1 \leqslant i \leqslant N}$ une famille de modules différentiels sur $\Qbar(z)$. Alors $$\sigma_s(\mathcal{M}_1 \otimes_{\Qbar(z)} \dots \otimes_{\Qbar(z)} \mathcal{M}_N) \lesssim N\max\big(\sigma_s(\mathcal{M}_1), \dots, \sigma(\mathcal{M}_N)\big).$$
En particulier, on a $\sigma_s\left(\mathcal{M}_1^{\otimes N}\right) \lesssim N \sigma_s(\mathcal{M}_1)$.
\end{prop} 

\begin{rqu}
 Dans le cas strict, on a (voir la remarque finale de la sous-section \ref{subsec:taillelarge}) : $$\sigma(\mathcal{M}_1 \otimes_{\Qbar(z)} \dots \otimes_{\Qbar(z)} \mathcal{M}_N) \leqslant (1+\log(N))\max\big(\sigma_s(\mathcal{M}_1), \dots, \sigma(\mathcal{M}_N)\big).$$ 
\end{rqu}

\begin{prop} \label{prop:andreproduitsommegfonctionslarge}
Soient $y_1$ et $y_2$ des solutions respectives de $G$-opérateurs \emph{au sens large} $L_1$ et $L_2$. On suppose que chaque $L_i$ est d'ordre minimal pour $y_i$. Alors
$y_1 y_2$ est solution d'un $G$-opérateur \emph{au sens large} $L_4$ et on peut trouver un tel $L_4$ vérifiant $$\sigma_s(L_4) \lesssim 2 \max\big(\sigma_s(L_1),\sigma_s(L_2)\big).$$
\end{prop}

\begin{rqu}
En utilisant la proposition \ref{prop:tailleetproduittensoriellarge} \textbf{e)}, on peut généraliser la proposition \ref{prop:andreproduitsommegfonctionslarge} à un produit quelconque $y_1 \dots y_N$ de solutions de $G$-opérateurs : si $L_0$ est l'opérateur minimal de $y_1 \dots y_N$ sur $\Qbar(z)$ et $L_i$ est l'opérateur minimal de $y_i$, on obtient $$\sigma_s(L_0) \lesssim N \max\big(\sigma_s(L_1),\dots,\sigma_s(L_N)\big).$$
\end{rqu}

\begin{dem}[de la proposition \ref{prop:andreproduitsommegfonctionslarge}]
Si $L \in \Qbar(z)[\partial]$ est l'opérateur minimal d'une fonction $y$, on remarque qu'il y a un isomorphisme naturel de modules différentiels entre le quotient $\Qbar(z)[\partial]/\Qbar(z)[\partial]L$ et $$\Qbar(z)[\partial](y) :=\big\{M(y) \mid M \in \Qbar(z)[\partial]\big\}$$ donné par $M \mod L \mapsto M(y)$. D'où $\sigma(L)=\sigma\left(\left(\Qbar(z)[\partial](y)\right)^*\right)$.

On définit aussi, pour $u,v$ solutions de $G$-opérateurs, $\Qbar(z)[\partial](u, v):= \left(\Qbar(z)[\partial](u)\right)[\partial](v)$ qui est l'ensemble des combinaisons linéaires à coefficients dans $\Qbar(z)$ des $\partial^k(u) \partial^{\ell}(v)$, quand $k, \ell \in \N$. L'existence d'équations différentielles à coefficients dans $\Qbar(z)$ satisfaites par $u$ et $v$ nous assure qu'il s'agit en effet d'un $\Qbar(z)$-espace vectoriel de dimension finie, de sorte qu'il peut être muni d'une structure de module différentiel.

On définit le morphisme $$\fonction{\psi}{\Qbar(z)[\partial](y_1) \otimes \Qbar(z)[\partial](y_2)}{\Qbar(z)[\partial](y_1, y_2)}{K(y_1) \otimes L(y_2)}{K(y_1)L(y_2).}$$ De plus, l'ensemble $\mathcal{M}=\Qbar(z)[\partial](y_1 y_2)$ est un sous-module différentiel de $\Qbar(z)[\partial](y_1, y_2)$, car si $L \in \Qbar(z)[\partial]$, $L=\sum\limits_{k=0}^{\mu} a_k \partial^k$, on a $$L(y_1 y_2)=\sum_{k=0}^{\mu} a_k \sum_{p=0}^k \binom{k}{p} \partial^p(y_1) \partial^{k-p}(y_2)=\psi\left(\sum_{k=0}^{\mu} a_k \sum_{p=0}^k \binom{k}{p} \partial^p(y_1) \otimes \partial^{k-p}(y_2)\right) \in \Qbar(z)[\partial](y_1, y_2).$$

Par définition, la restriction de $\psi$ à $\psi^{-1}(\mathcal{M})$ est une surjection $\widetilde{\psi} : \psi^{-1}(\mathcal{M}) \rightarrow \mathcal{M}$.

La factorisation de $\widetilde{\psi}$ par son noyau fournit un isomorphisme entre $\mathcal{M}$ et $\mathcal{N}/\ker(\widetilde{\psi})$, où $\mathcal{N}=\psi^{-1}(\mathcal{M})$ est un sous-module différentiel de $\Qbar(z)[\partial](y_1) \otimes \Qbar(z)[\partial](y_2)$. 

 En passant aux modules duaux, on obtient $\mathcal{M}^* \simeq \mathcal{N}^*/\Img(\widetilde{\psi}^*)$ donc par la proposition \ref{prop:4tailleetopmodules} \textbf{a)}, si $L_4$ est l'opérateur minimal non nul de $y_1 y_2$ sur $\Qbar(z)$, \begin{equation} \label{eq:4effectifandre3}\sigma_s(L_4) =\sigma_s(\mathcal{M}^*) \lesssim \sigma_s(\mathcal{N}^*).\end{equation}

Par ailleurs, le morphisme dual de $\mathcal{N} \hookrightarrow \Qbar(z)[\partial](y_1) \otimes \Qbar(z)[\partial](y_2)$ est une surjection $$\left(\Qbar(z)[\partial](y_1) \otimes \Qbar(z)[\partial](y_2)\right)^* \twoheadrightarrow \mathcal{N}^*$$ et on a $$\left(\Qbar(z)[\partial](y_1) \otimes \Qbar(z)[\partial](y_2)\right)^* \simeq \left(\Qbar(z)[\partial](y_1)\right)^* \otimes \left(\Qbar(z)[\partial](y_2)\right)^*\;,$$ de sorte que, par les propositions \ref{prop:4tailleetopmodules} \textbf{a)} et \ref{prop:tailleetproduittensoriellarge},
\begin{equation} \label{eq:4effectifandre4}\sigma_s(\mathcal{N}^*) \lesssim \sigma\left(\left(\Qbar(z)[\partial](y_1)\right)^* \otimes \left(\Qbar(z)[\partial](y_2)\right)^*\right) \lesssim 2 \max\big(\sigma_s(L_1),\sigma_s(L_2)\big).\end{equation}

Ainsi, la combinaison de \eqref{eq:4effectifandre3} et \eqref{eq:4effectifandre4} montre que $\sigma_s(L_4) \lesssim 2 \max\big(\sigma_s(L_1),\sigma_s(L_2)\big)$, si bien que $L_4$ est un $G$-opérateur \emph{au sens large}.
\end{dem}

Une autre application est que l'opérateur minimal d'une série Nilsson-Gevrey \emph{au sens large} de type arithmétique et holonome d'ordre $0$  est un $G$-opérateur au sens large. On définit la notion de série Nilsson-Gevrey \emph{au sens large} de type arithmétique d'après André, qui a étudié les propriétés de leurs analogues \emph{au sens strict} dans \cite{AndregevreyI}.

\begin{defi} \label{def:NGAlarge}
 \begin{itemizeth}[label=\textbullet]
    \item Soit $s \in \Q$. Les \emph{séries Nilsson-Gevrey \emph{au sens large} de type arithmétique d'ordre $s$} sont les \begin{equation}\label{eq:defNGholo}y(z)=\sum_{(\alpha, k, \ell) \in S} c_{\alpha,k,\ell} z^{\alpha} \left(\log z\right)^k y_{\alpha,k,\ell}(z),\end{equation} où $S \subset \Q \times \N^2$ est fini, $c_{\alpha,k,\ell} \in \C$ et $$y_{\alpha,k,\ell}(z)=\sum\limits_{n=0}^{\infty} n!^s a_{\alpha,k,\ell,n} z^n \in \Qbar\llbracket z \rrbracket$$  est tel que $(a_{\alpha,k,\ell,n})_{n \in \N}$ vérifie les conditions \textbf{b)} et \textbf{c)} de la définition \ref{def:gfonctionlarge}. On note $\mathrm{NGA}^{\ell}\{z\}_s$ l'ensemble de ces séries.

\item On dit de plus que $y(z)$ est \emph{de type holonome}
si les $y_{\alpha,k,\ell}(z)$ satisfont la condition \textbf{a)} de la définition \ref{def:gfonctionlarge}. On note $\mathrm{NGA}_h^{\ell}\{z\}_s$ l'ensemble des séries Nilsson-Gevrey \emph{au sens large} de type arithmétique et holonome d'ordre $s$. 
\end{itemizeth}
\end{defi}
\begin{rqu}
Dans la remarque p. 717 de \cite{AndregevreyI}, André semble indiquer que toutes les séries Nils\-son-Gevrey (\emph{au sens strict}) holonomes d'ordre $s$ de type arithmétique sont des séries de type holonome, mais nous ne voyons pas comment le prouver avec son raisonnement. En effet, il y est énoncé que si $y(z)$ est est de la forme \eqref{eq:defNGholo}, avec les $\alpha$ dans des classes distinctes modulo $\Z$ et $(a_{\alpha,k,\ell,n})_{n \in \N}$ vérifiant les conditions \textbf{b)} et \textbf{c)} de la définition \ref{def:gfonctionstrict}, alors les $y_{\alpha,k,\ell}(z)$ sont toutes holonomes, ce qui est faux comme le montre l'exemple de $0=y_0(z)-y_0(z)$, avec $y_0(z)$ non holonome.
\end{rqu}

Une conséquence immédiate de la proposition \ref{prop:andreproduitsommegfonctionslarge} et du corollaire \ref{prop:4andreproduitsommegfonctions} est la proposition suivante :

\begin{prop}\label{prop:4chudnovskynilssongevreyandreI}
Soit $S \subset \Q \times \N \times \N$ un ensemble fini, $(c_{\alpha, k,\ell})_{(\alpha, k,\ell) \in S} \in \left(\C^*\right)^{S}$ et une famille $(f_{\alpha, k,\ell}(z))_{(\alpha, k,\ell) \in S}$ de $G$-fonctions \emph{au sens large} non nulles. On considère $$f(z)=\sum\limits_{(\alpha, k, \ell) \in S}^{} c_{\alpha, k,\ell} z^{\alpha} \log(z)^k f_{\alpha, k,\ell}(z),$$ qui est une série Nilsson-Gevrey \emph{au sens large} de type arithmétique et de type holonome d'ordre $0$. Alors $f(z)$ est solution d'une équation différentielle linéaire non nulle à coefficients dans $\Qbar(z)$  et l'opérateur minimal $L$ de $f(z)$ sur $\Qbar(z)$ est un $G$-opérateur \emph{au sens large}.
\end{prop}

\begin{dem}
La fonction $f(z)$ s'écrit comme une somme de produits de solutions de $G$-opérateurs puisque toute $G$-fonction \emph{au sens large} est solution d'un $G$-opérateur au sens large. Donc selon la proposition \ref{prop:andreproduitsommegfonctionslarge} et le corollaire \ref{prop:4andreproduitsommegfonctions}, $f(z)$ est solution d'un $G$-opérateur au sens large $L$. Le point \textbf{a)} de la proposition \ref{prop:algpropGop} nous assure alors que l'opérateur minimal de $f(z)$ sur $\Qbar(z)$ est un $G$-opérateur \emph{au sens large}.
\end{dem}

\section[Conséquences du théorème 1]{Conséquences du théorème \ref{th:acklarge}} \label{sec:csqACK}

\subsection{$E$-opérateurs \emph{au sens large}} \label{subsec:Eoplarges}


Le but de cette partie est de définir la notion d'$E$-opérateur \emph{au sens large}, selon la terminologie d'André et de donner quelques propriétés de ces opérateurs.

Si $f(z)$ est une fonction suffisamment régulière, sa \emph{transformée de Laplace} est la fonction \begin{equation}\label{eq:defLaplace}\mathcal{F}(f) : x \mapsto \displaystyle\int_{0}^{\infty} e^{-xz} f(z) \mathrm{d}z. 
\end{equation} 

Dans le cas où $f(z)=\sum\limits_{n=0}^{\infty} \gf{a_n}{n!} z^n$ est une $E$-fonction \emph{au sens large}, $\mathcal{F}(f)$ est bien définie, car, comme mentionné dans l'introduction, la suite $(|a_n|)_n$ est en réalité majorée par une suite géométrique $(C^{n+1})_n$ à partir d'un certain rang, ce qui assure la convergence de l'intégrale de \eqref{eq:defLaplace} pour $\mathrm{Re}(x)>\gf{1}{C}$.  

Considérons $$ \begin{array}{cccc}
\mathcal{F}^* : & \Qbar\left[z,\gf{\mathrm{d}}{\mathrm{d}z} \right] & \longrightarrow & \Qbar\left[z,\gf{\mathrm{d}}{\mathrm{d}z} \right] \\ 
 & z & \longmapsto & \gf{\mathrm{d}}{\mathrm{d}z} \\ 
 & \gf{\mathrm{d}}{\mathrm{d}z} & \longmapsto & -z
\end{array} $$ la \emph{transformée de Fourier-Laplace} des opérateurs différentiels. Il s'agit d'un automorphisme d'anneaux non commutatifs de $\Qbar[z,\mathrm{d}/\mathrm{d}z]$ d'ordre 4 et d'inverse $\overline{\mathcal{F}^*}=s \circ \mathcal{F}^* \circ s$, où $s$ est la symétrie définie par  $\left(s(z),s(\mathrm{d}/\mathrm{d}z) \right)=\left(\mathrm{d}/\mathrm{d}z,z \right)$.

Les propriétés essentielles des transformées de Laplace et Fourier-Laplace sont rappelées dans \cite[p. 716]{AndregevreyI}. Signalons que ce qu'André entend par transformée de Fourier-Laplace est, avec nos notations, $\overline{\mathcal{F}^*}$, ce qui n'influe pas sur la définition de la notion d'$E$-opérateur (voir \cite[\S 4.1, p. 720]{AndregevreyI}). De plus, l'argument de \cite[\S 4.1, p. 720]{AndregevreyI} prouve \emph{mutatis mutandis} que $L$ est un $E$-opérateur \emph{au sens large} si et seulement si $\mathcal{F}^*(L)$ est un $G$-opérateur \emph{au sens large}.

Les faits suivants découlent des calculs de \cite[p. 25]{Beukers} : 
\begin{enumerateth}[label=\textbf{\textit{\roman*)}}]
\item Si $f$ est une $E$-fonction \emph{au sens large}, alors sa transformée de Laplace $g(x)=\mathcal{F}(f)(x)$ est une $G$-fonction \emph{au sens large} en $1/x$. La transformée de Laplace établit ainsi une correspondance entre $E$- et $G$-fonctions \emph{au sens large}. Ceci est également vrai au sens strict.
\item Si $L$ est un opérateur différentiel tel que $L(\mathcal{F}(f))=0$, alors on a $\mathcal{F}^*(L)(f)=0$, ce qui explique le nom donné à $\mathcal{F}^*$.
\end{enumerateth}

\medskip

André a défini un $E$-opérateur \emph{au sens strict} comme la transformée de Fourier-Laplace d'un $G$-opérateur \emph{au sens strict}. Par analogie, on définit une notion d'$E$-opérateur \emph{au sens large}.

\begin{defi} \label{def:Eoplarge}
On désigne par \emph{$E$-opérateur} au sens large la transformée de Fourier-Laplace d'un $G$-opérateur \emph{au sens large}.
\end{defi}

Les propriétés de morphisme d'anneaux de la transformée de Fourier-Laplace permet de déduire que l'ensemble des $E$-opérateurs \emph{au sens large} a les propriétés algébriques de l'ensemble des $G$-opérateurs \emph{au sens large} (proposition \ref{prop:algpropGop}). C'est ce qu'affirme la proposition suivante. 

\begin{prop} \label{coro:stabEop}
Soient $L$ et $M$ des $E$-opérateurs \emph{au sens large}. Alors :
\begin{enumerateth}[label=\textbf{\alph*)}]
  \item Tout diviseur à droite de $L$ dans $\Qbar(z)[\partial]$ est un $E$-opérateur \emph{au sens large}.
  \item L'opérateur adjoint $L^*$ de $L$ est un $E$-opérateur \emph{au sens large}. 
  \item Il existe un multiple commun à gauche de $L$ et $M$ qui est un $E$-opérateur \emph{au sens large} (propriété de Ore à gauche).
  \item Le produit $LM$ est un $E$-opérateur \emph{au sens large}.
  \end{enumerateth}
\end{prop}

Ceci découle des propositions  \ref{prop:algpropGop} et \ref{prop:tailleproduitgop} et du fait que $\mathcal{F}^*$ est un morphisme d'anneaux de l'ensemble des $G$-opérateurs \emph{au sens large} vers l'ensemble des $E$-opérateur \emph{au sens large}.

De plus, pour le point \textbf{b)}, on rappelle que l'opérateur adjoint $\left(\mathcal{F}^*(L)\right)^*$ est la transformée de Fourier-Laplace $\mathcal{F}^*(L^*)$ de $L^*$ (cf \cite[\S V.3.6]{Malgrange}).

Un résultat similaire vaut pour les $E$-opérateurs \emph{au sens strict}, comme cela avait été remarqué par André \cite[p. 720]{AndregevreyI}.

\medskip

L'appellation d'\og{}$E$-opérateur \fg{} \emph{au sens large} est justifiée par la proposition suivante, adaptation du résultat déjà connu dans le cas strict \cite[Théorème 4.2, p. 720]{AndregevreyI}. C'est une conséquence du théorème \ref{th:acklarge} de l'introduction.

\begin{prop} \label{prop:EfonannuleeEop}
Toute $E$-fonction \emph{au sens large} est solution d'une équation différentielle $L(y(z))=0$, où $L$ est un $E$-opérateur \emph{au sens large}.

De manière générale, soit une série Nilsson-Gevrey \emph{au sens large} de type arithmétique et holonome d'ordre $-1$, c'est-à-dire une fonction \begin{equation}\label{eq:serieNG} f(z)=\sum_{(\alpha, k,\ell) \in S} c_{\alpha,k,\ell} z^{\alpha} \log(z)^k f_{\alpha,k,\ell}(z),\end{equation} où $S \subset \Q \times \N^2$ est un ensemble fini, $c_{\alpha,k,\ell} \in \C$, et les $f_{\alpha,k,\ell}(z)$ sont des $E$-fonctions \emph{au sens large}. Alors $f$ est solution d'un $E$-opérateur \emph{au sens large}.
\end{prop}

Pour la démonstration, nous aurons besoin de la proposition \ref{prop:4chudnovskynilssongevreyandreI} ci-dessus, qui affirme que toute série Nilsson-Gevrey de type arithmétique et holonome d'ordre $0$ \emph{au sens large} est solution d'un $G$-opérateur \emph{au sens large}.

\begin{dem}[de la proposition \ref{prop:EfonannuleeEop}]
Soit $f(z)$ une $E$-fonction \emph{au sens large}. Alors il existe une $G$-fonction $g$ telle que $\mathcal{F}(f)(x)=g(1/x)$. 

Selon le théorème \ref{th:acklarge}, l'opérateur minimal $M_0$ de $g(u)$ sur $\Qbar(u)$ est un $G$-opérateur \emph{au sens large}. Donc comme la notion de $G$-opérateur \emph{au sens large} est invariante par changement de variable $u=1/x$ (Proposition \ref{prop:galochkinLargeInvChgeVar} \textbf{b)} ci-dessus), l'opérateur minimal $M$ de $\mathcal{F}(f)$ sur $\Qbar(x)$ est un $G$-opérateur \emph{au sens large}.

Or, l'assertion \textbf{\textit{ii)}} précédent la définition \ref{def:Eoplarge} implique que $M(\mathcal{F}(f))=0 \Rightarrow \mathcal{F}^*(M)(f)=0$, de sorte que $\mathcal{F}^*(M)$ est un $E$-opérateur \emph{au sens large} annulant $f(z)$.

\medskip

Traitons le cas plus général où $f(z)$ est de la forme \eqref{eq:serieNG}. 
Soit $L \in \Qbar\left[z, \ddz\right]$ tel que $L(f(z))=0$. Alors, en notant $g(z)=\mathcal{F}(f)(z)$, selon les formules (5.3.7) et (5.3.8) de \cite[pp. 729--730]{AndregevreyI}, on peut trouver $\rho \in \N$ tel que $\gf{\mathrm{d}^\rho}{\mathrm{d}z^{\rho}} \overline{\mathcal{F}^*}(L)(g(z))=0.$ Donc $h(z)$ est holonome. De plus, selon le lemme \ref{lem:laplacegevrey} ci-dessous, on peut écrire
$$g(z)=\sum\limits_{(\alpha,k,\ell) \in T} \lambda_{\alpha,k,\ell} z^{\alpha} \log(z)^k g_{\alpha,k,\ell}\left(\gf{1}{z}\right),$$ où $T \subset \Q \times \N^2$ est un ensemble fini, $\lambda_{\alpha,k,\ell} \in \C$, et les $g_{\alpha,k,\ell}(z)$ sont des $G$-fonctions au sens large. 

 Ainsi, selon la proposition \ref{prop:4chudnovskynilssongevreyandreI}, $g(1/z)$ est annulée par un $G$-opérateur \emph{au sens large}. En effectuant le changement de variable $u=1/z$, on obtient en conséquence de la proposition \ref{prop:galochkinLargeInvChgeVar} \textbf{b)} que $g(z)$ est annulée par un $G$-opérateur \emph{au sens large} $M$. On a alors, selon \cite[(5.3.7) et (5.3.8), p. 729]{AndregevreyI}, $\mathcal{F}^*(M)(f(z))=0$ et $\mathcal{F}^*(M)$ est un $E$-opérateur \emph{au sens large} par définition. Ceci est le résultat voulu.
\end{dem}

On peut traiter de manière plus directe le cas où $f(z)=z^{\alpha} F(z)$, où $F$ est une $E$-fonction \emph{au sens large}. On utilise pour cela la même formule que celle donnée dans \cite[p. 721]{AndregevreyI}, c'est-à-dire $$\mathcal{F}(f)(x)=\Gamma(\alpha+1) x^{-\alpha-1} \left[ \left(1-\gf{1}{x} \right)^{-\alpha} \star x \mathcal{F}(F)(x)\right],$$ où $\star$ désigne le produit de Hadamard. On voit que $\left(1-\gf{1}{x} \right)^{-\alpha} \star x \mathcal{F}(F)(x)$ est une $G$-fonction \emph{au sens large} en $1/x$, car $\alpha \in \Q$ et $F$ est une $E$-fonction \emph{au sens large}. 

On se sert alors du lemme suivant :
\begin{lem} \label{lem:Gopinvmultzpuissancea}
Soit $L$ un $G$-opérateur \emph{au sens large} (resp. \emph{au sens strict}), $a \in \Q$ et $L_{a}=z^{a} L z^{-a} \in \Qbar(z)[\mathrm{d}/\mathrm{d}z]$. Alors $L_a$ est un $G$-opérateur \emph{au sens large} (resp. \emph{au sens strict}).
\end{lem}

Selon ce lemme, $\mathcal{F}(f)(x)$ est donc solution d'un $G$-opérateur \emph{au sens large}. Le raisonnement ci-dessus s'applique donc encore à $\mathcal{F}(f)(x)$. 

La preuve du lemme \ref{lem:Gopinvmultzpuissancea} ci-dessous est plus directe que celle de la proposition \ref{prop:4chudnovskynilssongevreyandreI}.

\begin{dem}[du lemme \ref{lem:Gopinvmultzpuissancea}]
On fait la démonstration dans le cas large.

Les solutions de l'équation $L_a(w(z))=0$ sont les $z^a y(z)$ quand $L(y(z))=0$. On note $G$ (resp. $G_a$) la matrice compagnon de $L$ (resp. $L_a$). Soit $$Y(z)=\begin{pmatrix}
f_1 & \dots & f_{\mu} \\ 
\vdots &  & \vdots \\ 
f_1^{(\mu-1)} & \dots & f_{\mu}^{(\mu-1)}
\end{pmatrix}$$ une matrice fondamentale de solutions du système $y'=Gy$. En notant $g_i(z)=z^a f_i(z)$, on voit donc que $$Z(z)=\begin{pmatrix}
g_1 & \dots & g_{\mu} \\ 
\vdots &  & \vdots \\ 
g_1^{(\mu-1)} & \dots & g_{\mu}^{(\mu-1)}
\end{pmatrix}$$ est une matrice fondamentale de solutions du système $y'=G_a y$. De plus, on a $G_s=Y^{(s)} Y^{-1}$ et $(G_a)_s = Z^{(s)} Z^{(-1)}$ pour tout $s$.

Calculons la matrice $Z$ en fonction de $Y$. Selon la formule de Leibniz pour tout $\ell \in \{ 0, \dots, \mu-1 \}$, $$g_i^{(\ell)}(z)=\sum_{k=0} \dbinom{\ell}{k} (z^a)^{(\ell-k)} f_i^{(k)}(z)=\sum_{k=0}^{\ell} (a-(\ell-k)+1)_{\ell-k} z^{a-(\ell-k)} f_i^{(k)}(z).$$

Ainsi, $Z=UY$, où $U := \left(u_{\ell,k} \right)_{k,\ell} \in z^a \mathcal{M}_n(\Qbar(z))$, $u_{\ell,k} :=\dbinom{\ell}{k} (a-(\ell-k)+1)_{\ell-k} z^{a-(\ell-k)}$ est une matrice triangulaire inférieure inversible. En notant $v_{\ell,k}=\dbinom{\ell}{k} (a-(\ell-k)+1)_{\ell-k}$, on a, pour tout $k \leqslant \ell$, $$\gf{u_{\ell,k}^{(s)}}{s!}= v_{\ell,k} \gf{(z^{a-(\ell-k)})^{(s)}}{s!}=v_{\ell,k} \gf{(a-(\ell-k)-s+1)_s}{s!} \gf{(z^{a-(\ell-k)})^{(s)}}{s!}.$$ Or, pour tout $\gamma \in \Q$, le dénominateur commun de $(\gamma)_0/0!, \dots, (\gamma)_s/s!$ divise $\den(\gamma)^{2s}$ (voir \cite[Lemma 10, p. 334]{Lepetit2}). Donc on peut trouver une constante $c \in \N^*$ tel que $\den(a)^{2s} c \gf{U^{(s)}}{s!} \in z^{a} \mathcal{M}_n\left(\Oal_{\Qbar}[z,1/z]\right)$.

Donc par la formule de Leibniz, on a pour $s \in \N$, $$\gf{Z^{(s)} Z^{-1}}{s!}=(UY)^{(s)}=\sum_{k=0}^{s} \gf{1}{s!}\dbinom{s}{k} U^{(s-k)} Y^{(k)} Y^{-1} U^{-1}=\sum_{k=0}^{s} \gf{U^{(s-k)}}{(s-k)!} \gf{G_k}{k!} U^{-1}.$$

Considérons $T \in \Qbar[z]$ tel que $TG \in \mathcal{M}_n(\Qbar[z])$ et $q \in \N$ tel que $z^q U \in z^a \mathcal{M}_n(\Qbar[z])$. Comme $G$ vérifie la condition de Galochkin \emph{au sens large}, on peut trouver pour tout $s$ un entier $q_s$ non nul tel que $q_s \gf{T^k G_k}{k!} \in \mathcal{M}_n(\Oal_{\Qbar}[z])$ et vérifiant, pour tout $\varepsilon >0$, $q_s \leqslant s!^{\varepsilon}$ pour tout $s$ suffisamment grand relativement à $\varepsilon$. Prenons un polynôme $D$ tel que $DU^{-1} \in z^a \mathcal{M}_n(\Oal_{\Qbar}[z])$.

Ainsi, on obtient finalement $$z^{qs} D T^s q_s \den(a)^{2s} c \gf{Z^{(s)} Z^{-1}}{s!} \in \mathcal{M}_n(\Oal_{\Qbar}[z])$$ et $\tilde{q}_s :=z^{qs} D T^s q_s \den(a)^{2s}$ satisfait, pour tout $\varepsilon >0$, la condition $\tilde{q}_s \leqslant s!^{\varepsilon}$ pour tout $s$ suffisamment grand relativement à $\varepsilon$. Le cas strict se traite \emph{mutatis mutandis}. Ceci achève la preuve du lemme \ref{lem:Gopinvmultzpuissancea}
\end{dem}

\subsection{Structure des $E$-opérateurs \emph{au sens large}}

Le but de cette partie est de déduire du théorème de structure des $G$-opérateurs \emph{au sens large} (théorème \ref{th:complementACKlarge}) le théorème \ref{th:andrelarge} sur les équations différentielles satisfaites par les $E$-fonctions \emph{au sens large}, en adaptant la théorie d'André des $E$-opérateurs \emph{au sens strict} développée dans \cite[pp. 715--724]{AndregevreyI}.

\bigskip

Par analogie avec les \Zemla -fonctions introduites dans \cite[p. 713]{AndregevreyI}, on définit la famille de fonctions suivante, qui apparaît dans le théorème de structure des $E$-opérateurs \emph{au sens large} (théorème \ref{th:andrelarge} ci-dessous) :

\begin{defi} \label{def:Zfonctionlarge}
Une \Zemla -fonction \emph{au sens large} est une série $f(z)=\sum\limits_{n=0}^{\infty} a_n n! z^n \in \Qbar\llbracket z \rrbracket$ vérifiant la condition \textbf{a)} de la définition \ref{def:gfonctionlarge} et telle que les $a_n$ vérifient les conditions \textbf{b)} et \textbf{c)} de la définition \ref{def:gfonctionlarge}.
\end{defi}

Dans \cite[Théorème 4.3, p. 721]{AndregevreyI}, André a précisé la structure des $E$-opérateurs \emph{au sens strict}. Le théorème suivant est l'analogue de son résultat \emph{au sens large}.  La démonstration, adaptée de celle donnée par André dans le cas strict, fera l'objet de la partie \ref{subsec:preuveandrelarge}.

\begin{Th} \label{th:andrelarge}
Soit $L \in \Qbar(z)[\mathrm{d}/\mathrm{d}z]$ un $E$-opérateur \emph{au sens large} d'ordre $\nu$. Alors
\begin{enumerateth}[label=\textbf{\alph*)}]
 \item Les seules singularités de $L$ sont $0$ et l'infini, qui est en général un point irrégulier ; 
 \item L'opérateur $L$ est singulier régulier en $0$. Les exposants de $L$ en $0$ sont rationnels, ceux qui ne sont pas entiers sont (modulo $\Z$ et comptés sans multiplicité), les exposants en $\infty$ de $\mathcal{F}^*(L)$.
 \item Les pentes du polygone de Newton de $L$ en l'infini sont dans $\{0,1\}$. 
 \item Il existe une base de solutions en $0$ de l'équation $L(y(z))=0$ de la forme $$\left(F_1(z), \dots, F_{\nu}(z) \right) \cdot z^{\Gamma_0},$$ où les $F_j$ sont des $E$-fonctions \emph{au sens large} et $\Gamma_0 \in \mathcal{M}_{\nu}(\Q)$ est triangulaire supérieure.
 \item Il existe une base de solutions en $\infty$ de l'équation $L(y(z))=0$ de la forme $$ \left(f_1\left(\gf{1}{z}\right), \dots, f_{\nu}\left(\gf{1}{z} \right)\right) \cdot \left(\gf{1}{z}\right)^{\Gamma_{\infty}} e^{-\Delta z},$$ où les $f_j$ sont des \Zemla -fonctions \emph{au sens large}, $\Delta$ est la matrice diagonale dont les coefficients sont les singularités à distance finie de $\mathcal{F}^*(L)$ (comptées avec multiplicité), et $\Gamma_{\infty}$ désigne une matrice triangulaire supérieure à coefficients dans $\Q$ qui commute avec $\Delta$.
\end{enumerateth}
\end{Th}

En conséquence de la proposition \ref{prop:EfonannuleeEop} et du théorème \ref{th:andrelarge}, si $f(z)$ est une $E$-fonction \emph{au sens large}, toute les singularités différentes de $0$ et $\infty$ de son opérateur minimal non nul $L$ sur $\Qbar(z)$ sont apparentes. En effet, l'opérateur minimal de $f(z)$ est facteur à droite d'un opérateur n'ayant que $0$ et $\infty$ pour singularités. Ce fait a été énoncé par André dans \cite[p. 747]{AndregevreyII}. De plus, la singularité $0$ vérifie les propriétés \textbf{b)} et \textbf{d)} du théorème \ref{th:andrelarge} (cf \cite[Corollaire 4.4, p. 724]{AndregevreyI}).

\bigskip

\begin{rqus}
\begin{enumerateth}[label=\textit{\textbf{\roman*)}}]
 \item Le lemme 4, p. 518, de \cite{Gorelov2004} affirme que tout point $\xi \in \Qbar^*$ est un point régulier ou singularité régulière de l'équation différentielle minimale d'une $E$-fonction, ce qui est une conséquence du théorème \ref{th:andrelarge}. En effet, si $L \in \Qbar(z)[\mathrm{d}/\mathrm{d}z]$ est un opérateur singulier régulier en $a \in \Qbar$ et $M$ divise $L$ dans $\Qbar(z)[\mathrm{d}/\mathrm{d}z]$, alors $M$ est singulier régulier en $a$. Ceci est un résultat connu qui découle de la définition de la régularité des singularités. La preuve de Gorelov exploite directement les propriétés des $E$-fonctions \emph{au sens large} sans passer par la transformée de Fourier-Laplace.
\medskip
\item Gorelov \cite[Theorem 1, p. 514]{Gorelov2004} a obtenu le résultat suivant sur les $E$-fonctions d'ordre 2 : si $f(z)$ une $E$-fonction \emph{au sens large} dont l'opérateur minimal $L$ sur $\Qbar(z)$ est d'ordre $2$, on peut écrire
\begin{equation}\label{eq:probsiegelRR} f(z) = a(z)e^{\mu z} {}_1 F_1 (\alpha ; \beta ; \lambda z) + b(z)e^{\mu z} _1 F'_1 (\alpha ; \beta ; \lambda z),\end{equation} où $a(z), b(z) \in \Qbar[z]$, $\lambda \in \Q$, $\mu \in \Q$, et $\alpha \in \Q, \beta \in \Q \setminus \Z_{\leqslant 0}$.

Dans \cite{RivoalRoques}, Rivoal et Roques ont démontré, en s'appuyant sur la théorie d'André des $E$-opérateurs, un énoncé analogue à \eqref{eq:probsiegelRR} dans le cas où $f(z)$ est une $E$-fonction \emph{au sens strict}, en ajoutant que $\alpha - \beta \not\in \Z$ mais avec l'affirmation plus faible que $a(z), b(z) \in \Qbar(z)$.
  
  Le théorème \ref{th:andrelarge} ci-dessus permet d'obtenir l'analogue du résultat de Rivoal et Roques pour les $E$-fon\-ctions \emph{au sens large} en adaptant la démarche de \cite{RivoalRoques}. En effet, selon le théorème \ref{th:andrelarge}, l'opérateur minimal de $f(z)$ sur $\Qbar(z)$ vérifie les hypothèses de \cite[Theorem 5, p. 8]{RivoalRoques}, ainsi que les conditions (1), (2) et (3) de \cite[\S 5.2, p. 12--13]{RivoalRoques}. 
  
Le théorème de Gorelov constitue une réponse positive au \emph{problème de Siegel} pour les $E$-fonctions d'ordre 2. Le problème de Siegel consiste en la question suivante : \textit{toute $E$-fonction est-elle dans la $\Qbar[z]$-algèbre engendrée par les fonctions hypergéométriques de la forme $$\sum_{n=0}^{\infty} \gf{(a_1)_n \dots (a_p)_n}{(b_1)_n \dots (b_q)_n} \lambda^n z^{n(q-p)},$$ quand $0 \leqslant p<q$, $a_1, \dots, a_p \in \Q$, $b_1, \dots, b_q \in \Q \setminus \Z_{\leqslant 0}$ et $\lambda \in \Qbar$ ?}

Dans le cas général, on sait que la réponse à cette question est négative. En effet, dans \cite{FresanJossen}, Fresán et Jossen ont exhibé un exemple de $E$-fonction non hypergéométrique.

\end{enumerateth}
\end{rqus}

\subsection{Démonstration du théorème \ref{th:andrelarge}} \label{subsec:preuveandrelarge}

Cette section est consacrée à la preuve du théorème \ref{th:andrelarge}, en suivant et adaptant la démonstration du Théorème 3.4.1 de \cite{AndregevreyI} détaillée dans \cite[\S 5, pp. 725--735]{AndregevreyI}.

\medskip

Puisqu'un $G$-opérateur \emph{au sens large} est fuchsien selon le théorème \ref{th:acklarge}, le point \textbf{a)}, ainsi que la régularité de $L$ en $0$, découlent des calculs menés dans \cite[p. 25]{Beukers}. En réalité, pour avoir la régularité de $L$ en $0$, il suffit d'invoquer le fait que $\infty$ est une singularité régulière du $G$-opérateur \emph{au sens large} $\mathcal{F}^*(L)$.

Passons au point \textbf{b)}. Pour $M \in \Qbar(z)[\mathrm{d}/\mathrm{d}z]$, on note $\mathrm{NR}(M)$ le polygone de Newton-Ramis de $M$, tel que défini dans \cite[pp. 7--8]{Ramis}. On sait que $M$ est un opérateur singulier régulier en $0$ et $\infty$ si et seulement si $\mathrm{NR}(M)$ admet pour seule pente $0$, c'est-à-dire si on peut trouver $a,b,c \in \R$ tels que $$\mathrm{NR}(M)=\left\lbrace (u,v) \in \R^2, a \leqslant u \leqslant b \;\; \text{et} \;\; v \leqslant c \right\rbrace.$$

Or, selon \cite[p. 726]{AndregevreyI}, $\mathrm{NR}(\overline{\mathcal{F}^*}(M))=\tau(\mathrm{NR}(M))$, où $\tau: (u,v) \mapsto (u+v,-v)$. Comme $\tau$ est involutive, on en déduit que $\mathrm{NR}(\mathcal{F}^*(M))=\tau(\mathrm{NR}(M))$. De plus, si $(u,v) \in \R^2$ et $\tau(u,v)=(u',v')$, on a $$\begin{cases} a \leqslant v \leqslant b \\  u \leqslant c \end{cases} \ssi \begin{cases} -b \leqslant v' \leqslant -a \\ v' \leqslant c-u', \end{cases}$$ de sorte que $M$ est singulier régulier en $0$ et $\infty$ si et seulement si les seules pentes du polygone de Newton $\mathrm{NR}(\mathcal{F}^*(M))$ sont $0$ et $-1$. Puisque les pentes négatives correspondent aux pentes du polygone de Newton classique de $\mathcal{F}^*(M)$ en $\infty$, il en découle que $M$ est singulier régulier en $0$ et $\infty$ si et seulement si $\mathcal{F}^*(M)$ est singulier régulier en $0$ et irrégulier de pentes dans $\{0,1\}$ en $\infty$. 

En appliquant cette assertion au $G$-opérateur $M$ tel que $L=\mathcal{F}^*(M)$, qui est un opérateur fuchsien selon le théorème \ref{th:acklarge}, on obtient le point \textbf{c)}.

On peut répéter le même raisonnement au voisinage de tout point $a$ de $\Qbar$. Il découle du caractère fuchsien de $M$ que les pentes de $\mathcal{F}^*(M_a)$, où $M_a$ est l'opérateur obtenu par changement de variable $u=z-a$, sont dans $\{0,1\}$.

On obtient ainsi en reproduisant \emph{mutatis mutandis} l'argument de \cite[p. 720]{AndregevreyI}, qui utilise le théorème de Turritin-Levelt, une base de solutions en $\infty$ de l'équation $L(y(z))=0$ de la forme $$ \left(\hat{f}_1\left(\gf{1}{z}\right), \dots, \hat{f}_{\nu}\left(\gf{1}{z} \right)\right) \left(\gf{1}{z}\right)^{\Delta_{\infty}} e^{-\Delta z},$$ où les $\hat{f}_j$ sont des séries de Laurent à coefficients dans $\Qbar$, $\Delta$ est une matrice diagonale, et $\Gamma_{\infty}$ désigne une matrice triangulaire supérieure sous forme de Jordan qui commute avec $\Delta$.

\medskip

Par construction (voir \cite[p. 725]{AndregevreyI}), si $M \in \Qbar(z)[\mathrm{d}/\mathrm{d}z]$ est fuchsien, alors $\overline{\mathcal{F}^*}(M)$ est \emph{de type exponentiel} selon la terminologie de \cite[p. 195]{Malgrange} (ou \emph{exponentiel élémentaire} avec les mots d'André). On peut donc utiliser le raisonnement de \cite[\S 5.2, pp. 726--728]{AndregevreyI}, fondé sur des arguments généraux sur les microsolutions des opérateurs de type exponentiel.
Il nous assure alors que, puisque $\mathcal{F}^*(L)$ est fuchsien à exposants rationnels, les exposants non entiers de $L$ sont (modulo $\Z$ et comptés sans multiplicité) les exposants en $\infty$ de $\mathcal{F}^*(L)$. Ceci prouve le point \textbf{b)} du théorème \ref{th:andrelarge}. 

\bigskip

La partie 5.3 pp 728--730 de \cite{AndregevreyI} a permis de construire à partir de la transformée de Laplace deux applications injectives $\C$-linéaires \begin{equation}\label{eq:phietpsi} \overline{\C\llbracket z \rrbracket}[\log z] \mathrel{\mathop{\rightleftarrows}^{\varphi}_{\psi}} \overline{\C\left\llbracket \gf{1}{z} \right\rrbracket} \left[\log \gf{1}{z} \right],\end{equation} $\overline{\C\llbracket u \rrbracket}$ désignant l'anneau des séries de Puiseux en $u$. 

Un ingrédient essentiel de la preuve des points \textbf{d)} et \textbf{e)} du théorème \ref{th:andrelarge} est le lemme suivant, analogue de \cite[Proposition 5.4.1]{AndregevreyI}. Il est également utilisé dans la preuve de la proposition \ref{prop:EfonannuleeEop} ci-dessus. On rappelle que les ensembles $\mathrm{NGA}_h^{\ell}\{z\}_{s}$
ont été introduits dans la définition \ref{def:NGAlarge}.

\begin{lem} \label{lem:laplacegevrey}
Les applications définies en \eqref{eq:phietpsi} induisent des applications $\C$-linéaires injectives $$\mathrm{NGA}_h^{\ell}\{z\}_{-1} \rightleftarrows \mathrm{NGA}_h^{\ell}\left\{ \gf{1}{z} \right\}_0\;, \qquad \mathrm{NGA}_h^{\ell}\{z\}_{0} \rightleftarrows \mathrm{NGA}_h^{\ell}\left\{ \gf{1}{z} \right\}_1. $$
De plus, on a $$\varphi^{-1}\left(\mathrm{NGA}_h^{\ell}\left\{ \gf{1}{z} \right\}_0\right)=\mathrm{NGA}_h^{\ell}\{z\}_{-1}, \quad \psi^{-1}\left(\mathrm{NGA}_h^{\ell}\{z\}_{-1}\right)=\mathrm{NGA}_h^{\ell}\left\{ \gf{1}{z} \right\}_0$$ et de même pour $\mathrm{NGA}_h^{\ell}\{z\}_{0}$ et $ \mathrm{NGA}_h^{\ell}\left\{ \gf{1}{z} \right\}_1$.
\end{lem}

\begin{dem}
Il s'agit de montrer l'analogue \emph{au sens large} des assertions (5.4.2) et (5.4.3) de la démonstration d'André dans le cas strict \cite[pp. 731--732]{AndregevreyI}. Dans ce qui suit, on appelle \og{} suite holonome \fg{} toute suite $(a_n)_n$ telle que $\sum\limits_{n=0}^{\infty} a_n z^n$ est solution d'une équation différentielle non nulle à coefficients dans $\Qbar(z)$

Le seul point de l'argumentation d'André qui nécessite une adaptation \emph{au sens large} consiste à prouver que les formules (5.4.4) et (5.4.7) de \cite{AndregevreyI}, issus de calculs génériques sur la transformée de Laplace, définissent des suites holonomes vérifiant les conditions \textbf{b)} et \textbf{c)} de la définition \ref{def:gfonctionlarge}.

Soient donc $\underline{a}=(a_n)_{n \in \N}$ une suite de nombres algébriques holonome vérifiant les conditions \textbf{b)} et \textbf{c)} de la définition \ref{def:gfonctionlarge}, $\alpha \in \Q$ et $k \in \N$.
On définit respectivement  la $E$-fonction  \emph{au sens large} et la \Zemla-fonction \emph{au sens large} en la variable $1/z$ 
$$F_{\underline{a}}(z)=\sum\limits_{n=0}^{\infty} \gf{a_n}{n!} z^n \quad \text{et} \quad \mathfrak{f}_{\underline{a}}=\sum\limits_{n=0}^{\infty} a_n {n!} z^{-n}.$$ 
Alors André a prouvé que 
\begin{align}
    \mathcal{F}\left(z^{\alpha} \log^k(z) F_{\underline{a}}(z)\right)=\sum_{j=0}^{j_0} \sum_{n=0}^{\infty} b_{n,j} z^{-\alpha-1-n} \log^j z \label{eq:laplaceFa}\\
    \text{et} \;\; \mathcal{F}\left(z^{\alpha} \log^k(z) \mathfrak{f}_{\underline{a}}(z)\right)=\sum_{j=0}^{j_1} \sum_{n=0}^{\infty} c_{n,j} z^{-\alpha-1-n} \log^j z, \label{eq:laplacefa}
\end{align}
où, si $\alpha \not\in \Z_{-}$ (resp. $\alpha \not\in \N$) la suite $(b_{n,j_0})_n$ (resp. $(c_{n,j_1})_n$) correspondant à la plus haute puissance du logarithme apparaissant dans \eqref{eq:laplaceFa} (resp. \eqref{eq:laplacefa}) est définie par l'équation (5.4.4) de \cite{AndregevreyI} :
\begin{align*}
    b_{n,j_0}=b_{n,k} &=(-1)^k \Gamma(\alpha+1) \gf{(\alpha+1)_n}{n!} a_n  \\
    c_{n,j_1}=c_{n,k} &= (-1)^{k+n} \Gamma(\alpha+1)  \gf{n!}{(-\alpha-1)_n} a_n \qquad \text{si} \;\; \alpha \not\in \Z \\
    c_{n,j_1}=c_{n,k+1} &= (-1)^{\alpha+n} \gf{1}{(k+1) \dbinom{n}{-\alpha-a}} a_n \qquad \text{si} \;\; \alpha \in \Z, \; \alpha <0
\end{align*}
et dans le cas contraire, comme l'a remarqué André, on se ramène au cas précédent en retranchant de $F_{\underline{a}}(z)$ (resp. $\mathfrak{f}_{\underline{a}}(z)$) un élément de $\Qbar[1/z]$ (resp. $\Qbar[z]$).

Il s'agit donc de démontrer que les suites $(b_{n,k})_n$ et $(c_{n,k})_n$ vérifient les points \textbf{b)} et \textbf{c)} de la définition \ref{def:gfonctionlarge}. 
Comme le maximum des dénominateurs (resp. des tailles) d'un quotient de symboles de Pochhammer $(u)_0/(v)_0, \dots, (u)_n/(v)_n$ a une croissance au plus géométrique en $n$ (cf \cite[\S 9, pp. 54--58]{Siegel}), (5.4.4) de \cite{AndregevreyI} fournit bien des suites dont le dénominateur et la taille croissent au plus en $n!^{\varepsilon}$, pour tout $\varepsilon >0$ fixé.

De plus, comme toute série hypergéométrique $\sum\limits_{n=0}^{\infty} \gf{(\alpha_1)_n \dots (\alpha_p)_n}{(\beta_1)_n \dots (\beta_q)_n} z^n$ est holonome et qu'un produit ou une somme de suites holonomes est holonome, l'équation (5.4.4) de \cite{AndregevreyI} définit bien des multiples dans $\C$ de suites holonomes.

\medskip

Soit maintenant $j \in \{0, \dots, k-1 \}$. De la même manière, les suites apparaissant dans la formule (5.4.7) de \cite{AndregevreyI} :
\begin{align*}
    b_{n,j} &= \gf{n!}{(-\alpha-1)_n} \rho_{n,j} a_n \qquad (\alpha \not\in \Z_{-}^{*}) \\
    c_{n,j} &=(-1)^n \gf{n!}{(-\alpha-1)_n} \rho_{-n,j} a_n \qquad \text{si} \;\; \alpha \not\in \Z \\
    c_{n,j} &= \gf{1}{\dbinom{n}{-\alpha-1}} \rho_{-n,j} a_n \qquad \text{si} \;\; \alpha \in \Z_{-}^*
\end{align*}
sont combinaisons linéaires à coefficients dans $\C$ de suites satisfaisant les conditions \textbf{b)} et \textbf{c)} de la définition \ref{def:gfonctionlarge}. En effet, la formule (5.4.6) de \cite{AndregevreyI} implique l'existence de suites $(r_{m,i})_{m \in \Z}$ pour $0 \leqslant i \leqslant k+1$ et de nombres complexes $\rho_0, \dots, \rho_{k+1}$ tels que $$ \forall m \in \Z, \;\; r_{m,j}=\sum\limits_{i=0}^{k+1} r_{m,i} \rho_i,$$ où pour tout $i$, $(r_{m,i})_{m \in \N}$ et $(r_{-m,i})_{m \in \N}$ satisfont les conditions \textbf{b)} et \textbf{c)} de la définition \ref{def:gfonctionstrict}. De plus, l'holonomie des suites $(r_{m,i})_m$ et $(r_{-m,i})_m$ est assurée par la formule (5.4.6) de \cite{AndregevreyI} : $$\rho_{m,j}=\rho_{m-1,j}-\gf{j+1}{\alpha+m} \rho_{m-1,j+1}.$$ Ainsi, la formule (5.4.7) de \cite{AndregevreyI} fournit bien des suites $(b_{n,j})$ et $(c_{n,j})$ qui sont combinaisons linéaires à coefficients dans $\C$ de suites holonomes. \end{dem}
 On démontre  le point \textbf{d)} en remplaçant \og{} $E$-opérateur \fg{} par \og{} $E$-opérateur \emph{au sens large} \fg{}, et de même pour les $G$-opérateurs, dans \cite[p. 734]{AndregevreyI}. En effet, la preuve d'André utilise les résultats de \cite[\S 5.3, pp. 728--730] {AndregevreyI} qui ne concernent pas les $E$-fonctions et la proposition~5.4.1 d'André dont le lemme \ref{lem:laplacegevrey} est l'analogue \emph{au sens large}. De plus, on prouve les trois faits suivants :

\begin{itemize}[label=\textbullet]
    \item Si $\Phi$ est un $E$-opérateur \emph{au sens large} et $\rho \in \N$, alors $\Psi := \gf{\mathrm{d}^{\rho}}{\mathrm{d}z^{\rho}} \overline{\mathcal{F}^*}(\Phi)$ est un $G$-opérateur \emph{au sens large}, car la transformée de Fourier-Laplace d'un $E$-opérateur \emph{au sens large} est un $G$-opérateur \emph{au sens large} et l'ensemble des $G$-opérateurs \emph{au sens large} est stable par produit.
    \item Selon le théorème \ref{th:complementACKlarge} appliqué en $\infty$, si $y(z) \in \overline{\C\left\llbracket \gf{1}{z} \right\rrbracket} \left[\log \gf{1}{z} \right]$ est tel que $\Psi(y(z))=0$, alors $y(z) \in \mathrm{NGA}_h^{\ell}\left\{ \gf{1}{z} \right\}_0$.
    \item Si $\mathcal{F}(y)(z) \in \mathrm{NGA}_h^{\ell}\left\{ \gf{1}{z} \right\}_0$ alors, selon le lemme \ref{lem:laplacegevrey}, on a $y(z) \in \mathrm{NGA}_h^{\ell}\left\{ z \right\}_{-1}$.
\end{itemize}

On obtient ainsi une base de solutions de l'équation $L(y(z))=0$ de la forme $$(y_1(z), \dots, y_{\nu}(z))=(F_1(z), \dots, F_{\nu}(z)) \cdot z^{\Gamma_0} \in \mathrm{NGA}_h^{\ell}\{z\}_1,$$ où $\Gamma_0 \in \mathcal{M}_n(\Q)$ est triangulaire supérieure.

\bigskip

Enfin, la preuve du point \textbf{e)} issue de \cite[\S 5.6, pp. 733--734]{AndregevreyI} s'adapte également au sens large. Elle utilise une nouvelle fois les résultats de \cite[\S 5.3, pp. 728--730] {AndregevreyI} ainsi que le lemme~\ref{lem:laplacegevrey}.

Les points spécifiques aux $E$-opérateurs à adapter -- et les arguments qui permettent cette adaptation -- sont les suivants :

\begin{itemize}[label=\textbullet]
\item \emph{$\Phi \otimes e^{-\zeta_{j} z}$ [...] est encore de type $E$, car $\overline{\mathcal{F}^*}(\Phi \otimes e^{-\zeta_{j} z})$ est le $G$-opérateur translaté de $\mathcal{F}^* \Phi$ par $\zeta_j$}. Dans le cas large, ceci découle de ce que l'ensemble des $G$-opérateurs \emph{au sens large} est invariant par changement de variable $u=z-a$, 
$a \in \Qbar$ (Proposition \ref{prop:galochkinLargeInvChgeVar} \textbf{b)}). De plus, la transformée de Fourier-Laplace d'un $E$-opérateur \emph{au sens large} est un $G$-opérateur \emph{au sens large}.
\item \emph{$\tilde{y}_j^{+}(-z)$ est solution du $G$-opérateur $\mathcal{F}^* \Phi$ [...], on a donc  $\tilde{y}_j^{+} \in \mathrm{NGA}\{ z \}_0$} (ici, $y_j^{+}$ la transformée de Laplace de $y_j$). Dans le cas large, ceci découle du théorème \ref{th:complementACKlarge}. De plus, $y_j^{+}$ est de type holonome, toujours selon le théorème \ref{th:complementACKlarge}. 
\end{itemize}

La conclusion de la preuve du point \textbf{e)} du théorème \ref{th:andrelarge} repose alors sur le lemme suivant, analogue au sens large de \cite[Proposition 5.6.3, p. 734]{AndregevreyI} :

\begin{lem}
Soient $\hat{f} \in \C\left(\left(\gf{1}{z}\right)\right)$ et $\alpha \in \Q$. Alors $\hat{f}$ est solution d'un $E$-opérateur \emph{au sens large} si et seulement si $z^{\alpha} \hat{f}$ est solution d'un $E$-opérateur \emph{au sens large}.
\end{lem}

\begin{dem}
La preuve est la même \emph{mutatis mutandis} que celle de \cite[Proposition 5.6.3, p. 734]{AndregevreyI}. En effet :
\begin{itemize}[label=\textbullet]
\item Tout $E$-opérateur \emph{au sens large} $\Phi$ a ses pentes à l'infini dans $\{0,1\}$ selon le point \textbf{c)} du théorème \ref{th:andrelarge}. Si $1$ est effectivement une pente de $\Phi$ à l'infini, alors toute solution $g \in \C\llbracket 1/z \rrbracket$ de $\Phi(y(z))=0$ au voisinage de $\infty$ est $1$-sommable au sens de Ramis (voir \cite[p. 34]{RamisDiv}) dans toute direction non singulière. Ceci est a fortiori vrai si $0$ est la seule pente de $\Phi$ à l'infini, car $\infty$ est alors un point singulier régulier de $ \Phi$, ce qui assure la croissance modérée des solutions de $\Phi(y(z))=0$ dans les secteurs. C'est cette $1$-sommabilité qui est nécessaire dans la preuve d'André que nous adaptons.
\item Soit $V$ un secteur de $\C$ inclus dans un plan fendu $\C \setminus \mathcal{D}$, où $\mathcal{D}$ est une demi-droite d'origine $0$. Selon le théorème \ref{th:andrelarge} \textbf{d)}, toute fonction $\tilde{y}(z)$ vérifiant $L(\tilde{y}(z))=0$ dans $V$ est la restriction d'un certain élément $y(z)$ de $\mathrm{NGA}_h^{\ell}\{ z\}_{-1}$. 
\item Si $w(z) \in \mathrm{NGA}_h^{\ell} \{ z \}_{-1}$, alors la proposition \ref{prop:EfonannuleeEop} nous assure de l'existence d'un $E$-opérateur \emph{au sens large} $L'$ tel que $L'(w(z))=0$.
\end{itemize}
\end{dem}

\subsection{Nouvelle preuve d'un théorème d'André sur les $E$-fonctions \emph{au sens large}}

Le but de cette partie est de déduire du théorème de structure des $E$-opérateurs \emph{au sens large} (théorème \ref{th:andrelarge}) une nouvelle démonstration d'un résultat diophantien sur les valeurs des $E$-fonctions \emph{au sens large}, dû à André \cite{Andre2014}, qui généralise le théorème fondamental suivant de Siegel-Shidlovskii.

\begin{Th}[Siegel--Shidlovskii, 1929/1956, \cite{Shidlovskii}, p. 139] \label{th:siegelshidlovskii}
Soit $(f_1, \dots, f_n)$ une famille de $E$-fonctions \emph{au sens large}. Soit $\mathbf{f}=^{t} (f_1, \dots, f_n)$, supposons qu'il existe $A \in \mathcal{M}_n(\Qbar(z)$ tel que $\mathbf{f}'=A \mathbf{f}$. Prenons $\alpha \in \Qbar$ tel que $\alpha T(\alpha) \neq 0$, où $T(z) \in \Qbar[z]$ est tel que $T(z)A(z) \in \mathcal{M}_n(\Qbar[z])$.

Alors le degré de transcendance sur $\Qbar$ de $(f_1(\alpha), \dots, f_n(\alpha))$ est égal au degré de transcendance sur $\Qbar(z)$ de $(f_1(z), \dots, f_n(z)).$
\end{Th}

Ce théorème généralise le théorème de Lindemann-Weierstrass. En effet, si $(\alpha_1, \dots, \alpha_n) \in \Qbar^n$ est une famille libre sur $\Q$, alors $(e^{\alpha_1 z}, \dots, e^{\alpha_n  z})$ est une famille de $E$-fonctions \emph{au sens large} vérifiant les hypothèses du théorème dont les composantes sont algébriquement indépendantes sur $\Qbar(z)$, de sorte qu'en évaluant en $1$, le théorème \ref{th:siegelshidlovskii} nous assure que $e^{\alpha_1}, \dots, e^{\alpha_n}$ sont algébriquement indépendants sur $\Qbar$.

Beukers a ensuite démontré le théorème suivant dans \cite{Beukers2006}. Il constitue un raffinement du théorème \ref{th:siegelshidlovskii} dans le cas strict. En effet, le théorème \ref{th:siegelshidlovskii} est vrai quant à lui pour les $E$-fonctions \emph{au sens large}.

\begin{Th}[Beukers, 2006, \cite{Beukers2006}, Theorem 1.3, p. 370] \label{th:beukersefonctions}
Soit $(f_1, \dots, f_n)$ une famille de $E$-fonctions \emph{au sens strict} vérifiant les hypothèses du théorème \ref{th:siegelshidlovskii}. Alors pour tout polynôme homogène $P \in \Q[X_1, \dots, X_n]$ tel que $P(f_1(\alpha), \dots, f_n(\alpha))=0$, il existe un polynôme $Q \in \Q[Z, X_1, \dots, X_n]$ homogène en les variables $X_1, \dots, X_n$ tel que $$Q(\alpha, X_1, \dots, X_n)=P(X_1, \dots, X_n) \quad \mathrm{et} \quad Q(z, f_1(z), \dots, f_n(z))=0.$$
\end{Th}

Finalement, le théorème suivant a été prouvé par André dans l'article \cite{Andre2014} dans lequel il développe une généralisation de la correspondance de Galois différentielle. C'est une conséquence de \cite[Corollaire 1.7.1, p. 6]{Andre2014}, qui s'applique non seulement, sous certaines hypothèses, à une famille de fonctions de la forme $(y,y', \dots, y^{(n-1)})$ quand $y$ est solution d'une équation différentielle à coefficients dans $\Qbar(z)$, mais plus généralement à tout vecteur de fonctions $(f_1, \dots, f_n)$ solution d'un système différentiel $y'=Ay$, $A \in \mathcal{M}_n(\Qbar(z))$, comme Y. André nous l'a confirmé. Les hypothèses du corollaire 1.7.1 sont satisfaites si $(f_1, \dots, f_n)$ est un vecteur de $E$-fonctions  \emph{au sens large}, car, d'une part, toute $E$-fonction \emph{au sens large} non polynomiale est transcendante sur $\Qbar(z)$ puisque c'est une fonction entière, et d'autre part, le théorème \ref{th:siegelshidlovskii} nous assure que le degré de transcendance sur $\Qbar$ de $(f_1(\xi), \dots, f_n(\xi))$ est égal au degré de transcendance sur $\Qbar(z)$ de $(f_1(z), \dots, f_n(z))$.

\begin{Th}[André, 2014, \cite{Andre2014}] \label{th:beukersgeneralise}
Le théorème \ref{th:beukersefonctions} reste vrai si l'on remplace partout \og{} strict \fg{} par \og{} large \fg{}.
\end{Th}

Notre but dans cette partie est de fournir une nouvelle preuve du théorème \ref{th:beukersgeneralise} plus proche de l'esprit initial de la preuve de Beukers, à l'aide de l'étude des $E$-opérateurs \emph{au sens large} menée dans la sous-section \ref{subsec:Eoplarges}.

Le point central de la preuve est la proposition suivante, conséquence du théorème \ref{th:andrelarge}. Beukers en a prouvé l'analogue au sens strict dans \cite[p. 372]{Beukers2006}. 

\begin{prop} \label{prop:singulariteapparenteQ}
Soit $f$ une $E$-fonction \emph{au sens large} à coefficients rationnels. Alors si $f(1)=0$, l'équation différentielle minimale de $f$ admet $1$ pour singularité, et elle est apparente.
\end{prop}

\begin{dem}
Écrivons $f(z)=\sum\limits_{n=0}^{\infty} \gf{f_n}{n!} z^n=(1-z)g(z)$, où $g(z)=\sum\limits_{n=0}^{\infty} \gf{g_n}{n!} z^n$. Alors, puisque $f(1)=0$, on a $$\forall n \in \N, \quad \gf{g_n}{n!}=\sum\limits_{k=0}^n \gf{f_k}{k!}=-\sum\limits_{k=n+1}^{\infty} \gf{f_k}{k!}.$$
Soient $\varepsilon >0$ et $n_0(\varepsilon) \in \N$ tels que $\forall n \geqslant n_0(\varepsilon), |f_n| \leqslant (n!)^{\varepsilon}$. On déduit d'une majoration du reste de la série $\sum\limits_{n \in \N} (n!)^{\varepsilon-1}$ que $$\forall n \geqslant n_0(\varepsilon), \quad |g_n| \leqslant n! \left| \sum\limits_{k=n+1}^{\infty} \gf{1}{(k!)^{1-\varepsilon}} \right| \leqslant \gf{n!}{(n!)^{1-\varepsilon}} \leqslant (n!)^{\varepsilon},$$ et de plus si $d_n=\mathrm{den}(f_0, \dots, f_n)$, on a $\forall k \leqslant n, d_n g_k \in \Oal_{\Qbar}$, donc $\mathrm{den}(g_0, \dots, g_n) \leqslant d_n \leqslant (n!)^{\varepsilon}$ pour $n$ suffisamment grand. Ainsi, comme $g$ est solution de $L((1-z)g(z))=0$ si $f$ vérifie $L(f(z))=0$, $L \in \Qbar(z)\left[\mathrm{d}/\mathrm{d}z\right]$, on en conclut que $g$ est une $E$-fonction \emph{au sens large}.

Selon le théorème \ref{th:andrelarge} \textbf{a)}, si $M(y(z))=0$ est une équation différentielle minimale de $g$, $M \neq 0$, alors elle a une base de solutions holomorphes au voisinage de $1$, notée $(g_1(z), \dots, g_{\mu}(z))$. Donc si $L=M \circ (1-z) \in \Qbar(z)\left[\mathrm{d}/\mathrm{d}z\right]$, l'équation différentielle $L(y(z))=0$ est d'ordre minimal $\mu$ pour $f$ et possède au voisinage de $0$ une base de solutions $((1-z) g_1(z), \dots, (1-z) g_{\mu}(z))$ holomorphes s'annulant en $1$. Donc elle a une singularité en $1$, mais elle est apparente. 
\end{dem}

En répétant \emph{mutatis mutandis} la preuve, basée sur des arguments de théorie de Galois différentielle, de \cite[pp. 373--374]{Beukers2006}, on en déduit le théorème suivant. Les points clefs qui permettent de généraliser sa démonstration au cas large sont le fait qu'une $E$-fonction \emph{au sens large} est une fonction entière et la proposition \ref{prop:singulariteapparenteQ}, le reste relève d'arguments généraux sur les groupes algébriques.

\begin{Th} \label{th:singulariteapparentegenerale}
Soit $f(z)$ une $E$-fonction \emph{au sens large} d'équation différentielle minimale $L(y(z))=0$. Soit $\xi \in \Qbar^*$ tel que $f(\xi)=0$. Alors $L$ a une singularité en $\xi$, qui est apparente.
\end{Th}

Le théorème \ref{th:beukersgeneralise} s'ensuit par une preuve identique à celle de \cite[pp. 375--377]{Beukers2006}.

\medskip

\begin{rqu}
Avec les notations du théorème \ref{th:singulariteapparentegenerale}, Gorelov a montré par des moyens différents dans \cite[Lemma 6, p. 520]{Gorelov2004} que toute singularité $\xi \in \Qbar^*$ de $L$ est apparente, ce qui est une partie de ce théorème.
\end{rqu}

\subsection{D'autres résultats diophantiens}

Le théorème suivant a été aussi démontré dans \cite{Beukers2006}. Nous allons en prouver un analogue \emph{au sens large}.

\begin{Th}[Beukers, \cite{Beukers2006}, Theorem 1.5, p. 371] \label{th:beukerscomplement}
Soit $\mathbf{f}={}^t (f_1, \dots, f_n)$ une famille de $E$-fonctions \emph{au sens strict} libre sur $\Qbar(z)$ telle que $\mathbf{f}'=A\mathbf{f}$, où $A  \in \mathcal{M}_n(\Qbar(z))$. 

Alors il existe des $E$-fonctions \emph{au sens strict} $e_1, \dots, e_n$ et une matrice $M \in \mathcal{M}_n(\Qbar[z])$ telles que :
\begin{itemizeth}
 \item On a \, $\mathbf{f}=M {}^t (e_1, \dots, e_n)$ ;
 \item La famille $(e_1, \dots, e_n)$ est solution d'un système de $n$ équations différentielles homogènes du premier ordre à coefficients dans $\Qbar\left[z,z^{-1}\right]$.
\end{itemizeth}

\end{Th}
\bigskip
\begin{Th} \label{th:beukerscomplementlarge}
Le théorème \ref{th:beukerscomplement} reste vrai si l'on remplace partout \og strict \fg{} par \og large \fg{} .
\end{Th}

La preuve est une conséquence du théorème \ref{th:singulariteapparentegenerale}, elle consiste à décrire un algorithme de désingularisation d'un système différentiel du premier ordre, en montrant qu'à chaque étape les fonctions obtenues restent des $E$-fonctions \emph{au sens large}. C'est ce que la proposition suivante permet de prouver. 

\begin{prop}\label{prop:Efonction/z-xi}
Soient $f(z) \in \Qbar\llbracket z\rrbracket$ une $E$-fonction \emph{au sens large} et $\xi \in \Qbar^*$ tel que $f(\xi)=0$. Alors $\gf{f(z)}{z-\xi}$ est une $E$-fonction \emph{au sens large}.
\end{prop}

On peut trouver une version plus faible de cette proposition dans \cite[Lemma 9, p. 522]{Gorelov2004}, qui donne le même résultat sous l'hypothèse supplémentaire que $f(z)$, $f'(z)$, $\dots$, $f^{(m-1)}(z)$ sont algébriquement indépendantes sur $\Qbar(z)$, $m$ étant l'ordre de l'équation différentielle minimale de $f(z)$.

\begin{dem}[de la proposition \ref{prop:Efonction/z-xi}]
Quitte à remplacer $f(z)$ par $f(\xi z)$, qui est une $E$-fonction \emph{au sens large}, on peut supposer que $\xi=1$.  En effet, si $h(z)=f(\xi z)$, on a $$\gf{f(z)}{z-\xi}=\gf{1}{\xi} \gf{h(\xi^{-1}z)}{\xi^{-1} z-1}.$$

Selon le théorème \ref{th:singulariteapparentegenerale}, comme $f(1)=0$, toutes les solutions de l'équation différentielle minimale $L(y)=0$ de $f$ sont holomorphes et s'annulent en $1$. 

Soit $\sigma \in \Gal(\Qbar/\Q)$. En définissant $L^{\sigma}$ comme dans la preuve de la proposition \ref{prop:baseGfonctionspointordinaire}, on voit que $L^{\sigma}(y^{\sigma})=L(y)^{\sigma}$ pour tout $y \in \Qbar\llbracket z\rrbracket$, donc l'espace des solutions de $L^{\sigma}(y)=0$ est l'image par $\sigma$ des solutions de $L(y)=0$. En particulier, $L^{\sigma}$ a une base de solutions holomorphes s'annulant en $1=\sigma(1)$. D'où $f^{\sigma}(1)=0$. En répétant l'argument utilisé dans la preuve de la proposition \ref{prop:singulariteapparenteQ}, on obtient que si $$g(z):=\gf{f(z)}{1-z}=\sum\limits_{n=0}^{\infty} \gf{g_n}{n!} z^n,$$ alors pour tout $\varepsilon >0$, $|\sigma(g_n)| \leqslant (n!)^{\varepsilon}$ pour $n \geqslant n_0(\sigma, \varepsilon)$. Comme $g(z)$ est à coefficients dans un corps de nombres, il n'y a qu'un nombre fini de $\sigma$ à considérer, si bien que $\house{g_n} \leqslant (n!)^{\varepsilon}$ pour $n \geqslant n_1(\varepsilon)$.

La condition sur les dénominateurs étant vérifiée de la même manière que dans la preuve de la proposition \ref{prop:singulariteapparenteQ}, on en déduit que $g$ est une $E$-fonction \emph{au sens large}.
\end{dem}

Une fois cette proposition prouvée, on obtient le théorème \ref{th:beukerscomplementlarge} en répétant \emph{mutatis mutandis} la démonstration de \cite[p. 378]{Beukers2006}, à l'aide également du théorème \ref{th:andrelarge} \textbf{a)} et du théorème~\ref{th:singulariteapparentegenerale}. 

\medskip

Dans \cite{AdamczewskiRivoal}, Adamczewski et Rivoal ont construit, à partir de l'algorithme de la preuve du théorème \ref{th:beukerscomplement}, un algorithme permettant de déterminer les points auxquels une $E$-fonction \emph{au sens strict} donnée prend des valeurs algébriques.

\begin{Th}[Adamczewski, Rivoal]
Il existe un algorithme qui, étant donnée une $E$-fonction \emph{au sens strict} $f(z)$, indique si $f(z)$ est transcendante sur $\Qbar(z)$ ou non et, si elle est transcendante, fournit la liste finie des nombres algébriques $\alpha$ tels que $f(\alpha)$ est algébrique, ainsi que la liste des valeurs $f(\alpha)$ correspondantes.
\end{Th}

\begin{prop}
L'algorithme du théorème précédent fonctionne si $f(z)$ est une $E$-fonction \emph{au sens large}.
\end{prop}

En effet, l'algorithme en question est décomposé en trois étapes, les deux premières n'étant pas particulières aux $E$-fonctions.
\begin{itemize}[label=\textbullet]
    \item D'abord, il détermine une équation différentielle homogène minimale $M(y(z))=0$ sur $\Qbar(z)$ satisfaite par $f(z)$. La preuve de la terminaison de cette étape vient essentiellement d'une borne sur le degré en $z$ de $M$ issue d'un théorème de Grigoriev \cite[p. 8]{Grigoriev}, rendu pleinement explicite par Bostan, Rivoal et Salvy \cite[Theorem 1, p. 2]{BostanRivoalSalvy}, d'une borne sur les modules des exposants de $M$ dûe à Bertrand, Chirskii et Yebbou \cite[p. 246, p. 252]{BertrandChirskiiYebbou} et enfin d'un lemme de zéros de Bertrand et Beukers \cite[p. 182]{BertrandBeukers}. Ces trois résultats sont valables sur les opérateurs différentiels linéaires à coefficients dans $\Qbar(z)$ en général.
    \item Puis il détermine, à partir de l'équation de l'étape précédente, une équation différentielle minimale inhomogène sur $\Qbar[z]$, $u_0(z) y^{(s)}(z)+ \dots+u_s(z) y(z)+u_{s+1}(z)=0$, avec $u_0 \not\equiv 0$, dont $f(z)$ est solution. Là encore, cette procédure est valable pour n'importe quel opérateur différentiel.
    \item La troisième étape détermine quels éléments $\alpha$ parmi les racines du polynôme $u_0(z)$ vérifient $f(\alpha) \in \Qbar$. Le fait que les autres nombres algébriques sont exclus est une conséquence du théorème \ref{th:beukersefonctions}, qui reste vrai au sens large selon le théorème \ref{th:beukersgeneralise}. Cette étape consiste principalement à appliquer l'algorithme de désingularisation de Beukers décrit et utilisé dans la preuve du théorème \ref{th:beukerscomplement}, qui est valable au sens large selon le théorème \ref{th:beukerscomplementlarge}. 
\end{itemize}

\printbibliography

\bigskip
G. Lepetit, Université Grenoble Alpes, CNRS, Institut Fourier, 38000 Grenoble, France.

\url{gabriel.lepetit@live.fr}. 

\bigskip

\emph{Keywords} : $E$- and $G$-functions, $E$- and $G$-operators, Chudnovsky's Theorem.

\bigskip

\emph{2020 Mathematics Subject Classification}. Primary 11J91 ; Secondary 34M03, 34M35

\end{document}